\documentclass{amsart}
\input{epsf}
\usepackage{xypic,verbatim,amsmath,latexsym,amssymb,amsbsy}

{\end{list}}
{
   \newtheorem{theorem}{Theorem}[subsection]                     
   \newtheorem{proposition}[theorem]{Proposition}     
   \newtheorem{lemma}[theorem]{Lemma}

   \newtheorem{corollary}[theorem]{Corollary}
   \newtheorem{conjecture}[theorem]{Conjecture}

}
{\theoremstyle{definition}
   
   \newtheorem{example}[theorem]{Example}
   \newtheorem{definition}[theorem]{Definition}
}
{\theoremstyle{remark}
   \newtheorem*{remark}{Remark}
   
}
\newcommand{\QQ}{{\mathbb{Q}}}
\newcommand{\NN}{{\mathbb{N}}}
\newcommand{\PP}{{\mathbb{P}}}
\newcommand{\ZZ}{{\mathbb{Z}}}

\newcommand{\bbA}{{\mathbb{A}}}

\newcommand{\cC}{{\mathcal C}}
\newcommand{\cI}{{\mathcal I}}
\newcommand{\cJ}{{\mathcal J}}

\newcommand{\cO}{{\mathcal O}}

\newcommand{\cZ}{{\mathcal Z}}

\newcommand{\prim}{{\operatorname{prim}}}

\newcommand{\reg}{{\operatorname{reg}}}
\newcommand{\grad}{{\operatorname{grad}}}

\newcommand{\val}{{\operatorname{val}}}

\newcommand{\cl}{{\operatorname{cl}}}

\newcommand{\Ver}{{\operatorname{Vert}}}

\newcommand{\Tan}{{\operatorname{Tan}}}
\newcommand{\inte}{{\operatorname{int}}}
\newcommand{\lin}{\operatorname{lin}}

\newcommand{\Gl}{\operatorname{Gl}}
\newcommand{\Spec}{\operatorname{Spec}}
\newcommand{\Proj}{\operatorname{Proj}}

\newcommand{\Star}{\operatorname{Star}}

\newcommand{\Mid}{{\operatorname{Mid}}}
\newcommand{\Ctr}{{\operatorname{Ctr}}}

\newcommand{\inv}{{\operatorname{inv}}}

\newcommand{\id}{{\operatorname{id}}}

\newcommand{\gr}{{\operatorname{gr}}}
\newcommand{\das}{\dashrightarrow}

\newcommand{\dra}{\dashrightarrow}
\newcommand{\Gal}{{\operatorname{Gal}}}

\newcommand{\para}{{\operatorname{par}}}

\newcommand{\cW}{\cite{Wlodarczyk2}}
\newcommand{\Otau}{\widetilde{O_\tau}}
\newcommand{\ks}{K^*}
\newcommand{\Ga}{\Gamma}

\newcommand{\bZ}{\mathbb Z}
\newcommand{\bQ}{\mathbb Q}
\newcommand{\bP}{\mathbb P}
\newcommand{\bN}{\mathbb N}
\newcommand{\bR}{\mathbb R}

\newcommand{\sig}{\sigma}
\newcommand{\la}{\langle}
\newcommand{\ra}{\rangle}
\newcommand{\Sig}{\Sigma}
\newcommand{\Del}{\Delta} 
\newcommand{\del}{\delta} 
\newcommand{\wtx}{\widetilde X}

\newcommand{\ox}{\overline x}

\newcommand{\whx}{\widehat X}    
\newcommand{\vr}{\varrho}

\newcommand{\ds}{\displaystyle}
\newcommand{\sgn}{{\operatorname{sgn}}}
\newcommand{\varp}{\varphi}

\begin{document}


\title[Simple Constructive Weak Factorization]{Simple Constructive  Weak Factorization}
\author{Jaros{\l}aw W{\l}odarczyk}
\thanks{The author was supported in part by  NSF grant DMS-0500659 and Polish KBN
grant GR-1784 .}
\address{Department of Mathematics\\Purdue University\\West
Lafayette, IN 47907\\USA}

\email{wlodar@math.purdue.edu}
\date{\today}
\begin{abstract}
We give a simplified algorithm of the functorial weak factorization of birational morphisms of nonsingular varieties over a field  of characteristic
zero into a composite of blow-ups and blow-downs with smooth centers.  
\end{abstract}
\maketitle

\tableofcontents
\addtocounter{section}{-1}

\section{Introduction}

In this paper we give a simplified version of our proof of the following theorem:
\begin{theorem} {\bf The Weak Factorization Theorem}\label{th:1}
 \begin{enumerate}
\item  Lef $f:X\dashrightarrow Y$ be a 
birational map  of
smooth complete varieties 
over a field of characteristic zero, which is an
isomorphism over an open set $U$. Then $f$ can be factored
 as
$$X=X_0\buildrel f_0 \over \dashrightarrow  X_1
\buildrel f_1 \over \dashrightarrow \ldots \buildrel f_{n-1} \over
\dashrightarrow X_n=Y  ,$$
where each $X_i$ is a smooth complete variety and $f_i$ is a blow-up
or blow-down at a smooth center which is an isomorphism
over $U$.
\item Moreover, if $X\setminus U$ and $Y\setminus U$ are divisors
with simple normal crossings, then each $D_i:=X_i\setminus U$ 
is a divisor with simple normal crossings and $f_i$ is a blow-up
or blow-down at a smooth center which has normal
crossings with components of $D_i$.

\item There is an index $1\leq r\leq n$  such that for  all $i\leq r$ the induced  birational map $X_i\buildrel f_0 \over\dashrightarrow X$ is a projective morphism and   for all $r\leq i \leq n$ the map  $X_i \buildrel f_0 \over\dashrightarrow Y$ is projective morphism.

\item The above factorization is functorial in the following sense:

\noindent Let $\phi_X$, $\phi_Y$  and $\phi_K$ be automorphisms of $X$ , $Y$ and $\Spec(K)$ such that $f\circ \phi_X=\phi_Y\circ f$, $j_X\circ\phi_X=j_Y\circ\phi_Y=\phi_K$ ,where $j_X: X\to\Spec(K)$ and $j_Y:Y\to \Spec(K)$  are  the natural morphisms. Then
the induced birational transformations   $\phi_i: X_i\dashrightarrow X_i$  are automorphisms of $X_i$ commuting with $f_i:X_i\to X_{i+1}$ and $j_{X_i}:X_i\to \Spec(K)$.
  Moreover if $\phi_X(D_X)=D_X$ and $\phi_Y(D_Y)=D_Y$ then for all $i$, we have  $\phi_i(D_i)=D_i$.

\item The factorization commute with field extensions $K\subset L$.
\end{enumerate}
\end{theorem}
The theorem was proven in \cite{Wlodarczyk3} and in \cite{AKMW} in a more general version. The above formulation essentially reflects the statement of the Theorem in \cite{AKMW}.

The weak factorization theorem extends a theorem of Zariski, which states that
any birational map between two smooth complete
surfaces can be factored into a succession of blow-ups at
points followed by a succession of blow-downs at points.  
A stronger version of the above theorem, called the strong factorization 
conjecture, remains open. 

\begin{conjecture} {\bf Strong Factorization Conjecture}.
Any birational map $f:X\dashrightarrow Y$ of  smooth
complete varieties can be factored into a succession of blow-ups
at
smooth centers followed by a succession of 
blow-downs at
smooth centers.
\end{conjecture}

Note that both statements are equivalent in dimension 2.
One can find the formulation of the relevant  conjectures in many papers.  Hironaka \cite{Hironaka1} formulated the strong
factorization conjecture. The weak factorization
problem was stated by Miyake and Oda \cite{Miyake-Oda}. The toric versions of the
strong and weak factorizations were also conjectured by
Miyake and Oda \cite{Miyake-Oda} and are called the strong and weak Oda
conjectures. The $3$-dimensional toric version of the weak form was
established by Danilov \cite{Danilov2} (see also Ewald \cite{Ewald}). 
The weak  toric conjecture 
in arbitrary dimensions was proved in   \cite{Wlodarczyk1} and later
independently by Morelli \cite{Morelli1},  who also claimed to have a proof of the
strong factorization conjecture 
(see also Morelli \cite{Morelli2}).
  Morelli's
proof of the weak Oda conjecture was completed, revised and
generalized to the toroidal case by Abramovich, Matsuki and Rashid
in \cite{Abramovich-Matsuki-Rashid}.  A gap in Morelli's proof of the strong Oda conjecture,
 which went
unnnoticed in \cite{Abramovich-Matsuki-Rashid}, was later
found by K. Karu.

The local version of the strong factorization problem was posed
by Abhyankar in dimension 2 and by Christensen  in general; Christensen has solved it for 3-dimensional toric varieties \cite{Christensen}.
The local version of the weak factorization problem (in
characteristic 0) was 
solved by Cutkosky  \cite{Cutkosky1}, who also 
showed that Oda's strong conjecture implies the local
version of the strong conjecture for proper  birational
morphisms \cite{Cutkosky2} and proved the local strong factorization conjecture
in dimension 3 \cite{Cutkosky2} via Christensen's theorem.
Finally Karu generalized Christensen's result to any dimension and completed the argument for the local strong factorization \cite{Karu}.

The proofs in  \cite{Wlodarczyk3} and \cite{AKMW} are both build upon the idea of cobordisms which was developed in
\cW \,   and was inspired by  Morelli's theory of polyhedral cobordisms \cite{Morelli1}.   The main idea of \cW\,  is to construct a space with a $K^*$-action for a given birational map. 
The space called a birational cobordism resembles  the idea of Morse cobordism and determines  a decomposition of the birational map into elementary transformations (see Remark \ref{cobordism}).
This gives a factorization into a sequence of 
weighted blow-ups and blow-downs. One can view   the factorization determined by the cobordisms also in terms of VGIT  developed in papers of Thaddeuss and Dolgachev-Hu.
As  shown in \cW\,  the weighted blow-ups which occur in the factorization have a
natural local toric description which is crucial for their further regularization.

The two existing methods of regularizing centers of this factorization are $\pi$-desingularization of cobordisms as in \cite{Wlodarczyk3} and local torification of the action  as in   \cite{AKMW}.

The present proof is essentially the same as in \cite{Wlodarczyk3}.
Instead of working in full generality and developing the suitable language for toroidal varieties we  focus on 
applying the general ideas to a particular construction of a smooth cobordism.
The $\pi$-desingularization is a desingularization of geometric quotients
of a $K^*$-action. This can be done locally and the procedure  can be globalized in the functorial and even canonical way. The $\pi$-desingularization makes all the intermediate varieties (which are geometric quotients) smooth, and also the connecting blow-ups have smooth centers. 
 
 The proof of Abramovich, Karu, Matsuki and the author \cite{AKMW}
 relies on a subtle analysis of differences between locally toric and toroidal structures defined by the action of $K^*$. The Abramovich-de Jong idea of torification is roughly speaking to construct the ideal sheaves whose blow-ups (or principalizations) introduce the structure of toroidal varieties in neighborhoods of fixed points of the action. This allows one to pass from birational maps between intermediate varieties in the neighborhood of fixed points  to toroidal maps. The latter can be factored into a sequence of smooth blow-ups by using the same combinatorial methods as for toric varieties. Combining all the  local factorizations together we get a global factorization. 

In the presentation of birational cobordisms below we base on \cW\, with some improvements in \cite{AKMW}. In particular we use Hironaka flattening for factorization into projective morphisms, and elements of GIT to show existence of quotients.
The presentation of the paper is  self contained. In particular the toric version of the weak factorization is proven in Section \ref{toric} to illustarate some of the ideas of the proof.


\section{Birational cobordisms}
\subsection{Definition of a birational cobordism}

Recall some basic definitions from Mumford's GIT theory.

\noindent  \begin{definition} Let $K^*$ act on $X$. By a {\it good
quotient}  we mean a variety $Y=X//K^*$ together with
a morphism $\pi:X\rightarrow Y$ which is constant on
$G$-orbits such that for any affine open subset
$U\subset Y$ the inverse image $\pi^{-1}(U)$ is affine and
$\pi^*:O_Y(U)\rightarrow O_X(\pi^{-1}(U))^{K^*}$ is an isomorphism.
If additionally for any closed point $y\in Y$ its inverse limit $\pi^{-1}(x)$
is a single orbit we call $Y:=X/K^*$ together with
$\pi:X\rightarrow Y$ a {\it geometric quotient}.
 \end{definition}

\begin{remark} A geometric quotient is a space of orbits while a good quotient is a space of equivalence classes of orbits generated by the relation that two 
orbits
are equivalent if their closures intersect.
\end{remark} 
\begin{definition}
Let $K^*$ act on $X$. We say that $ \lim_{{t}\to 0} \, {t}x $ exists (respectively $\lim_{{t}\to \infty} \, { t}x$  exists ) if the morphism
$\Spec(K^*)\to X$ given by $t\mapsto {tx}$ extends to $\Spec(K^*)\subset\bbA^1\to X$ (or respectively $\Spec(K^*)\subset\bP^1\setminus\{0\}\to X$).
\end{definition}
 \begin{definition} (\cW) 
Let $X_1 $ and $X_2$ be two birationally equivalent
 normal varieties.
A {\it birational cobordism} or simply a {\it cobordism} 
$B:=B(X_1,X_2)$ between them is a
normal variety $B$ with an algebraic action of $K^*$ such that the
sets
\[\begin{array}{rccc}
&B_-:=\{x \in B\mid  \, \lim_{{ t}\to 0} \, { t}x \,\, \mbox{does
not exist}\}& \, \mbox{and}& \\
&B_+:=\{x \in B\mid \, \lim_{{ t}\to\infty} \, { t}x \,\, 
\mbox{does not exist}\} && 
\end{array}\]
are nonempty and open and there exist 
 geometric quotients
$B_-/ K^*$ and $B_+\// K^*$ such that $B_+/ K^*\simeq X_1$ and $B_-\//
K^*\simeq X_2$ and the birational map $X_1 \dashrightarrow
X_2$ is given by the above
isomorphisms and the open embeddings of
$B_+ \cap B_-\// K^*$ into
$B_+\//K^*$ and $B_-\//K^*$ respectively.
\end{definition}

\begin{remark} An analogous notion of cobordism of fans of
toric varieties was introduced by Morelli in \cite{Morelli1}.
\end{remark}
\begin{remark}\label{cobordism}
The above definition can  also be considered as an analog of the notion of
cobordism in Morse theory. Let $W_0$ be a cobordism in Morse
theory of two
differentiable manifolds $X$ and $X'$ and $f:W_0\rightarrow
[a,b]\subset \bR$ be a
Morse function such that $f^{-1}(a)=X$ and $f^{-1}(b)=X'$.
Then $X$ and $X'$ have open neighborhoods $X\subseteq
V\subseteq W_0$ and $X'\subseteq
V'\subseteq W$ such that  $V\simeq X\times
[a,a+\epsilon )$  and $V'\simeq X'\times (b-\epsilon ,b]$ for which
$f_{\mid V}:V\simeq X \times [a,a+\epsilon)\rightarrow [a,b]$ and
$f_{\mid V'}:V'\simeq X' \times (b-\epsilon ,b] \rightarrow [a,b]$ are
the natural projections on the second coordinate.  Let
$W:=W_0\cup_V X\times
(-\infty,a+\varepsilon) \cup_{V'}
X'\times (b-\varepsilon,+\infty )$. One can easily see
that $W$ is isomorphic to $W_0\setminus X\setminus
X'=\{x\in W_0\mid  a < f(x) < b\}$. Let $f':W\rightarrow \bR$ be the
map defined by glueing the function $f$ and the natural projection on the second
coordinate. Then $\grad (f')$ defines an action on $W$ of
a $1$-parameter group
 $T\simeq {\bR}\simeq {\bR}_{>0}^*$ of diffeomorphisms. The last group isomorphism is given by the exponential.

 Then one can see that $W_-:=\{x \in W\mid  \lim_{t\to 0}
\,  tx  \  \mbox{does not exist}\}$ and $W_+:=\{x
\in W\mid \lim_{t\rightarrow\infty} \, tx \  \mbox{does
not exist}\}$ are open and  $X$
 and $X'$ can be considered as quotients of these sets by $T$.
 The critical points of the Morse function are $T$-fixed points. ``Passing through the fixed points'' of the action
 induces a simple birational transformation similar to spherical modification in Morse theory (see Example \ref{main}).
\end{remark}
\begin{figure}[ht]

\epsfysize=1.5in

	\epsffile{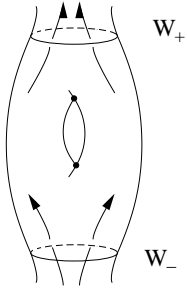}

	\caption{Cobordism in Morse theory}\label{Fi:1}
\end{figure}

\begin{example} 
\label{main}

Let $K^*$ act on $B:=\bbA^{l+m+r}_K$ by
$$ t(x_1,\ldots,x_l,y_1,\ldots,y_m,z_1,\ldots,z_r)=(t^{a_1}\cdot x_1,\ldots,t^{a_l} \cdot
x_l,t^{-b_1}\cdot y_1,\ldots,t^{-b_m}\cdot y_m, z_1,\ldots,z_r),$$
where $a_1,\ldots,a_l,b_1,\ldots,b_m>0$.
Set $\overline{x}=(x_1,\ldots,x_l) , \overline{y}=(y_1,\ldots,y_m) ,
\, \overline{z}=( z_1,\ldots,z_r)$.
Then
$$\displaylines{ B_-= \{p=(\overline{x},\overline{y},\overline{z})\in
\bbA^{l+m+r}_K \mid  \, \overline{y}\neq 0\},\cr
B_+= \{p=(\overline{x}, \overline{y},\overline{z})\in
\bbA^{l+m+r}_K \mid  \, \overline{x}\neq 0 \}.\cr}$$
%
%
\noindent{\bf Case 1.}  $a_i=b_i=1$, $r=0$ (Atiyah, Reid).

 One can easily see that $B//K^*$ is the affine cone over
the Segre embedding ${\PP}^{l-1}\times {\PP}^{m-1}\rightarrow
{\PP}^{l\cdot m -1}$, and $B_+/K^*$ and $B_-/K^*$ are smooth.

The relevant birational map $\phi : B_-/K^* \dashrightarrow B_+/K^*$ is a
flip
for $l, m\geq 2$ replacing ${\PP}^{l-1}\subset B_-/K^* $
with ${\PP}^{m-1}\subset B_+/K^* $.
For $l=1, m\geq2$, $\phi$ is a blow-down, and for
$l\geq2, m=1$ it is a blow-up. If $l=m=1$ then  $\phi$ is the identity.
One can show that $\phi : B_-/K^* \dashrightarrow B_+/K^*$ factors into the blow-up of
  ${\PP}^{l-1}\subset  B_-/K^* $  followed by the blow-down of ${\PP}^{m-1}\subset B_+/K^* $.
  %
%
%
%

\noindent {\bf  Case 2}. General case. 

For $l=1$, $m\geq2$, $\phi$ is a toric blow-up whose
exceptional fibers are weighted projective spaces. For
$l\geq2$, $m=1$, $\phi$ is a toric blow-down. If $l=m=1$ then
$\phi$ is the identity. The birational map $\phi : B_-/K^* \dashrightarrow B_+/K^*$ factors into  
a weighted blow-up and a weighted blow-down.

\noindent {\bf  Case 3}. $l=0$ and $m\neq 0$ (or $l\neq 0$ and $m=0$).

In this case we have only negative and zero weights (respectively positive and zero weights.)
Then $B=\bbA^{l+m+r}$ is not a cobordism.
In particular $B_+=\emptyset$ .The morphism  $B_-/K^*=\bP(\bbA^m)\times \bbA^r\to B//K^*=\bbA^r$ is the standard projection, where $\bP(\bbA^m)$ is the weighted projective space defined by the action of $K^*$ on $\bbA^m$.

\end{example}
\begin{figure}[ht]

\epsfysize=.8in

	\epsffile{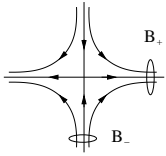}

	\caption{Affine Cobordism}\label{Fi:2 }
\end{figure}

\begin{figure}[ht]

\epsfysize=.8in

	\epsffile{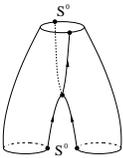}

	\caption{Spherical modifications}\label{Fi:3}
\end{figure}

\begin{remark} In Morse theory we have an
analogous situation. In  cobordisms with one critical point we replace
$S^{l-1}$ by $S^{m-1}$. (See Figure \ref{Fi:3})
\end{remark}

\subsection{Fixed points of the action} 
Let $X$ be a variety with an action of $K^*$. Denote by $X^{K^*}$ the set of fixed points of the action and by $\cC(X^{K^*})$ the set of its irreducible fixed components. For any  $F\in \cC(X^{K^*})$ set $$F^+(X)=F^+=\{x\in X\mid 
\, \lim_{{t}\to 0} {t}x \in F\},  \quad F^-(X)=F^-=\{x\in X\mid 
\, \lim_{{t}\to \infty} {t}x \in F\}.$$




\begin{example} In Example \ref{main}, $$F=\{p \in B\mid \overline{x}=\overline{y}= 0\}, \quad\quad F^-=\{p\in B\mid \overline{x}= 0\}, \quad\quad F^+=\{p\in B\mid \overline{y}= 0\}.$$
\end{example}
\begin{lemma}\label{fix} If $F$ is  the fixed point set  of an affine variety $U$ then $F$, $F^+$ and $F^-$ are closed in $U$. Moreover the ideals $I_{F^+}, I_{F^-}\subset K[V]$ are  generated by all semiinvariant functions with positive (respectively negative) weights.
\end{lemma}
\begin{proof} Embed $U$ equivariantly into affine space $\bbA^n$ with linear action and use the example above.
\end{proof} 
\subsection{Existence of a smooth birational cobordism}\label{se: construction}
The following result is a consequence of Hironaka flatenning theorem \cite{Hironaka4}.
\begin{proposition}{fact}
Let $\phi:X\dashrightarrow Y$ be a birational map between smooth projective
varieties. Then $\phi$ factors as $X\leftarrow Z\to Y$, where $Z\to X$ and $Z\to Y$ are birational morphisms  from a smooth projective variety $Z$. 
The above factorization is functorial. Moreover  there exist functorial divisors  $D_X$ and $D_Y$ on $Z$ which are relatively ample over $X$ and $Y$ respectively. 
If $\phi$ is an isomorphism over $U$ and  the complements $X\setminus U$ and $Y\setminus U$ are simple normal crossing divisors then $Z\setminus U$ is a simple normal crossing divisor.

\end{proposition}
\begin{proof} Let $\Ga(X,Y)\subset X\times Y$ be the graph of $\phi$ and $Z_0$ be its canonical resolution of singularities \cite{Hironaka3}. 
 If  $X$ and $Y$ are projective we take  simply $Z=Z_0$. If $X$ and $Y$ are arbitrary we can apply Hironaka flattening to $Z_0\to Y$ to find a projective factorization 
$\phi: Z_0\leftarrow Z_Y\to Y$, where $Z_Y\to Y$ is a composition of blow-ups at smooth centers and $Z_Y\to Z_0$ is a composition of blow-ups which are pull-backs of these blow-ups (\cite{Hironaka4},\cite{Raynaud-Gruson}). Next we apply  Hironaka flattening
to $Z_Y\to X$ to obtain a factorization $Z_Y\leftarrow Z_X\to X$.
Finally, $Z\to Z_X$ is a canonical prinicipalization of $\cI_D$, where $D$ is the complement of $U$ on $Z_X$. The divisors $D_X$ are $D_Y$ are  constructed as a combination of components of the exceptional divisors of $Z\to X$ and $Z\to Y$ respectively.
\end{proof}


It suffices to construct the cobordism and factorization for the projective morphism $Z\to X$.

\begin{proposition} (\cW,\cite{AKMW}) \label{construction} Let $\varphi: Z\to X$ be a birational projective morphism of smooth complete varieties with the exceptional divisor $D$. Let $U\subset X, Z$ be
an open subset where $\varphi$ is an isomorphism. 
There exists a smooth complete variety $\overline{B}$ with a $K^*$-action, which contains fixed point components isomorphic to $X$ and $Z$ such that
\begin{enumerate}
\item
\begin{itemize}
\item $B=B(X,Z):= \bar{B} \setminus(X \cup Z)$ is a cobordism  between $X$ and $Z$.
\item $U\times K^*\subset B_-\cap B_+\subset B$.
\item There are $K^*$-equivariant isomorphisms $X^-\simeq X\times ({\bP^1\setminus\{0\}})$ and $Z^+\simeq\cO_Z(D)$. 
\item $X^-\setminus X=B_+$ and $Z^+\setminus Z=B_-$


\item There exists a $K^*$-equivariant projective morphism  $\pi_B:\overline{B}\to X$  such that $i_X\pi_B=\id_X$ and $i_Z\pi_B=f$, where $i_X: X\hookrightarrow \overline{B}$ and $i_Z: Z\hookrightarrow \overline{B}$ are embeddings of $X$ and $Z$.
Here the action of $K^*$ on $X$ is trivial.
\item There is  a relatively ample divisor for $\pi_B$ which is functorial and in particular $K^*$-invariant.
 
\end{itemize}

\item If $D_X:=X\setminus U$ and $D_Z:=Z\setminus U$ are divisors
with simple normal crossings then there exists 
a smooth cobordism $\tilde{B}\subset\overline{\tilde{B}}$ between $\tilde{X}$ and
$\tilde{Z}$ as in (1) such that 
\begin{itemize}
\item $\tilde{X}$ and $\tilde{Z}$ are obtained from $X$ and $Z$ by a sequence of
blow-ups at centers which have normal crossings with
components of the total transforms of $D_X$ and $D_Z$ respectively.
\item $U\times \bP^1 \subset \tilde{B} $ and
$\overline{\tilde{B}}\setminus (U\times \bP^1)$ is a
divisor with simple normal crossings.
\end{itemize}
\end{enumerate}

\end{proposition}
In further considerations we shall refer to $\bar{B}$ as a {\it compactified cobordism.}
\begin{figure}[ht]

\epsfysize=1in

	\epsffile{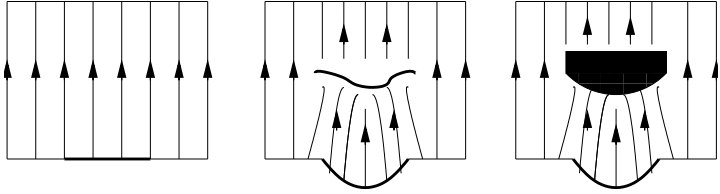}

	\caption{Compactified cobordism}\label{Fi:4}
\end{figure}

\begin{proof} (1) We follow here the Abramovich construction of cobordism.
Let $\cI\subset \cO_X$ be a sheaf of ideals such that
 $Z=Bl_\cI{X}$  is obtained from $X$ by blowing up of $\cI$. Let $z$ denote the standard coordinate on ${\PP}^1$ and
let $\cI_0$ be the ideal of the point $z=0$ on ${\PP}^1$.
 Set
$W:=X\times {\PP}^1$ and denote by $\pi_1:W\rightarrow X$,
$\pi_2:W\rightarrow {\PP}^1$  the standard projections.
 Then
$\cJ:={\pi_1}^*(\cI)+{\pi_2}^*(\cI_0)$ is an ideal supported on $X\times \{0\}$.
Set $W':=Bl_{\cJ}W$.
The proper transform
of $X\times
\{0\}$ is isomorphic to $Z$
and we identify it with $Z$. Let us describe $Z$ locally.
Let $f_1,\ldots,f_k$ generate the ideal $\cI$ on some open affine set $U\subset X$. Then after the blow-up $Z\to X$ at $\cI$ the inverse image of $U$ is a union of open charts $U_i\subset Z$, where $$K[U_i]=K[U][f_i,f_1/f_i,\ldots, f_k/f_i].$$
 Now the functions  $f_1,\ldots,f_k,z$ generate the ideal $\cJ$ on
$U\times \bbA^1\subset W$. After the blow-up $W'\to W$ at $\cJ$, the
inverse  image of $U\times \bbA^1$ is a union of open charts $V_i\supset Z$, where $$K[V_i]=K[U][f_i, f_1/f_i,\ldots, f_k/f_i,z/f_i]=K[U_i][z/f_i]$$  and the relevant $V_z$ which does not intersect $Z$. Then $V_i=U_i^+\simeq U_i\times \bbA^1$ where $ z':=z/f_i$ is the standard coordinate  on $\bbA^1$. The  action of $K^*$
 on the factor $U$ is trivial while on $\bbA^1$ it is standard given by $t(z')=tz'$.
Thus the open subset $Z^+=\bigcup U_i^+=\bigcup V_i\subset W'$ is a line bundle over $Z$ with the standard action of $K^*$. On the other hand the neighborhood $X^-:=X\times (\bP^1\setminus\{0\})$ of $X\subset W$ remains unchanged after the blow-up of $\cJ$. We identify $X$ with $X\times\{\infty\}$.
We define $\overline{B}$ to be the canonical desingularization of $W$.
Then $B:=\overline{B}\setminus X\setminus Z$. We get $B_-/K^*=
(Z^+\setminus Z)/K^*=Z$, while $B_+/K^*=(X^+\setminus X)/K^*=X
$. 
The relatively very ample divisor is the relevant combination of the divisor $X$ and the exceptional divisors with negative coefficients.

(2) The sets $Z^+$ and $X^-$ are line bundles with projections $\pi_+:
Z^+ \to Z$ and $\pi_-: X^-\to X$. Let $Z:=\overline{\tilde{B}}\setminus (U\times \bP^1)$. Then $Z\cap Z^+$ and
$Z\cap X^-$ are simple normal crossing divisors and $\pi_+(Z\cap
Z^+)=D_Z$.
 and $\pi_-(Z\cap
X^-)=D_X$. Let $f:\overline{\tilde{B}}\to \bar{B}$ be a canonical principalization of
$\cI_Z$ (see Hironaka \cite{Hironaka3},
Villamayor \cite{Villamayor} and Bierstone-Milman \cite{Bierstone-Milman}). Let $f_+:f^{-1}(Z^+)\to Z^+$ be the restrion of $f$. By functoriality  
$f_+$  is  a canonical
principalization of $\cI_{Z|Z^+}=\pi_+^*(\cI_{D_Z})$ on $Z^+$ (resp. $X^-$) which commutes with $\pi_+$.
Then $f_+$ is a pull-back of  the canonical principalization $\tilde{Z}\to Z$ of $\cI_{D_Z}$ on $Z$. 
In particular  $f^{-1}(Z^+)=\tilde{Z}^+$ and all centers of blow-ups are $K^*$-invariant and  of the form $\pi_+^{-1}(C)$, where $C$ has
normal crossings with components of the total transform of $D_Z$. Analogously for $X^-$.
\end{proof}

\begin{remark} The Abramovich construction can be considered as a generalization of the Fulton-Macpherson example of the deformation to the normal cone. If we let $\cI=\cI_C$ be the ideal sheaf of the smooth center then the relevant blow-up is already smooth.
On the other hand this a particular case of the very first construction of a cobordism in (\cW\,Proposition 2. p 438)
 which is a $K^*$-equivariant completion of the space
$$L(Z,D;X,0):=O_Z(D)\cup_{U\times K^*} X\times (\bP^1\setminus \{0\}).$$ 
Another  variant of our construction  is given by
Hu and Keel in \cite{Hu-Keel}. 
\end{remark}

\subsection{Collapsibility}
In the following $\bar B\subset B$  denotes compactified cobordism between $X$ and $Z$ subject to the conditions from Proposition \ref{construction} but not necessarily smooth.
\begin{definition} (\cW).\label{de: order} Let $X$ be a
cobordism or any variety with a $K^*$-action.
\begin{enumerate}
\item We say that $F\in\cC(X^{K^*})$
 is an {\it  immediate predecessor} of
$F'\in \cC(X^{K^*})$ if there exists a nonfixed point $x$ such that
$\lim_{{\bf t}\to 0} {\bf t}x\in F$ and
$\lim_{{\bf t}\to
\infty} {\bf t}x\in F'$. 
\item
We say that $F$  {\it precedes}  $F'$ and
write $F<F'$ if there exists a sequence of connected fixed
point set components $F_0=F ,F_1,\ldots,F_l=F'$ such that
$F_{i-1}$ is an immediate predecessor of $F_i$ 
(see \cite{BB-S}).
\item
We call a cobordism (or a variety with a $K^*$-action) {\it collapsible} (see
also Morelli \cite{Morelli1}) if the relation $<$ on its set of
connected components of the fixed point set  
is an order. (Here an
order is just required to be transitive.)
\end{enumerate}
\end{definition}

\begin{remark} One can show (\cW) that a projective cobordism is collapsible. The collapsibility follows from the existence of a $K^*$-equivariant embedding into a projective space and direct computations
for the projective space. A similar technique works for a relatively projective cobordism.

\end{remark}
\begin{definition} A  function $\chi:\cC(X^{K^*})\to \ZZ$ is  {\it strictly increasing} if  $\chi(F)<\chi(F')$ whenever $F< F'$.
\end{definition}

\subsection{Existence of a strictly increasing function for $\bP^k$}

 The space $\bP^k=\bP(\bbA^{k+1})$ splits according to the weights as
\[\bP^k=\bP(\bbA^{k+1})=\bP(\bbA_{a_1}\oplus\cdots\oplus \bbA_{a_r}) \]
where $\ks$ acts on $\bbA_{a_i}$
 with the weight $a_i$. Assume that $a_1<\cdots< a_r$. Let $\ox_{a_i}=[x_{i,1},\dots, x_{i,{r_i}}]$
be the coordinates on $\bbA_{a_i}$. The action of $K^*$ is given by
\[t[\ox_{a_1},\dots,\ox_{a_r}]=[t^{a_1}\ox_{a_1},\dots,t^{a_r}\ox_{a_r}].\]
It follows that the fixed point components of $(\bP^k)^{K^*}$ are $\bP(\bbA_{a_i})$. We define a strictly increasing function $\chi_\bP:\cC(\bP^{K^*})\to\bZ$ by $$\chi_{\bP}(\bP(\bbA_{a_i}))=a_i.$$
We see that for $x=[x_{a_0},\dots, x_{a_r}]$, $\underset{t\to 0}{\lim} tx\in \bP(\bbA_{a_{\min}}), \underset{t\to \infty}{\lim} tx \in \bP(\bbA_{a_{\max}})$, where 
$${}a_{\max}=\max\{a\mid \ox_a\neq 0\},\quad
a_{\min}=\min\{a\mid\ox_a\neq 0\}.$$
Then $\bP(\bbA_{a_i})<\bP(\bbA_{a_j})$ iff $a_i<a_j$. 


\subsection{Existence of a strictly increasing function for a  compactified cobordism $\overline{B}$.}
Let $E$ be a $K^*$-invariant relatively very ample divisor for $\pi_B:\overline{B}\to X$.

For any $x\in F$, where $F\in\cC(\bar{B}^{\ks})$, we  find a semiinvariant function $f$ describing $-E=(f)$ in the neighborhood of $x\in B$. Then we put $\chi_E(F)=a$ to be the weight of the action $t(f)=t^af$ of $\ks$ on $f$. Note that $\chi_E:\cC(B^{K^*})\to\ZZ$ is locally constant so it is independent of the choice of $x\in F$. Let $f'$ be another function describing $E$ at $x$, with weight $a'$. Then the function $f'/f$ is invertible at $x$ so $t(f'/f)(x))=t^{a'-a}f'/f(x)=f'/f(tx)=f'/f(x)$. Then $a'-a=0$ and $a'=a$.

For any open affine set $U\subset X$ there exist $K^*$-semiinvariant  sections $s_0,\dots,s_K\in\Ga(\cO_{\overline B_U}(E))$ corresponding to rational $K^*$-semiinvariant functions $f_i$ (with the same weight) such that $(f_i)+E\geq 0$, which define a closed embedding $$\varp_U: \overline B_U\hookrightarrow \bP^k_U=\bP^k\times U,$$   where $\overline B_U=\pi_B^{-1}(U)$.
Every fixed point component $F$ on $\bar{B}_U$ is contained in $\bP_a\times U$.
%
%
  For any  $x\in F$ there exists a section $s_i$ such that $(s_i)=(f_i)+E=0$ in the neighborhood of $F$. Thus the section $s_i$ with weight $a_i$ is invertible at $x$. This implies that $F\cap B_U\subset \bP(\bbA_{a_i})\times U$. On the other hand, $(f_i)=-E$ and the weight of $f_i$ is $a_i$. Thus we get 
$\chi_E(F)=\chi_\bP(\bP(\bbA_{a_i}))=a_i.$
%
%
The function $\chi_\PP$ is strictly increasing, and the intersection of every component $F\in \bar{B}^{K^*}$  with $\bar{B}_U$ is contained in $\bP(A_{a})\times U$, where $\chi_E(F)=\chi_\bP(F)=a$. In particular we get
%
%
$\chi_E(F)<\chi_E(F')$ if $F<F'$ so $\chi_E$ is a strictly increasing function on $\bar{B}$. This implies
\begin{lemma} A compactified cobordism $\overline{B}$ is collapsible.
\end{lemma}

\subsection{Decomposition of a birational cobordism}
\begin{definition}(\cite{AKMW}, \cW) A cobordism  $B$ is {\it elementary} if for any $F\in \cC(B^{K^*})$ the sets $F^+$ and $F^-$ are closed. (In particular  any two distinct component $F, F'\in \cC(B^{K^*})$  are incomparable with respect to $>$.)\end{definition}

The  function $\chi_F$ defines a decomposition of  $\cC(B^{K^*})$ into elementary cobordisms
$$B_{a_i}:=B\setminus (\bigcup_{\chi_B(F)<a_i}F^-\cup \bigcup_{\chi_B(F)>a_i}F^+),$$
%
%
 where $a_1<\cdots<a_r$ are the values of $\chi_B$.
This yields
\begin{lemma}
\begin{enumerate}
\item
$(B_{a_1})_-=B_-$, $(B_{a_r})_+=B_+$.

\item $(B_{a_{i+1}})_-=(B_{a_i})_+= B\setminus (\bigcup_{\chi_B(F)\leq a_i}F^-\cup \bigcup_{\chi_B(F)\geq a_{i+1}}F^+)$.
\item $\chi(F)=a_i$ for any $F\in \cC(B_{a_i})$.
\item $(B_{a_i})_-=B_{a_i}\setminus  (\bigcup_{\chi_B(F)= a_i}F^+)$,  $(B_{a_i})_+=B_{a_i}\setminus  (\bigcup_{\chi_B(F)= a_i}F^-)$

\end{enumerate} 

\end{lemma}

\begin{figure}[ht]

\epsfysize=1in

	\epsffile{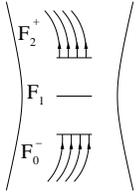}

	\caption{Elementary birational cobordism}\label{Fi:5}
\end{figure}

\begin{figure}[ht]

\epsfysize=1in

	\epsffile{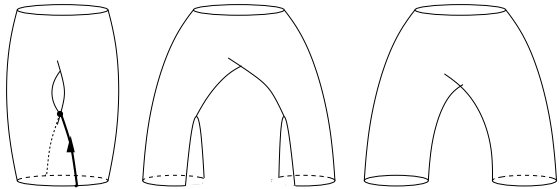}

	\caption{"Handle"-elemenatry cobordism in Morse Theory}\label{Fi:6}
\end{figure}

\subsection{Decomposition of $\bP^k$}
Set $\bbA_{\geq a_i}:=\bbA_{a_i}\oplus\cdots\oplus \bbA_{a_r}$,
$\bbA_{>a_i}:=\bbA_{a_{i+1}}\oplus\cdots\oplus \bbA_{a_r}$, and define $\bbA_{< a_i}$, $\bbA_{\leq a_i}$  analogously.

\begin{lemma}   $\bP(\bbA_{a_i})^+=\bP(\bbA_{\geq a_i})$ and
$\bP(\bbA_{a_i})^-=\bP(\bbA_{\leq a_i})$. Moreover if $F=\bP(F_A)\subset  \bP(\bbA_{a_i})$ is a closed subset then $F^+=\bP( F_A\oplus A_{a_i}\oplus\cdots\oplus \bbA_{a_r})$ is closed.
\end{lemma}

\begin{lemma} Set  $\bP_{a_i}:=\bP^k\setminus (\bigcup_{\chi_\bP(F)<a_i}F^-\cup \bigcup_{\chi_\bP(F)>a_i}F^+).$ 
%
%
 Then

${}\bP_{a_i}=\bP^k \setminus \bP(\bbA_{> a_i}) \setminus \bP(\bbA_{< a_i}),\quad\quad 
(\bP_{a_i})_+=\bP^k \setminus \bP(\bbA_{\geq a_i}) \setminus \bP(\bbA_{< a_i}),\quad\quad 
(\bP_{a_i})_-=\bP^k \setminus \bP(\bbA_{> a_i}) \setminus \bP(\bbA_{\leq a_i})$. 
\end{lemma}
\begin{lemma}$ \phi^{-1}_U(\bP_{a_i}\times U)=(B_U)_{a_i},$\quad $ \phi^{-1}_U((\bP_{a_i})_+ \times U)=((B_U)_{a_i})_+,$ \quad $ \phi^{-1}_U((\bP_{a_i})_- \times U)=((B_U)_{a_i})_-,$

\end{lemma} 
Combining these results gives us
\begin{lemma} 
\label{open}
The sets  $(B_a)_-$, $(B_a)_+$ and  $B_a$ are open in $B$. For any $F\in \cC(B_a)^{K^*}$, the sets $F^+, F^-$ are  closed.
\end{lemma}

\subsection{ GIT and existence of quotients for $\bP^k$.}

 The sets $\bP_{a_i}$ can  be interpreted in terms of Mumford's GIT theory. Any lifting of the action of $K^*$ on $\bP^k$ to $\bbA^{k+1}$ is called a {\it linearization}. Consider the twisted action on $\bbA^{k+1}$,
\[t_r(x)=t^{-r}\cdot t(x).\]
The twisting  does not change the action on $\bP(\bbA^{k+1})$ and defines different linearizations. 
If we compose the action with a group monomorphism $t\mapsto t^k$ the weights of the new action $t^k(x)$ will be multiplied by $k$. The good and geometric quotients for $t(x)$ and $t^k(x)$ are the same.
Keeping this in mind it is convenient to allow linerizations with rational
weights.
\begin{definition} The point $x\in\bP^{k}$ is {\it semistable} with respect to $t_r$, written $x\in(\bP^{k},t_r)^{ss}$,
%
%
  if there exists an invariant section $s\in \Ga(\cO_{{\bP}^{k+1}}(n)^{t_r})$ such that $s(x)\neq 0$.
\end{definition}

\begin{lemma}(\cite{AKMW}) $\bP_{a_i}=(\bP^k, t_{a_i})^{ss}$, $(\bP_{a_i})_-=(\bP^k, t_{a_i-\frac{1}{2}})^{ss}$, $(\bP_{a_i})_+=(\bP^k, t_{a_i+\frac{1}{2}})^{ss}$.
\end{lemma}

\begin{proof} $x\in \bP_{a_i}$ iff either $\ox_{a_i}\neq 0$  or $\ox_{a_{j_1}}\neq 0$ and $\ox_{a_{j_2}}\neq 0$ for $a_{j_1}<r=a_i<a_{j_2}$.
%
%
 In both situations we find a nonzero $t_r$-invariant section $s_i=x_{i}$ or $s_{j_1j_2}=x_{j_1}^{b_1}x_{j_2}^{b_2}$ for suitable coprime $b_1$ and $b_2$. 

$x\in (\bP_{a_i})_-$ iff    $\ox_{a_{j_1}}\neq 0$ and $\ox_{a_{j_2}}\neq 0$ for $a_{j_1}<a_i\leq a_{j_2}$ (or equivalently $a_{j_1}<r=a_i-1/2< a_{j_2}$). As before there is a nonzero $t_r$-invariant section  $x_{j_1}^{b_1}x_{j_2}^{b_2}$ for suitable coprime $b_1$ and $b_2$.
\end{proof}

It follows from GIT theory that $(\bP^k,t_r)^{ss}/\/\ks$ exists and it is a projective variety. Moreover we get

\begin{lemma} \label{invariant} Let $X_r:=(\bP^k,t_r)^{ss}\subset \bP^k$. For a sufficiently divisible $n\in\bN$,  the invariant sections $x^{\alpha_0},\dots,=x^{\alpha_\ell}\in \Ga(\cO_{X}(n)^{t_r})$ define the morphism $\psi: X\to\bP^\ell=\Proj(K[s_0,\ldots,s_\ell])$, by putting $s_i\mapsto x^{\alpha_i}$. Then
$X_r//\ks\cong\psi(X_r)\subset \bP^\ell$. Moreover if $\pi: X_r\to X_r//K^*$ is the quotient morphism (determined by $\psi$) then the push-forward $\pi_*(\cO_X(n)^{t_r})$ of the sheaf $\cO_X(n)^{t_r}$ of $\cO_{X_r}^{K^*}$-modules on $X_r$ is a very ample line bundle on $X_r//K^*$.
\end{lemma}

\begin{proof}
{(1)} Let $x_{i}$ be a coordinate on $\bP^k$ with $t_r$-weight $0$. The section  $x_i^n\in\Ga(\cO(n)^{t_r})$ corresponds to a coordinate $s$ on $\bP^\ell$ and let $\bbA^\ell_s:=\{p\in \bP^\ell\mid s(p)\neq 0\}\subset \bP^\ell$ be the open affine subset. Thus the inverse image $\psi^{-1}(\bbA^\ell_s)$ equals to $U_i:=\{x\mid x_{i}\neq 0\}$ and the morphism
$\psi_{|U_i}: U_i\to \bbA^\ell_s$ is given by $s_j/s\mapsto x^{\alpha_j}/x_i^n$. If $n\in\bN$ is sufficiently divisible then $K[U_i]^{t_r}$ is generated by all the monomials $x^{\alpha_j}/x_i^n$. We have a surjection $K[ \bbA^\ell_s]\to K[U_i]^{t_r}\subset K[U_i]$ corresponding to the embedding of the quotient $U_i\to U_i//{\ks}\subset \bbA^\ell_s$.

{(2)} Let $x_{j_1}, x_{j_2}$ be two coordinates on $\bP^k$ whose $t_r$-weights have opposite signs. Then  the $t_r$-invariant section $(x_{j_1}^{b_1}\cdot x_{j_2}^{b_2})^{n/{(b_1+b_2)}}\in \Ga(\cO(n)^{t_r})$ for suitable coprime $b_1,b_2$ corresponds to a coordinate $s$ on $\bP^\ell$. The inverse image  $\psi^{-1}(\bbA^\ell_s)$ is given by $U_{j_1j_2}:=\{x\mid x_{j_1}\cdot x_{j_2}\neq 0\}$. We get the quotient morphism again: 
$\psi_{|U_{j_1j_2}}:  U_{j_1j_2}\to U_{j_1j_2}//{\ks}=\Spec  K[U_{j_1j_2}]^{t_r}\subset \bbA^\ell_s.$

Note that the monomials $x^{\alpha_i}\in \Ga(\cO(n)^{t_r})$ are the products of $x_i$ as in (1) and $x_{j_1}^{b_1}\cdot x_{j_2}^{b_2}$ as in (2). Thus $U_i$ and  $U_{j_1j_2}$ cover the image of $\psi$ so $\psi((\bP^k, t_{r})^{ss})\subset \bP^\ell_s$ is a closed subset  and $\psi$  defines a quotient morphism.
\end{proof}

\begin{corollary}
There exist quotients $\pi_{a_i}:\bP_{a_i}\to \bP_{a_i}//{\ks}$, $\pi_{a_{i-}}=\bP_{a_{i-}}\to (\bP_{a_{i}})_-/{\ks}$, $\pi_{a_{i+}}=(\bP_{a_{i}})_+\to (\bP_{a_{i}})_+//{\ks}$. 
\end{corollary}

\subsection{Existence of quotients for  $\bar{B}$}
\begin{lemma} There exist quotients $\pi_{a_i}:B_{a_i}\to B_{a_i}//{\ks}$, $\pi_{a_{i-}}=B_{a_{i-}}\to B_{a_{i-}}/{\ks}$, $\pi_{a_{i+}}=B_{a_{i+}}\to B_{a_{i-}}//{\ks}$. Moreover the induced morphisms $B_{a_i}//{\ks}\to X,  B_{a_{i-}}/{\ks}\to X,  B_{a_{i+}}/{\ks}\to X$  are projective
\end{lemma}
 \begin{proof}  Since $B_U\subset\bP^k\times U$ and $(B_U)_a\subset \bP_a\times U$ are closed subvarieties  the quotients $(B_U)_{a_i}//K^*$, $(B_U)_{a+}/{\ks}$, ${(B_U)_{a-}}/{\ks}$ exist for any open affine $U\subset X$. Glueing these together  defines the global quotients $(B)_{a_i}//K^*$, $(B)_{a+}/{\ks}$, ${(B)_{a-}}/{\ks}$. 
 Consider twisted action $t_r$ on the line bundle $E$.
  By Lemma \ref{invariant}, $(\pi_a)_*(\cO_{B_a}(nE)^{t_a})$ is relatively very ample on $B_{a}// {\ks} \to X$. Analogously for 
$B_{a_+}// {\ks}$ and $B_{a_-}// {\ks}$. 
 \end{proof}
\begin{lemma} The open embeddings $(B_a)_-,(B_a)_+\subset B_a$ define the factorization $$ \begin{array}{ccccc}
 ({B}_{a_i})_-/K^* & & \stackrel{\varphi_i}{\das} & &({B}_{a_i})_+/K^* \\
 & \searrow & & \swarrow & \\
& & {B}_{a_i} /\!/ K^* & & \end{array}$$
which is an isomorphism over $\pi_a(B_a^{K^*})\simeq B_a^{K^*}\subset  B_a//K^*$.
\end{lemma}

As a corollary from the above we get
\begin{proposition} \cW \,\label{deco} There is a factorization of the projective morphism $\phi: Z\to X$ given by
$$Z=(B_{a_1})_-/K^* \dashrightarrow (B_{a_1})_+/K^*= (B_{a_2})_-/K^* \dashrightarrow
\ldots (B_{a_{k-1}})_+/K^*= (B_{a_k})_+/K^* \dashrightarrow
({B_{a_k}})_+/K^*=X.$$ 
\end{proposition}

\subsection{Local description of elementary cobordisms}
\begin{proposition} (\cW) \label{local}
Let $B_{a}$ be a smooth elementary cobordism. Then for
any $x\in B^{\ks}_a$ there exists an invariant neighborhood
$V_x$ of $x$ and a
$K^*$-equivariant \'etale morphism  (i.e. locally analytic isomorphism) $\phi :V_x\rightarrow
\Tan_{x,B} $,
where $\Tan_{x,B}\simeq \bbA_K^n$
%
%
 is the tangent space with the
induced linear $K^*$-action, such that in the diagram
\[\begin{array}{rccccccc}
&{(B_a)}_-/K^*& \supset &V_x//K^* \times_{ \Tan_{x,B}//K^*}
{(\Tan_{x,B})}_-/K^*& \simeq &{V_x}_-/K^*& \rightarrow
&{(\Tan_{x,B})}_-/K^* \\

& \downarrow&&&& \downarrow & &\downarrow \\

&B_a//K^*    &&\supset&&V_x//K^*    & \rightarrow  &\Tan_{x,B}//K^* \\
&\uparrow&&&&\uparrow& &\uparrow \\
&{(B_a)}_+/K^*&\supset&V_x//K^*  \times_  { \Tan_{x,B}//K^*}
{(\Tan_{x,B})}_+/K^* & \simeq &{V_x}_+/K^*&\rightarrow & {(\Tan_{x,B})}_+/K^*\\

\end{array}\]
%
%
 the vertical arrows are defined by open embeddings
and the horizontal morphisms are defined by $\phi$ and are \'etale.
\end{proposition}

\begin{proof} By Lemma \ref{open}, for any irreducible component $F\in \cC(B_a^{K^*})$ the sets $F^+$ and $F^-$ are closed.  Let $U$ be an open $K^*$-equivariant neighborhood of $x\in \in B_a^{K^*}$, disjoint from the  closed sets  $F^+$ and $F^-$, where $F\in \cC(B_a^{K^*})$,  does  not pass through $x$.
By taking local semiinvariant parameters at the point $x$ one
can construct an equivariant morphism $\phi:U_x\rightarrow
\Tan_{x,B}\simeq \bbA^n_k$ from some
open affine invariant neighborhood $U_x\subset U$ such that $\phi$ is \'etale
at $x$.
By Luna's Lemma (see [Lu], Lemme 3 (Lemme Fondamental))
there exists an invariant affine
neighborhood $V_x\subseteq U_x$ of the point $x$ such that
$\phi_{\mid V_x}$ is \'etale, the induced map
$\phi_{\mid V_x/K^*}:V_x//K^* \rightarrow \Tan_{x,B}//K^*$ is \'etale
and $V_x\simeq V_x//K^* \times_{\Tan_{x,B}//K^*} \Tan_{x,B}$.
This defines the isomorphisms $V_x//K^* \times_{ \Tan_{x,B}//K^*}
{(\Tan_{x,B})}_-/K^* \simeq {V_x}_-/K^*$.
It follows that the irreducible components of $B_a^{K^*}$ are smooth and disjoint and the sets $F^+$ and $F^-$, where $F\in \cC(B_a^{K^*})$, are irreducible. Let $F_0\in \cC(B_a^{K^*})$ be the component through $x$.
 Note that $V_x\cap (\bigcup_{F\in \cC(B_a^{K^*})}F^+)=V_x\cap F_0^+= (V_x\cap F_0)^+$. (Both sets are closed and irreducible (Lemma \ref{fix}.)) Thus $(V_x)_-=V_x\cap (B_{a})_-$ and we get the horizontal inclusions.
\end{proof}
\begin{remark} It follows from the above thet the birational maps $(B_a)_-/K^*\dashrightarrow (B_a)_+/K^*$ are 
locally described by Example \ref{main}. Both spaces have cyclic singularities and differ by  the composite of a weighted blow-up and a weighted
blow-down. To achieve the factorization we need to desingularize quotients as in for instance case 1 of the example. It is hopeless to modify weights
by  birational modification of smooth varieties. Instead we want to view Example \ref{main} from the perspective of toric varieties. 
\end{remark}
\section{Toric varieties}

\subsection{Fans and toric varieties}
Let $N\simeq
{\ZZ}^k$ be a lattice contained in 
the vector space $N^{\QQ}:=N\otimes {\QQ}\supset N$.  
\begin{definition} (\cite{Danilov1},
\cite{Oda})\label{de: fan} By a {\it fan}
$\Sigma $ in
$N^{\QQ}$
 we mean a finite collection of finitely 
generated strictly convex cones $\sigma$ in $N^{\QQ}$ such 
that 

$\bullet$ any face of a cone in $\Sigma $ belongs to $\Sigma$,

$\bullet$ any two cones of $\Sigma $ intersect in a common face. 
\end{definition}

If $\sigma$ is a face of $\sigma'$ we shall write $\sigma\preceq\sigma'$.

We say that a cone
$\sigma$ in $N^{\QQ}$ is {\it nonsingular} if it is generated by a part of a basis of the lattice
$e_1,\ldots,e_k\in N$, written $\sigma=\langle e_1,\ldots,e_k\rangle$.
%
%
 A cone $\sigma$ is {\it simplicial} if it 
is generated over $\bQ$ by linearly
independent integral vectors $v_1,\ldots,v_k$, written $\sigma =\langle v_1,\ldots,v_k\rangle $

\begin{definition}\label{de: star} Let $\Sigma$ be a fan
and $\tau \in \Sigma$. The {\it star} of the
cone $\tau$ and the {\it closed star} of $\tau$ are
defined as follows:
 $${\rm Star}(\tau ,\Sigma):=\{\sigma \in \Sigma\mid 
\tau\preceq \sigma\},$$ 
$$\overline{{\rm Star}}(\tau ,\Sigma):=\{\sigma \in
\Sigma\mid  \sigma'\preceq \sigma  \mbox{ for some }  \sigma'\in
{\rm Star}(\tau ,\Sigma)\}.$$ 
\end{definition} 

 To a fan $\Sigma $ there is associated a toric variety
$X_{\Sigma}\supset {\bf }T$, i.e. a normal variety on which a torus
$T$ acts effectively with an open dense orbit 
(see \cite{KKMS}, \cite{Danilov2}, \cite{Oda}, \cite{Fulton}). To 
each cone $\sigma\in \Sigma$
corresponds an  open affine invariant subset
$X_{\sigma}$ and its unique closed orbit $O_{\sigma}$. The
orbits in the
closure of the  orbit $O_\sigma$ correspond to the cones of 
${\rm Star}(\sigma ,\Sigma)$.  In particular, $\tau\preceq\sigma$ iff $\overline{O_{\tau}}\supset O_\sig$. (We shall also denote the closure of $O_\sig$ by $\cl(O_\sig)$.) 

The fan  $\Sigma $ is {\it nonsingular} (resp. {\it simplicial}) if all its cones are nonsingular (resp.  simplicial). Nonsingular fans correspond to nonsingular varieties. 

Denote by $$M:={\rm Hom}_{alg.gr.}(T,K^*)$$  the lattice of
group homomorphisms to $K^*$, i.e.  characters of $T$.
The dual lattice $
{\rm Hom}_{alg.gr.}(K^*,T)$ of $1$-parameter subgroups of $T$ can be 
identified with the lattice $N$.
Then the
vector space $M^{\QQ}:=M\otimes{\QQ}$ is dual to $N^{\QQ}=
N\otimes{\QQ}$.

 The elements $F\in M=N^*$ are functionals on $N$ and integral functionals on $N^{\QQ}$.
For any $\sigma\subset N^{\QQ}$ we denote by
$$\sigma^\vee:=\{F\in M \mid F(v) \geq 0 \,\, {\rm
for\,\, 
any} \,\,\,    v\in \sigma\}$$ 
 the set of integral vectors of the dual cone to $\sigma$. Then the
ring of regular functions $K[X_\sigma]$ is $K[\sigma^\vee]$.  

We call a vector $v\in N$ {\it primitive} if it generates the 
sublattice ${\QQ}_{\geq 0}v\cap N$. Primitive vectors correspond to $1$-parameter monomorphisms.

For any $\sigma\subset N^{\QQ}$ set
$$\sigma^\perp:=\{m\in M \mid (v,m) = 0 \,\, {\rm
for\,\, 
any} \,\,\,    v\in \sigma\}.$$ 
 The latter set represents all  invertible characters
on $X_\sigma$. All noninvertible characters  are in $\sigma^\vee\setminus\sigma^\perp $ and vanish on $O_\sigma$.
The ring of regular functions on $O_\sigma\subset X_\sigma$ can be
written as $K[O_\sigma]=K[\sigma^\perp]\subset
K[\sigma^\vee]$.

\subsection{Star subdivisions and blow-ups}

\begin{definition}(\cite{KKMS}, \cite{Oda},
\cite{Danilov2}, \cite{Fulton}) A
{\it birational toric morphism} or simply a {\it toric morphism} of toric
varieties $X_\Sigma \to X_{\Sigma'}$ is a 
morphism identical on $T\subset X_\Sigma, X_{\Sigma'}$. 
\end{definition}

By the {\it support} of a fan $\Sigma$ we mean the union of all
its faces, 
$|\Sigma|=\bigcup_{\sigma\in \Sigma}\sigma$.

\begin{definition} (\cite{KKMS}, \cite{Oda},
\cite{Danilov2}, \cite{Fulton}) 
A {\it subdivision} of a fan
$\Sigma$ is a fan $\Delta$ such that $|\Delta|=|\Sigma|$
and any cone $\sigma\in
\Sigma $ is a union of cones $\delta\in
\Delta$. 
\end{definition}

\begin{definition}\label{de: star subdivision} Let $\Sigma$ be a fan and
$\varrho$ be a ray passing in the
relative interior of $\tau\in\Sigma$. Then the {\it star
subdivision} $\varrho\cdot\Sigma$ of $\Sigma$ with respect to
$\varrho$ is defined to be
$$\varrho\cdot\Sigma=(\Sigma\setminus {\rm Star}(\tau ,\Sigma) )\cup
\{\varrho+\sigma\mid   \sigma\in \overline{\rm Star}(\tau
,\Sigma)\setminus {\rm Star}(\tau
,\Sigma)\}.$$ If $\Sigma$ is nonsingular, i.e. all its cones are
nonsingular, $\tau=\langle v_1,\ldots,v_l\rangle $ and 
$\varrho=\langle v_1+\cdots+v_l\rangle $
then we call the star
subdivision $\varrho\cdot\Sigma$ {\it nonsingular}. 
\end{definition}

\begin{proposition} (\cite{KKMS}, \cite{Danilov2},
\cite{Oda}, \cite{Fulton}) Let
$X_\Sigma$ be a toric variety. There is a 1-1 correspondence
between subdivisions of the fan $\Sigma$ and proper toric
morphisms $X_{\Sigma'} \to X_{\Sigma}$.
\end{proposition}

\begin{remark} Regular star subdivisions from
\ref{de: star subdivision} correspond to blow-ups of smooth varieties
at closures of orbits (\cite{Oda}, \cite{Fulton}). Arbitrary
star subdivisions correspond to blow-ups of some ideals
associated to valuations  (see Lemma \ref{blow-up}).
\end{remark}

\section{Polyhedral cobordisms of Morelli}
\subsection{Preliminaries}
%
%
By $N^{\QQ+}$
%
%
 we shall denote a vector space $N^{\QQ+}\approx \QQ^k$ containing a
lattice $N^+\simeq \bZ^k$, together with a primitive vector ${v_0}\in N^+$ and the
canonical projection
$$\pi:N^{\QQ+}\to N^{\QQ}\simeq N^{\QQ+} / \bQ\cdot{v_0}.$$

\begin{definition}(\cite{Morelli1})
A cone $\sigma\subset N^{\QQ+}$ is {\it{$\pi$-strictly convex}} if $\pi(\sigma)$ is strictly convex (contains no line). A fan $\Sigma$ is $\pi$-strictly convex if it consists of $\pi$-strictly convex cones.
\end{definition}

In the following all the  cones in $N^{\QQ+}$ are assumed to be $\pi$-strictly convex and simplicial. The
$\pi$-strictly convex cones $\sigma$ in $N^{\QQ+}$ split into two categories.
\begin{definition}
A cone $\sigma\subset N^{\QQ+}$ is called  {\it independent}
if the restriction of $\pi$ to $\sig$ is a linear isomorphism
(equivalently $v_0\not\in{\rm span}(\sigma)$).
%
%
A cone $\sigma\subset N^{\QQ+}$  is called  {\it
dependent} if the restriction of $\pi$ to $\sig$ is a lattice submersion which is not an isomorphism
(equivalently $v_0\in{\rm span}(\sig)$).

A dependent cone is called a
{\it circuit} if all its proper faces are independent. 
\end{definition}
\begin{lemma}\label{circuit} Any dependent cone $\sigma$ contains a unique circuit $\delta$.
\end{lemma}

\subsection{$K^*$-actions and $N^{\QQ+}$}
The vector ${v_0}=(a_1,\dots,a_k)\in N^{\QQ+}$ defines a 1-parameter subgroup $t^{v_0}:=t^{a_1}_1\cdots t_k^{a_k}$
%
%
 acting on $T$ and all toric varieties $X\supset T$. Denote by $M^+$ the lattice dual to $N^+$. Then the lattice $N:={N^+}/\bZ\cdot{v_0}\subset N^\QQ:=N\otimes \bQ$
 is dual to the lattice $M:=\{a\in M^+\mid(a,{v_0})=0\}\subset M^\QQ:=M\otimes\bQ$ of all the characters invariant with respect to the group action. The natural projection of cones $\pi :\sigma\to \pi(\sigma)\subset N^\QQ$ defines the good quotient morphism
$$X_\sigma=\Spec  K[\sigma^\vee]\to X_\sigma//K^*
=\Spec  K[\sigma^\vee\cap M]
=\Spec  K[{\pi(\sigma)}^\vee]=X_{\pi(\sigma)}.$$

\begin{lemma}\label{fixed3} 
A cone $\sigma$ is independent iff the geometric quotient $X_\sigma\to X_{\sigma}/K^*$ exists or alternatively if $X_\sigma$ contains no fixed points. The cone $\sigma$ is dependent if $O_\sigma$ is a fixed point set.
\end{lemma}
\begin{proof} Note that the set $X_\sigma^{K^*}$ is closed and if it is nonempty then it contains  $O_\sigma$. Then a point  $p\in O_\sigma$ is  fixed, i.e.  $t^{v_0} p=p$, iff for all functionals $F\in \sigma^\perp $   (i.e. $x^F(p)\neq 0)$ we have $x^F(p)=x^F (t^{v_0}p)=t^{F(v_0)}x^F(p)$.

 Then for all $F\in\sig^\perp\subset\rm{span}(\sig)^\perp$ we have 
 $F({v_0})=0$ so  ${v_0}\in \rm{span} (\sigma)$.
\end{proof}
\begin{remark} In Example \ref{main} $X_\sigma=\bbA^n$ is a cobordism iff $\sig$ is $\pi$-strictly convex. 
\end{remark}
\begin{corollary}\label{fixed2} A cone $\delta\in \Sigma$ is a circuit if and only if $O_\delta$ is the generic orbit of some  $F\in \cC(X_\Sigma^{K^*})$. \end{corollary}
\begin{proof}$O_\sigma$ is fixed with respect to the action of $K^*$ if $\sigma$ is dependent. Thus $O_\sigma\subset \overline{ O_\delta}$ where $\delta$ is the unique circuit in $\sigma$ (Lemma \ref{circuit}).
\end{proof}

\subsection{Morelli cobordisms}

\begin{definition}(Morelli \cite{Morelli1},
\cite{Abramovich-Matsuki-Rashid})\label{de: Morelli} 
A fan $\Sigma$ in $N^{{\QQ}+}\supset N^+$
 is
called a {\it polyhedral cobordism} or simply a cobordism if
the sets of cones 
 $$\partial_-(\Sigma): =\{\sigma \in\Sigma\mid \mbox{there  exists } p\in \inte(\sig)\mbox{ such that }
 p-\epsilon\cdot v_0\not\in\ |\Sigma| \ \mbox{for all small} \  \epsilon>0\}, 
$$  $$\partial_+(\Sigma): =\{\sigma \in\Sigma\mid \mbox{there  exists } p\in \inte(\sig)\mbox{ such that }
 p+\epsilon\cdot v_0\not\in\ |\Sigma| \ \mbox{for all small} \ \epsilon>0\} \quad\quad\quad{}$$ are
subfans of $\Sigma$ and
$\pi(\partial_-(\Sigma)):=\{\pi(\tau)\mid\tau\in \partial_-(\Sigma)\}$ and 
$\pi(\partial_+(\Sigma)):=\{\pi(\tau)\mid\tau\in \partial_+(\Sigma)\}$
are fans in $N^{\QQ}$.

\end{definition}

\bigskip
\subsection{Dependence relation}\label{dep}

Let $\sigma=\langle v_1,\dots,v_k\rangle$ be a dependent (simplicial) cone. Then, by definition $v_0\in{\rm{span}} (v_1,\dots,v_k)$
 where $v_1,\dots,v_k$ are linearly independent. There exists a unique up to rescaling integral relation
\[r_1v_1+\dots +r_kv_k=av_0,\quad \mbox{where}\quad a>0.\quad (*) \] 
\begin{definition}(\cite{Morelli1}) The rays of $\sigma$ are called {\it positive, negative} and {\it null } vectors, according to the sign of the coefficient in the defining relation.
\end{definition}
\begin{remark}
Note that the relation $(*)$ defines a unique relation $$r'_1w_1+\dots+r'_kw_k=0\quad\quad (**)$$ where $w_i$ are  generating vectors in the rays $\pi\left(\langle v_i\rangle\right)$, $r'_iw_i=r_i\pi(v_i)$. 
In particular $r'_i/r_i>0$.
\end{remark}
\begin{lemma}\label{sign} Let $\sigma=\langle v_1,\dots, v_k\rangle$ be a dependent cone. Then an independent face $\tau$ is in $\partial_+(\sigma)$ (resp. $\tau\in\partial_+(\sigma)$) if 
$\tau$ is a face of $\langle v_1,\ldots, \check{v_i},\ldots, v_k \rangle$  for some index $i$ such that $r_i<0$ (resp. $r_i>0$).
\end{lemma}
\begin{proof} By definition $\tau\in\partial_+(\sigma)$  there exists $p\in\rm{int} (\tau)$ such that for any sufficiently small $\epsilon>0$, 
$p+\epsilon {v_0}\notin\sigma$.
Write $p=\sum\alpha_iv_i=\sum_{r_i>0}\alpha_iv_i+\sum_{r_i<0}\alpha_iv_i+\sum_{r_i=0}\alpha_iv_i$, where $\alpha_i\geq 0.$ Then one of the coefficients in $$p+\epsilon {v_0}=\sum_{r_i>0}(\alpha_i+r_i\epsilon)v_i+\sum_{r_i<0}(\alpha_i+r_i\epsilon)v_i+\sum_{r_i=0} (\alpha_i+r_i\epsilon)v_i.$$
is negative for small $\epsilon>0$. This  is possible if $\alpha_i=0$ for some index $i$ with $r_i<0$.
\end{proof}

 \begin{lemma} A cone  $\tau$ is in $\partial_+(\sigma)$ iff  there exists ${F\in\sigma^\vee\cap \tau^\perp}$ such that $F({v_0})<0$.
\end{lemma}
\begin{proof}If $\tau\in\partial_+(\sigma)$ then there exists $p\in\rm{int} (\tau)$ for which $ p+\epsilon{v_0}\notin\sigma$.
 Hence there exists $F\in\sigma^\vee$ such that $F(p+\epsilon{v_0})<0$ for small $\epsilon>0$. Then $F(p)=0$ and $F({v_0})<0$. Since $p\in\rm{int} (\tau)$ we have $F_{|\tau}=0$.
\end{proof}

\begin{corollary}$\partial_+(\sigma)$ (resp. $\partial_-(\sigma)$) is a fan.
\end{corollary}
\begin{proof} By the lemma above, if $\tau\in\sigma^+$ then  every face 
$\tau'$ of $\tau$ is in $\sigma^+$.
\end{proof}
\begin{lemma}\label{sigma2}Let $\sigma$ be a dependent cone in $N^{\QQ+}$. Then $B:=X_\sigma$ is a birational cobordism such that
\item{$\bullet$} $(X_\sigma)_+=X_{\partial_-(\sigma)}$, $(X_\sigma)_-=X_{\partial_+(\sigma)}$.
\item{$\bullet$} $(X_\sigma)_+/K^*\cong X_{\pi(\partial_-(\sigma))}$, $(X_\sigma)_-/K^*\cong X_{\pi(\partial_+(\sigma))}$.
%
%
\item{$\bullet$}  $\pi(\partial_-(\sigma))$ and $\pi(\partial_+(\sigma))$ are both decompositions of $\pi(\sigma)$.
\item{$\bullet$}  There is a factorization into a sequence of proper morphisms  $(X_\sigma)_+/K^*\to(X_\sigma)//K^*\leftarrow (X_\sigma)_-/K^*$.
\end{lemma}
\begin{proof} We have $p\in O_\tau$ where $O_\tau\subset(X_\sigma)_-$ iff  $ \lim t^{v_0} p\notin X_\sigma$. This is equivalent to existence of a functional ${F\in\sigma^{\vee}}$ for which $ x^F(t^{v_0} p)=t^{F({v_0})}x^F(p)$ has a pole at $t=0$. This means exactly that $x^F(p)\neq 0$ and $F({v_0})<0$. The last condition says $ F_{|\tau}=0$ and $F({v_0})<0$, which is equivalent to $\tau\in\partial_+(\sigma)$.

Let $x\in\pi(\sigma)$. Then $\pi^{-1}(x)\cap\sigma$ is a line segment or a point. Let $y=\sup \{\pi^{-1}(x)\cap\sigma\}$. Then $y\in\inte(\tau)$, where $\tau\prec\sigma$ and $y+\epsilon{v_0}\notin\sigma$, which implies that $\tau\in\partial_+(\sigma)$. Thus every point in $\pi(\sigma)$ belongs to a relative interior of a unique cone $\pi(\tau)\in \pi(\partial_+(\sigma))$. Since $\pi_{|\tau}$ is a linear isomorphism and $\partial_+(\sigma)$ is a fan, all faces of $\pi(\tau)$ are in $\pi(\partial_-(\sigma))$. Finally, $\pi(\partial_+(\sigma))$ and $\pi(\partial_-(\sigma))$ are both decompositions of $\pi(\sigma)$ corresponding to toric varieties $(X_\sigma)_-/K^*= X_{\pi(\partial_+(\sigma))}$ and $(X_\sigma)_+/K^*= X_{\pi(\partial_-(\sigma))}$.
%
%
\end{proof} 
The Lemmas \ref{sign} and \ref{sigma2} yield
\begin{lemma} \label{sigma} $B=X_\sigma$ is an elementary cobordism with a single
fixed point component $F:=\overline{O_\delta}$, where $\del=\langle v_i\mid r_i\neq 0\rangle $ is a circuit. Moreover $(X_\sigma)_+=X_{\partial_-(\sigma)}=X_\sigma\setminus \overline{O_{\sig_+}}$, where $$\sig_+:=\langle v_i\mid r_i> 0\rangle, \quad \sig_-:=\langle v_i\mid r_i< 0\rangle.$$  In particular
$F^+=(\overline{O_\delta})^+=\overline{O_{\sig_+}}$,\quad  $F^-=(\overline{O_\delta})^-=\overline{O_{\sig_-}}$.
\end{lemma}

\subsection{Example \ref{main} revisited}
The cobordism $X_\sigma$ from the lemma generalizes the cobordism $B= \bbA^{l+m+r}_K\supset T=(K^*)^{l+m+r}$  from Example \ref{main}. 
The action of $K^*$  determines
a $1$-parameter subgroup of $T$ which corresponds to a vector $v_0=[a_1,\ldots,a_l,-b_1,\ldots,-b_m,0,\ldots,0].$
The cobordism $B$  is associated with  a nonsingular
cone $\Delta\subset N_{\QQ}$, while $B_-$ and $B_+$ correspond to the fans $\partial_+(\Delta)$ and $\partial_-(\Delta)$ consisting of the
faces of $\Delta$ visible from above   and 
 below respectively.

 The  quotients $B_+/K^*$ , $B_-/K^*$ and $B//K^*$
 are toric varieties corresponding  to the fans
$\pi(\partial_+(\Delta))=\{\pi(\sigma)\mid \sigma\in\partial_+(\Delta)\}$,
$\pi(\partial_-(\Delta))=\{\pi(\sigma)\mid \sigma\in\Delta^-\}$ and $\pi(\Delta)$
respectively, where $\pi$ is the projection defined by $v_0$.

The relevant birational map $\phi : B_-/K^* \mathrel{{-}\,{\rightarrow}} B_+/K^*$
for $l,m\geq 2$ is a
toric flip associated with a bistellar operation replacing
the triangulation $\pi(\partial_-(\Delta))$ of the cone $\pi(\Delta)$
with  $\pi(\partial_+(\Delta))$.


\begin{figure}[ht]

\epsfysize=1.2in

	\epsffile{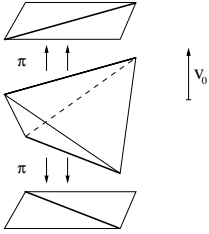}

	\caption{Morelli cobordism}\label{Fi:7}
\end{figure}\subsection{ $\pi$-nonsingular cones}
 \begin{definition}(Morelli)\label{mo}
 An independent cone $\tau$  is $\pi$-{\it nonsingular} if
 $\pi(\tau)$ is nonsingular. A fan $\Sigma$ is   
$\pi$-{\it nonsingular} if all independent cones in $\Sigma$ 
 are $\pi$-nonsingular. In particular a dependent cone $\sigma$ is 
$\pi$-{\it nonsingular} if all its independent faces are $\pi$-nonsingular.
\end{definition}
\begin{lemma}
 Let $\sigma=\langle v_1,\ldots,v_k\rangle$ be a dependent cone and  $w_i$ be primitive generators of the rays $\pi(v_i)$. Let $\sum r'_i w_i=0$ be the unique relation (**) between vectors $w_i$.
Then the  ray $\vr:=\pi(\sigma_+)\cap \pi(\sigma_-)$ is generated by the vector $\sum\limits_{r'_i>0}r'_i w_i=\sum\limits_{r'_i<0}-r'_iw_i$ and $\vr\cdot\pi(\partial_+(\sigma))=\vr\cdot\pi(\partial_-(\sigma))$. 
If $\sigma$ is a $\pi$-nonsingular dependent cone then the ray $\vr$ defines nonsingular star subdivisions of $\pi(\partial_+(\sigma))$ and $\pi(\partial_-(\sigma))$.

\end{lemma}
\begin{proof}
Note that $\pi(\partial_+(\sigma))\setminus \pi(\partial_-(\sigma))$ are exactly the cones containing $\pi(\sigma_+)$. That is, $\pi(\partial_+(\sigma)) \setminus \pi(\partial_-(\sigma))=\rm{Star} (\pi(\sigma_+),\pi(\partial_+(\sigma)))$.
%
%
This gives $\vr\cdot\pi(\partial_+(\sigma))=(\pi(\sigma_+)\cap\pi(\sigma_-))\cup \{\vr+\tau \mid  \tau\in\pi(\sigma_+)\cap\pi(\sigma_-)\}=\vr\cdot\pi(\sigma_-)$. Assume now that $\sig$  is $\pi$-nonsingular and all the coefficients $r'_i$ are coprime. By Lemma \ref{sign} and the $\pi$-nonsingularity the set of vectors $w_1,\ldots,\check{w_i},\ldots,w_k$ where $r'_i\neq 0$ is a basis of the lattice $\pi(\sigma)\cap N$. Thus
every vector $w_i$, where $r'_i\neq 0$, can be written as  an integral  combination of others. Since the relation (**) is unique it follows that the coefficient $r'_i$ is equal to $\pm 1$.
 Thus $\vr$ is generated by the vector $\sum\limits_{r'_i>0} w_i=\sum\limits_{r'_i<0}w_i$ and determines nonsingular star subdivisions.
\end{proof}
\begin{corollary} \label{fact} If  $\sigma$ is dependent then
there exists  a factorization  
$$(X_\sig)_-/K^* \buildrel \phi_-\over\longleftarrow \Ga((X_\sig)_-/K^*,(X_\sig)_+/K^*)\buildrel \phi_+\over\longrightarrow (X_\sig)_+/K^*,$$
 where $\Ga((X_\sig)_-/K^*,(X_\sig)_+/K^*)$ is the normalization of the graph of $(X_\sig)_-/K^*\to (X_\sig)_+/K^*$. If $\sig$ is  $\pi$-nonsingular the morphisms
$ \phi_-, \phi_+$ are blow-ups of smooth centers.
\end{corollary}
 \begin{proof} By definition $\Ga((X_\sig)_-/K^*,(X_\sig)_+/K^*)$ is a toric variety. By the universal property of the graph (dominating component of the fiber product)  it corresponds to the coarsest simultaneous subdivision of  both  $\pi(\sigma_-)$ and  $\pi(\sigma_+)$, that is, to the fan $\{\tau_1\cap\tau_2\mid\tau_1\in \pi(\sigma_-), \tau_2\in\pi(\sigma_+)\}=\vr\cdot\pi(\sigma_-)=\vr\cdot\pi(\sigma_+)$.
\end{proof}

\subsection{The $\pi$-desingularization lemma of Morelli and centers of blow-ups}
For any simplicial cone
$\sigma =\langle v_1,\ldots,v_k\rangle $ in $N$ set $${\rm par}(\sigma ):=\{ v\in
\sigma\cap N_{\sigma}\mid  v=\alpha_1v_1+\cdots+\alpha_kv_k,
\mbox{where}\,\, 0\leq\alpha_i< 1\},$$  
$$\overline{{\rm par}(\sigma )}:=\{ v\in
\sigma\cap N_{\sigma}\mid  v=\alpha_1v_1+\cdots+\alpha_kv_k,
\mbox{where}\,\, 0\leq\alpha_i\leq 1\}.$$

 We associate with a dependent cone $\sig$ and an  integral vector
  $ v \in
\pi(\sigma)$ a vector ${\rm Mid}(v,\sigma):=\pi_{|\partial_-(\sigma)}^{-1}(v)+\pi_{|\partial_+(\sigma)}^{-1}(v)\in \sigma$ (\cite{Morelli1}), where
$\pi_{|\partial_-(\sigma)}$ and
$\pi_{|\partial_+(\sigma)}$ are the
restrictions of $\pi$ to $\partial_-(\sigma)$ and
$\partial_+(\sigma)$.

We also set ${\rm Ctr}_-(\sig):=\sum_{r_i<0}w_i$,  ${\rm Ctr}_+(\sig):=\sum_{r_i>0}w_i$.
 
\begin{lemma}(Morelli \cite{Morelli1}, \cite{Morelli2},
\cite{Abramovich-Matsuki-Rashid}) \label{le: Morelli} 
 Let $\Sigma$ be a simplicial cobordism in $N^+$. Then there exists a
simplicial cobordism $\Delta$ obtained from $\Sigma$ by a sequence of star
subdivisions such that $\Delta$ is $\pi$-nonsingular.
Moreover, the sequence can be taken so that any independent
and already $\pi$-nonsingular face of $\Sigma$ remains
unaffected during the process. 
All  the centers of the star
%
%
 subdivisions are of
the form $\pi_{|\tau}^{-1}({\rm par}(\pi(\tau)))$
%
%
 where $\tau$ is independent, and ${\rm Mid}({\rm Ctr}_\pm(\sig),\sig)$,  where $\sig$ is dependent. 
\end{lemma}

\begin{remark} It follows from Lemma \ref{le: Morelli} that $\pi$-desingularization  can be done for an open affine neighborhood  of a point $x$ of $F\in \cC(B^{K^*})$ on the smooth cobordism $B$ which is \'etale isomorphic with the tangent space $\Tan_{x,B}$.
We need to show how to globalize this procedure in a coherent and possibly  canonical way. This will replace the tangent space $\Tan_{x,B}$ in the local description of flips defined by elementary cobordisms  (as in Proposition \ref{local}) with $\pi$-nonsingular $X_\sig$. 

By Corollary \ref{fact} we get a factorization into a blow-up and a blow-down at smooth centers: $(B_a)_-/K^* \buildrel \phi_-\over\longleftarrow \Ga((B_a)_-/K^*,(B_a)_+/K^*)\buildrel \phi_+\over\longrightarrow (B_a)_+/K^*$. 
\end{remark}
\section{Proof of the $\pi$-desingularization  lemma}
\subsection{Dependence relation revisited}
\begin{lemma} \label{le: normal} (\cite{Wlodarczyk1}, Lemma
10, \cite{Morelli1}) Let $w_1,\ldots,w_{k+1}$ be integral vectors in
${\bf Z}^{k}\subset {\bf Q}^{k}$ 
which are not contained in a proper vector subspace of ${\bf Q}^{k}$.
Then $$\sum_{i=1}^{k+1}(-1)^i\det(w_1,\ldots,\check{w_i},\ldots,w_{k+1})\cdot w_i=0 $$ is the
unique (up to proportionality) linear relation between $w_1,\ldots,w_{k+1}$. 
\end{lemma}
\noindent{\bf Proof.}
 Let $v:=
\sum_{i=1}^{k+1} (-1)^i\det(w_1,\ldots,\check{w_i},\ldots,w_{k+1})\cdot w_i$.
Then  for any $i<j$, 
\[\begin{array}{rc}
&\det(w_1,\ldots,\check{w_i},\ldots,\check{w_j},\ldots,w_{k+1},v)=
\det(w_1,\ldots,\check{w_i},\ldots,w_{k+1})\cdot 
\det(w_1,\ldots,\check{w_i},\ldots,\check{w_j},\ldots,w_{k+1},w_i)+
\\
&\det(w_1,\ldots,\check{w_j},\ldots,w_{k+1})\cdot 
\det(w_1,\ldots,\check{w_i},\ldots,\check{w_j},\ldots,w_{k+1},w_j)=\\
&(-1)^i(-1)^{k-i}\det(w_1,\ldots,\check{w_j},\ldots,w_{k+1})
\cdot\det(w_1,\ldots,\check{w_i},\ldots,w_{k+1})+\\
&(-1)^j(-1)^{k-j+1}\det(w_1,\ldots,\check{w_i},\ldots,w_{k+1})
\cdot\det(w_1,\ldots,\check{w_j},\ldots,w_{k+1})=0.
\end{array}\]
Therefore $v\in \bigcap_{i,j}
\lin\{w_1,\ldots,\check{w_i},\ldots,\check{w_j},\ldots,w_{k+1}\}=\{0\}$.

\begin{definition} Let  $\delta=\langle v_1,\ldots,v_k\rangle$ be a
dependent cone and $w_i:=\prim(\pi(\langle v_i \rangle))$. 
Then we shall call a relation  $\sum_{i=1}^k r_iw_i=0$ {\it
normal} if it is a positive multiple of the relation (**) from Section \ref{dep} and
$|r_i|=|\det(w_1,\ldots,\check{w_i},\ldots,w_{k})|$ for $i=1,\ldots,k$. 

\subsection{Dependent $n$-cones}

Assume for simplicity that the normal relation is of the form
$$r_1w_1+\ldots +r_kw_k+r_{k+1}w_{k+1}+\ldots +r_{k+l}w_{k+l}+0\cdot w_{k+l+1}+0\cdot w_{k+l+2}+\ldots=0 \eqno(0),$$
where $r_1\geq r_2\geq\ldots r_k>0$ and $-r_{k+1}\geq -r_{k+2}\geq\ldots -r_{k+l}>0$.
We can represent it by two decreasing sequences of positive numbers:
$$r(\sig)=(r_1,r_2,\ldots,r_k;-r_{k+1},-r_{k+2},\ldots,-r_{k+l})$$
Set $\sgn(\sig)=+$ if  either $r_1>-r_{k+1}$ or $r_1=-r_{k+1}$ and $l\geq 2$. 

$\sgn(\sig)=-$ if  either $r_1<-r_{k+1}$ or $r_1=-r_{k+1}$ and $l=1$
\end{definition}

\begin{definition} An independent cone $\sigma$ is called an {\it $n$-cone}
if $|\det(\sigma)|=n$. A dependent cone $\sigma$ 
is called an {\it $n$-cone} if one of its independent faces 
is an $n$-cone and the others are $m$-cones, where $m\leq n$. 

We shall distinguish 5 types of dependent $n$-cones.
\begin{enumerate}
\item $(n,*;n,*)$, if $r_1=-r_{k+1}=n$ and $k,l\geq 2$.
\item $(n; n,*)$ , if $r_1=-r_{k+1}=n$ and either $k=1$ and $l\geq 2$ ($\sgn(\sig)=+$), or by symmetry $k\geq 2$ and $k=1$  ($\sgn(\sig)=-$). 
\item $(n,*;*)$ if $r_1=n>-r_{k+1}$ and $k\geq 2$  ($\sgn(\sig)=+$) or by symmetry $r_1<-r_{k+1}=n$ and $l\geq 2$ ($\sgn(\sig)=-$)
\item $(n;n)$ if  $r_1=-r_{k+1}=n$ and $k=l=1$ 
\item  $(n;*)$  if  $r_1=n>-r_{k+1}$, $k=1$, $l\geq 2$ ($\sgn(\sig)=+$), or by symmetry $r_1<-r_{k+1}=n$, $k\geq1$ and $l=1$ ($\sgn(\sig)=+$).
\end{enumerate}
\end{definition}
We can assign invariant to any dependent $n$-cones of type $(i)$, where $i=1,\ldots,5$, to be  $$\inv(\sig):=(n,-i).$$
These invariants are ordered lexicogrpahically.

\subsection{Star subdivision at $\Mid(\Ctr_{\sgn(\sig)}(\sig),\sig)$}
\begin{lemma} Let $\del$ be  a maximal dependent cone with a normal relation $(0)$. Then $v=\Mid(\Ctr_+(\delta),\delta)$ is in the relative interior of the circuit $\del_0:=
\langle v_i\mid r_i\neq 0\rangle$. The star subdivision at $\la v\ra$ affects the cones  $\del'\in\Star(\del_0,\Sig)$ only. All the normal relations for the cones $\del'$ are proportional to the normal relation for $\del$.
\end{lemma}
\begin{proof}
$w_1+\ldots+w_k=w_1+\ldots+w_k-\epsilon(r_1w_1+\ldots+r_kw_k+r_{k+1}w_{k+1}+\ldots+r_{k+l}w_{k+l})$ is a combination of $w_1,\ldots,w_{k+l}$ with positive coefficients for a sufficiently small $\epsilon>0$ .
\end{proof}
\begin{lemma} \label{le: normal2} Let $\delta=\langle v_1,\ldots,v_{k+l},\ldots,v_r\rangle$ be a maximal dependent cone with a normal relation $(0)$.
 Let $v=\Mid(\Ctr_+(\delta),\delta)\in\inte
\langle v_1,\ldots,v_{k+l}\rangle$. Let $m_w\geq 1$ be an integer
such that the vector
$$w=\frac{1}{m_w}(w_1+\ldots+w_k)$$ is primitive. Then the maximal
dependent cones in $\langle v \rangle\cdot\delta$ are of the
form $\delta_{i_0}=\langle
v_1,\ldots,\check{v}_{i_0},\ldots,v_{k+l},\ldots, v_r,v\rangle$, where $i_0\leq k+l$.
\begin{enumerate}

\item  Let $r_{i_0}>0$ i.e $1\leq i_0\leq k$. Then for the maximal dependent cone
$\delta_{i_0}=\langle
v_1,\ldots,\check{v_{i_0}},\ldots,v_k,v\rangle$ in 
$\langle v \rangle\cdot\delta$, the  normal relation
is given (up to sign) by $$\sum_{r_i> 0, i\neq i_0}\frac{r_i-r_{i_0}}{m_w}w_i+
\sum_{r_i< 0} \frac{r_i}{m_w}w_i+r_{i_0}w=0. \eqno(1a)$$

\item Let $r_{i_0}<0$ i.e $k+1\leq i_0\leq l+k$. Then for the maximal dependent cone
$\delta_{i_0}=\langle
v_1,\ldots,\check{v}_{i_0},\ldots,v_k,v\rangle$ in $\langle
v \rangle\cdot\delta$, the  normal relation
is given (up to sign) by $$\sum_{i\neq i_0}-\frac{r_{i_0}}{m_w}w_i+
r_{i_0}w=0. \eqno(1b)$$
\end{enumerate}
\end{lemma}
\begin{proof} It is straightforward to see that the above equalities hold. We only
need to show that the relations considered are normal. For that
it suffices to show that one of the coefficients is equal
(up to sign) to the corresponding determinant. 

 Comparing the coefficients of $w$ in the above
relations with the normal
relations from Lemma \ref{le: normal} we get

1a. $\det(w_1,\ldots,\check{w}_{i_0},\ldots,w_k)=
r_{i_0}$.

1b. The coefficient of $w$ is equal to
$\det(w_1,\ldots,\check{w}_{i_0},\ldots,w_k)=r_{i_0}$. 
\end{proof}

\begin{corollary} Let $\del$ be a dependent $n$-cone and $v=\Mid(\Ctr_{\sgn(\del)}(\del),\del)$. 
\begin{enumerate}
\item If $\del$ is of type $(1)$ or $(3)$ then $\la v\ra\cdot\del$ consists of $n$-cones of smaller type. 
\item If $\del$ is of type $(2)$ or $(5)$ then $\la v\ra\cdot\del$ consists of $n$-cones of one cone of the same type as $\del$ and cones of smaller type. 
\end{enumerate}
\end{corollary}
\begin{proof} Without loss of generality assume that $\sgn(\del)=+$. Then $r_1=n$ . 

(1) In the case when $\del$ is of type (1) or $(3)$, the index $k\geq 2$.

If $r_{i_0}=n$ then in the relation $(1b)$, $r_i-r_{i_0}\leq 0$ and only the coefficient $r_{i_0}=n$ is positive.  The cone $\del_{i_0}$ is 
$n$-cone  $(n;n,*)$ (in case $\del$ is of type (1)) and $(n;*)$ in case $\del$ is of  type (3) corresponding to $r_{i_0}=n$.

If $n>r_{i_0}>0$ then  in the relation $(1b)$, all positive coefficients $r_i-r_{i_0}\leq 0$ and  $r_{i_0}$ are smaller than $n$. The cone $\del_{i_0}$ is an
$n$-cone  $(n;*)$ (in case $\del$ is of type (1)) and $(*;*)$ in case $\del$ is of  type (3) corresponding to $r_{i_0}=n$.

If $r_{i_0}=-n$ then $\del$ is of type $(1)$ and $\del_{i_0}$ is an
$n$-cone  $(n;n)$ of type (4).

If $-n<r_{i_0}<0$ then $\del_{i_0}$ is an
$m$-cone  $(m;m)$, $m<n$. 

(2) In the case when $\del$ is of type (2) or $(5)$, the index $k=1$ and we have one positive ray only $r_1=n$.

If $r_{i_0}=r_1=n$ then in the relation $(1b)$, $r_i-r_{i_0}\leq 0$ and only the coefficient $r_{i_0}=n$ is positive.  The cone $\del_{i_0}$ is 
$n$-cone  $(n;n,*)$ (in case $\del$ is of type (2)) and $(n;*)$ in case $\del$ is of  type (5).

If $r_{i_0}=-n$ then $\del$ is of type $(2)$ and $\del_{i_0}$ is 
$n$-cone  $(n;n)$.

If $-n<r_{i_0}<0$ then $\del_{i_0}$ is 
$m$-cone  $(m;m)$, $m<n$. 
\end{proof}
A direct consequence of the above is the following
\begin{corollary} \label{def} Let $n=\max\{|\det(\pi(\tau)|\mid \tau\quad  \mbox{is independent in} \quad \Sig\}$.
Let $\del$ be a dependent $n$-cone with a circuit $\del_0$ and $v=\Mid(\Ctr_{\sgn(\del)}(\del),\del)$. If $\del$ is of type 
$(i)$ then all the dependent cones  $\del'\in\Star(\del_0,\Sig)$ are $m$-cones, with $m\leq n$ of type $(i)$. Moreover:
\begin{enumerate}
\item If $\del$ is of type $(1)$ or $(3)$ then either $\la v\ra\cdot\Sig$ contains
a smaller number  of dependent cones with maximal invariant $\inv(\sig)$ or 
the maximal invariant $\inv(\sig)$ drops. 
\item If $\del$ is of type $(2)$ or $(5)$ then $\la v\ra\cdot\Sig$  contains 
unchanged number of dependent cones with maximal invariant $\inv$.
\end{enumerate}
\end{corollary}

\subsection{Codefinite faces}
\begin{definition} \label{def0}(\cite{Morelli1}) An independent  face $\tau$ of a dependent cone $\del$ is called {\it codefinite} iff it does not contain both negative and positive rays.
\end{definition}
Every dependent cone $\del$ contains two  maximal codefinite  faces $$\del^+:=\la v_i\mid r_i\geq 0\ra \quad\quad \del^-:=\la v_i\mid r_i\leq 0 \ra.$$
\begin{corollary} Any   independent cone $\tau\in\Sig$ can be made a codefinite face of all dependent cones containing it. The process use star  subdivisions at $\Mid(\Ctr_{\sgn(\sig)}(\sig),\sig)$ applied to dependent cones for which $\tau$ is not codefinite. Moreover using the procedure we do not increase a number of cones with maximal invariant.
\end{corollary}
\begin{proof}
First apply the procedure to all dependent cones of types (1)  and (3) for which $\tau$ is not codefinite.This procedure terminates since by the previous lemma the invariant drops until we arrive at the situation where all cones for which $\tau$ is not codefinite are of type
(2) or  (5). Note also that cones of type (3) have only one positive and one negative ray and thus all their independent faces are codefinite.
Next apply the star subdivision  at $v=\Mid(\Ctr_{\sgn(\sig)}(\sig),\sig)$ to all cones of type (2) or (5).
After the star subdivision the cone $\del=\la v_1,\ldots, v_k\ra $ with the normal relation
$$r_1w_1+r_{k+1}w_{k+1}+\ldots +r_{k+l}w_{k+l}=0,$$ where $k=1$, we create   a cone
$\del_1=\la v',\ldots, v_k\ra $ with exactly the same relation
$$r_1w'+r_2w_2+\ldots +r_{k+l}w_{k+l}=0$$
Since $v_1$ was only negative ray an it was replaced with the center of subdivision $v$, the cone $\tau$ contains only negative rays of $\del_1$. Other cones $\del_j$, where $j\geq 2$ are of the type (3) and again $\tau$ is not their  codefinite  face.
\end{proof}

\subsection{Star subdivisions at $v\in\pi_{|\tau}^{-1}(\para(\pi(\tau)))$}

\begin{lemma}
 Let
$w=\sum_{i\in I}\alpha_iw_i\in
\para(\pi(\del^+))$. Let $\tau=\langle v_i\mid i\in I\rangle\preceq \del^+$ be a
codefinite face of
$\delta$ containing $v=\pi_{|\del^+}^{-1}(w)$ in its relative interior. Then the maximal dependent cones in 
 $\langle v \rangle\cdot\delta$
are of the form $\delta_{i_0}=
\langle v_1,\ldots,\check{v}_{i_0},\ldots,v_k,v\rangle$,
where $i_0\in I$.

2a. Let $i_0\in I$ and $r_{i_0}>0$. Then for the maximal
dependent cone $\delta_{i_0}=
\langle v_1,\ldots,\check{v}_{i_0},\ldots,v_k,v\rangle$ in
$\langle v  \rangle\cdot\delta$, 
 the  normal relation
is given (up to sign) by 
$$\sum_{i\in I\setminus \{i_0\},r_i>0}
(\alpha_{i_0}r_{i}-
\alpha_{i}r_{i_0})w_{i}+\sum_{i\not\in I, r_i>0}
\alpha_{i_0}r_{i}w_i+\sum_{
i\in I,r_i=0} -\alpha_{i}r_{i_0}w_{i}+ 
\sum_{r_i<0}\alpha_{i_0} r_iw_{i}+r_{i_0}w=0. \eqno(2a)$$ 

2b. Let  $i_0\in I$ and $r_{i_0}=0$. For  the maximal dependent cone
$\delta_{i_0}=\langle
v_1,\ldots,\check{v}_{i_0},\ldots,v_k,v\rangle$,
  the  normal relation
is given (up to sign) by
$$\sum_{i\neq i_0} \,\,\alpha_{i_0}r_iw_i+0w=0. \eqno(2b)$$

\end{lemma}
\begin{proof}
2a. The coefficient of $w$ is equal to 
$\det(w_1,\ldots,\check{w}_{i_0},\ldots,w_k)=r_{i_0}$. 

2b. The coefficient of $w_i$, where $r_i>0$, is equal to \\ 
$\det(w_1,\ldots,
\check{w}_{i_0},\ldots,\check{w_{i}},\cdots,w_k,w)=\alpha_{i_0}
\det(w_1,\ldots,
\check{w}_{i_0},\ldots,\check{w_{i}},\cdots,w_k,w_{i_0})=\\
(-1)^{k-i_0}\alpha_{i_0}
\det(w_1,\ldots,\check{w}_{i},\cdots,w_k,w)=
(-1)^{k-i_0} \alpha_{i_0}r_i$.
\end{proof}
\begin{lemma}\label{def1} Let  $\del$ be  a dependent $n$-cone
$w\in \para(\pi(\del^+))$ and $v=\pi_{|\del^+}^{-1}(w)\in\del^{\pm}$.
Then
\begin{enumerate}
\item If $\del$ is of type (1) or (2) or (4) then $\la v\ra\cdot\del$ contains  $n$-cones of smaller type. 
\item If $\del$ is of type (3) then $\la v\ra\cdot\del$ may contain  $n$-cones of type (3) and smaller type. 
\item  If $\del$ is of type (5) then $\la v\ra\cdot\del$ may contain at most one $n$-cone of type (5) and smaller types.

\end{enumerate}
\end{lemma}
\begin{proof}
(1)  If $\del$ is of type (1) or (2) that is $(n,*:n,*)$ or $(n:n,*)$. After the star subdivision we create $n$-cones for $r_{i_0}=n$ with only one coefficient $n$. These are $n$-cones of type $(3)$ or $5$.
If we subdivide cone of type $(4)$ that is $(n;n)$ we create only one $n$-cone of type $(5)$.

(2) If $\del$ is of type (3) that is $(n,*:*)$. After the star subdivision we create $n$-cones for $r_{i_0}=n$ with only one coefficient $n$. These can be $n$-cones of type $(3)$ or $(5)$.

(3) If If $\del$ is of type (5) that is $(n:*)$. After the star subdivision we create  an $n$-cone  of type (5) for $r_{i_0}=n$. It has only one positive coefficient $n$ and other coefficients are negative $>-n$. \end{proof}

\begin{lemma}\label{def2} Let $n>1$ and  $\del$ be  a dependent $n$ cone of type $(2)$, $(4)$ or $(5)$.  Then $\del^{-\sgn(\del)}$ is a maximal independent face and $|\det(\pi(\del^{-\sgn(\del)})|=n.$ There exists
$w\in
\para(\pi(\del^{-\sgn(\del)}))$ and the corresponding $v=\pi_{|\del^{-\sgn(\del)}}^{-1}(w)\in\del^{-\sgn(\del)}$.
The subdivision $\la v\ra\cdot\del$ contains  $n$-cones of smaller type. 
\end{lemma}
\begin{proof} Without loss of generality we assume that $\sgn(\del)$ is negative and we take a star subdivision at
$v\in\del^+$. 
If $\del$ is of type $(2)$ then $r_1=-r_{k+1}=n$, so we have one negative ray with coefficient $-n$ and $k\geq 2$ positive rays with coefficients $\leq n$. 
If $\del$ is of type $(4)$ then  we have one positive and  one negative ray with coefficients $n$ and $-n$.
If $\del$ is of type $(5)$ then $r_1<-r_{k+1}=n$, so we have one negative ray with coefficient $-n$ and $k\geq 2$ positive rays.
After the subdivision at $\la v\ra$ we create cones with negative coefiicients $>-n$ and positive rays   with coefficients $\leq r_1$.

If  $\del$ is of type $(2)$ the new dependent $n$-cones are of type (3), (4), (5).

If  $\del$ is of type $(4)$ the new dependent $n$-cones are of type (5).
In the normal relation for a new cone $\del_{i_0}$, where $r_{i_0}=r_1=n $ there is one positive ray with coefficients $n$ and negative coefficients $>-n$.

If  $\del$ is of type $(5)$ we create only dependent $m$-cones with $m<n$.
\end{proof}

\subsection{Resolution algorithm}
The $\pi$-desingularization algorithm consists of
eliminating all dependent $n$-cones, where $n>1$ in the following order. 

{\bf Step 1}. {\it Eliminating all dependent $n$-cones $\del$ of type $(1)$}  by applying the
star subdivision at 
$\la \Mid(\Ctr_{\sgn(\del)}(\del),\del)\ra$. (Corollary
\ref{def}.)

\bigskip
{\bf Step 2}. {\it Eliminating all dependent  $n$-cones $\del$ of type $(2)$} .

{\bf Step 2a}. By definition $\del^{-\sgn(\del)}$ is a maximal independent face and $|\det(\pi(\del^{-\sgn(\del)})|=n.$

Let $v\in\pi^{-1}_{|\del^{-\sgn(\del)}}(\para(\pi(\del^{-\sgn(\del)}))$. Then $v\in\inte(\tau)$ for some independent face $\tau\preceq \del^{-\sgn(\del)}$. We
make  $\tau$ codefinite with respect to all
dependent cones containing it. By Lemma \ref{def0}, this process will not increase a number of $n$-cones of type (2).

{\bf Step 2b}. Apply the star subdivision at $\la v\ra$. We change all the cones in $\Star(\tau,\Sig)$. The cone $\tau$ is codefinite with respect to all faces from $\Star(\tau,\Sig)$. Moreover
by definition $\tau\preceq \del^{-\sgn(\del)}$.
By Lemmas \ref{def1}, \ref{def2},
the process will decrease  the number of $n$-cones of type (2).

\bigskip
{\bf Step 3}.  {\it Eliminating all dependent  $n$-cones $\del$ of type $(3)$}  by applying
star subdivision at $\la \Mid(\Ctr_{\sgn(\del)}(\del),\del)\ra$.

\bigskip
{\bf Step 4}. {\it Eliminating all  dependent $n$-cones of type  (4) by using the two steps procedure as in Step 2.}

\bigskip
{\bf Step 5}. {\it Eliminating all dependent  $n$-cones of type  (5) by using the two steps procedure as in Step 2.}

\bigskip
{\bf Step 6}. {\it Eliminating all independent $n$-cones $\tau$ which are not faces of some dependent cones}.

{\bf Step 6a}. 
Let $v\in\pi^{-1}_{\tau}(\para(\pi(\tau))$. Then $v\in\inte(\tau_0)$ for some independent face $\tau_0\preceq \tau$. We
make  $\tau_0$ codefinite with respect to all
dependent cones containing it.

{\bf Step 6b}. Apply the star subdivision at $\la v\ra$. Determinant of all independent faces $\tau'$ containing $\tau_0$ drops.
$$|\det\pi(\tau')|=|\det(w_1\ldots\check{w_i},\ldots,w_k,w)|=\alpha_i |\det(w_1\ldots,w_k)|<
|\det(w_1\ldots,w_k),|$$
where $\pi(\tau')=\la w_1,\ldots,w_k\ra$,  $w=\sum_i\alpha_iw_i$, $0\leq\alpha_i<1$, $v=\pi_{|\tau_0}^{-1}(w)$.
\begin{remark} 
 The strategy of this algorithm of using the above centers for the corresponding $n$-cones was first applied in \cite{Wlodarczyk1} in the proof of regularization of toric factorization (see  (\cite{Wlodarczyk1}), Lemmas 11-12 pages 403-410). Then it was used directly in the context of $\pi$-desingularization in \cite{Abramovich-Matsuki-Rashid} and in the revision of the Morelli's original algorithm in \cite{Morelli2}.
\end{remark}
\subsection{The Weak Factorization of toric morphisms}\label{toric}
\begin{theorem}(\cite{Wlodarczyk1},\cite{Morelli1})  Lef $f:X\dashrightarrow Y$ be a 
birational toric map  of
smooth complete toric varieties. 
Then $f$ can be factored
 as
$X=X_0\buildrel f_0 \over \dashrightarrow  X_1
\buildrel f_1 \over \dashrightarrow \ldots \buildrel f_{n-1} \over
\dashrightarrow X_n=Y  ,$
where each $X_i$ is a smooth  complete toric variety and $f_i$ is a  blow-up
or blow-down at a smooth invariant center.
 \end{theorem}
 \begin{proof} By Proposition \ref{fact}, there is a a smooth toric variety $Z$ and a factorization of $f$ into $X\leftarrow Z\to Y$ where $Z\to X$ and $Z\to Y$ are projective  toric morphisms. By Proposition \ref{construction} there is a toric variety $\overline{B}\supset B\supset T\times K^*$ which is a compactified cobordism defined for a projective  toric morphism $Z\to X$. The variety $B$ corresponds to a strictly $\pi$-convex nonsingular fan $\Del$. Its $\pi$-desingularization determines a $\pi$-nonsingular fan $\Del^\pi$ corresponding to a toric variety $B^\pi$ projective over $B$. The open subsets $B_- $ and $B_+$ have smooth quotients  $B_- /K^*$ and $B_+/K^*$. They correspond to $\pi$-nonsingular subfans $\Del_+$ and $\Del_-$ of $\Del$ and are not affected by $\pi$-desingularization.
 That is $B_-=B^\pi_-$ and $B_+=B^\pi_+$ and $B^\pi$ is a cobordism  between $X$ and $Z$ and admits a compactification $\bar{B}^\pi=B^\pi\cup X\cup Z=B^\pi\cup \cO(Z)\cup X\times (\bP^1\setminus\{0\})$ (see  Proposition \ref{construction}). By Proposition \ref{deco}, the  compactified cobordism
 $\bar{B}^\pi\supset B^\pi$ determines a 
 a decomposition into elementary cobordisms $B^\pi_a$ and the toric factorization into maps $(B^\pi_a)_-/K^*\dashrightarrow (B^\pi_a)_+/K^*$.
 If $B^\pi_a$ is an elementary cobordism corresponding to the fan $\Del_a$ then $(B^\pi_a)_-$, $(B^\pi_a)_+$ and $(B^\pi_a)_-\cap (B^\pi_a)_+\subset B^\pi_a$ correspon to subfans $\Del_-$, $\Del_+$ and $\Del_0:=\Del_+\cap\Del_-$  respectively, consisting of independent cones. Every toric orbit in $B^\pi_a\setminus ((B^\pi_a)_-\cap (B^\pi_a)_+)=F^+\cup F^-$ contains a fixed point orbit in its closure corresponding to a dependent cone. Thus $\Del_a\setminus \Del_0$ is  a collection of dependent cones and some of their independent faces. 
If $\sig\in \Del_a$ is a dependent cone then $X_\sig$ intersects a unique fixed
point component $\bar{O}_\del$, where $\del$ is a unique circuit in $\sig$. All orbits in $X_\sig$ contain in its closure a fixed orbit $O_\tau\subset \bar{O}_\del$. Thus $X_\sig$ is disjoint from other closed sets
$F^+$ and $F^-$, where $F\neq \bar{O}_\del\in\cC(B_a^\pi)^{\ks}$. We get that $(X_\sig)_-=X_\sig\setminus (\bar{O}_\del)^+=X_{(\partial_+(\sig)}\subset B^\pi_-$ and $(X_\sig)_+=X_{(\partial_-(\sig)} \subset B^\pi_+$. It follows from the above that 
 $\pi(\Del_+)$ and $\pi(\Del_-)$ are two nonsingular subdivisions of the fan $\pi(\Del_a)$ which coincide on $\pi(\Del_0)$ and which define two different decomopsitions for all projections $\pi(\sig)$ of dependent cones: $\pi(\Del_+)_{|\pi(\sig)}=\pi(\partial_-(\sig))$ and  $\pi(\Del_-)_{|\pi(\sig)}=\pi(\partial_+(\sig))$. If $\del$ is a circuit in $\sig$, such that $\pi(\del)=\la w_1,\ldots,w_k\ra$ then  by Lemma \ref{mo}, the unique relation is given by $\sum r_iw_i=0$ where $r_i=\pm 1$. Let $w_\del=\sum_{r_i=1} w_i=-\sum_{r_i=-1} w_i$. Then
 the ray $\la w_\del\ra$ determines nonsingular star subdivisions of 
$\pi(\Del_+)$, $\pi(\Del_+)$ and
$\la w_\del\ra\cdot\pi(\Del_+)_{|\pi(\sig)}=\la v_\del\ra\cdot\pi(\Del_-)_{|\pi(\sig)}$. If $\del_1\ldots,\del_r$ be all ciruits in $\Del_a$ then the stars $\Star(\pi(\del_i),\pi(\Del_a))$ are disjoint and $\la w_{\del_1}\ra \cdot\dots\cdot \la w_{\del_r}\ra\cdot \pi(\Del_+)=\la w_{\del_1}\ra \cdot\dots\cdot \la w_{\del_r}\ra\cdot\pi(\Del_-)$ and consequently
$(B^\pi_a)_-/K^*\dashrightarrow (B^\pi_a)_+/K^*$ factors into a sequence of blow-ups about smooth toric centers followed by  a sequence of blow-downs about the smooth centers.
\end{proof}
\section{$\pi$-desingularization of birational cobordisms}
\subsection{Stratification by isotropy groups on a smooth cobordism}
\label{stratification}
\medskip\newcommand{\os}{\overline s}
\noindent
Let $B$ be a smooth cobordism of dimension $n$. Denote by $\Ga_x$ the isotropy group of a point $x\in B$. Let $D$ be a $\ks$-invariant divisor on $B$ with simple normal crossings. 

Define the stratum $s=s_x$ through $x$ to be an irreducible component of the set  $\{p\in B\mid\Ga_x=\Ga_p\}.$
 We can find  $\Gamma_x$-semiinvariant parameters in the affine open neighborhood $U$ of $x$ such that
\begin{enumerate}
\item$\Ga_x$ acts nontrivially on
$u_1,\dots, u_k$ and trivially on $u_{k+1},\dots,u_n$.
\item Any  component of $D$ through $x$ is described by a parameter $u_i$ for some $i$.
\end{enumerate}
After suitable shrinking of $U$ the parameters define an \'etale $\Ga_x$-equivariant
morphism $\varp: U\to \Tan_{x,B}=\bbA^n$. By definition the stratum $s$ is locally described by
$u_1=\ldots =u_k=0$. The parameters $u_1, \dots,u_k$ determine a $\Ga_x$-equivariant
smooth morphism $$\psi: U \to \Tan_{x,B}/{\Tan_{x,s}}=\bbA^k.$$  We shall view
$\bbA^k=X_\sig$ as a toric variety with a torus $T_\sig$ and refer to $\psi$ as
a {\it toric chart}.
This assigns to a stratum $s$ the cone $\sigma$ and the relevant group $\Gamma_\sigma$ acting on $X_\sigma$.
Then Luna's \cite{Luna} fundamental lemma implies that the morphisms  $\phi$ and $\psi$ 
preserve stabilizers, the induced morphism $\psi_\Gamma:U//\Gamma_x\to X_\sigma//\Ga_\sig$ is smooth and $U\simeq U//{\Ga_x}\times_{\bbA^k//\Ga_x}\bbA^k$. Note that for a toric charts on $B$ we require that inverse images of toric divisors have simple normal crossings with components of $D$. We refer to this property as {\it compatibility with D}.

The invariant  $\Ga_x$ can be defined for $X_\sig=\bbA^k$ and determine
the relevant $T_\sig$-invariant stratification $S_\sigma$ on $X_\sig$.
By shrinking $U$ we may assume that the strata on $U$ are inverse images of the
strata on $X_\sig$. Any stratum $s_y$ on $U$ through $y$ after a suitable
rearrangement of $u_1,\ldots, u_k$ is described in the neighborhood $U'\subset U$ of $y$ by $u_1=\ldots =u_\ell=0$,
 where  $\Ga_y\leq \Ga_x$ acts nontrivially on
$u_1,\dots,u_\ell$ and trivially on $u_{\ell+1}\ldots,u_k,u_{k+1},\ldots,u_n$. The remaining $\Ga_y$-invariant parameters  at $y$ are
$u_{\ell+1}-u_{\ell+1}(y),\dots,u_{n}- u_{n}(y)$. Then the closure of
$\os_y$
%
%
 is described on $U$ by $u_1=\ldots =u_\ell=0$ and contains $s_x$. This shows \begin{lemma} The
closure of any stratum is a union of strata. 
\end{lemma}
 We can introduce an order on the strata by setting
$$s'\leq s\quad\quad \mbox{iff} \quad\quad\os'\subseteq s.$$

\begin{lemma} \label{inclusion} If $s'\leq s$ then there exists an inclusion $i_{\sig'\sig}:\sig'\hookrightarrow\sig$ onto a face of $\sig$. 
The inclusion $i_{\sig'\sig}$ defines a $\Ga_{\sig'}$-equivariant morphism of toric varieties
$X_{\sig'}\to X_{\sig'}\times 1\hookrightarrow X_{\sig'}\times T\subset X_{\sig},$
where $T_{\sig'}\times T= T_{\sig}$ and  $\Ga_{\sig'}\subset T_{\sig'}$. 
Moreover we
can write $X_\sig\cong X_{\sig'}\times \bbA^r$ where $\Ga_{\sig'}$ acts trivially on
$\bbA^r$ and nontrivially on all coordinates of $X_{\sig'}\simeq\bbA^\ell$.
\end{lemma}

\noindent
 In the above situation  we shall write
$$ \sig'\leq \sig.$$ 
The lemma above implies immediately
\begin{lemma} \label{maximal} If $\tau<\sig$ (that is, $\tau\leq\sig$, $\tau\neq\sig$)  then $\Ga_\tau\subsetneq \Ga_\sig$.
\end{lemma}
Consider the stratification $S_\sig$ on $X_\sig$. Every stratum $s_\tau\in S_\sig$, where $\tau\leq \sig$,  is a union of orbits $O_{\tau'}$. Set $$\bar\tau:=\{\tau'\mid O_{\tau'}\subset s_\tau\}.$$ 
\begin{lemma} Any cone from the set $\tau'\in \overline{\tau}$ can be expressed as $\tau'\simeq \tau\times\la e_1,\ldots,e_r\ra\subset\sig$, and $X_{\tau'}=X_\tau\times \bbA^s\times T^{r-s}$ where $\Ga_\tau$ acts trivially on $\bbA^r\times T^{r-s}$.
\end{lemma}

\begin{lemma} For any $\tau'\in \overline{\tau}$, we have
$\Gamma_\tau=\Gamma_{\tau'}:=\{g\in\Ga_\sig\mid \underset{x\in O_{\sig'}}{\forall} gx=x\}$.
\end{lemma}

\subsection{Local projections}
\begin{definition} A cone $\sig$ in $N^{\QQ}$ is \textit{of maximal dimension} if $\dim \sig=\dim N^{\QQ}$.
\end{definition}
Every cone $\sig$ in $N^{\QQ}$ defines a cone of maximal dimension in $N^{\QQ}\cap \rm{span} \{\sig\}$ with lattice $N\cap\rm{span} \{\sig\}$. We denote it by $\underline\sig$. There is a noncanonical isomorphism
$$X_\sig=X_{\underline\sig}\times O_\sig.$$
The vector space $\mbox{span}\,\{\sig\}\subset N^{\QQ}$ corresponds to a subtorus $T_{\underline\sig}\subset T_\sig$ defined as $T_{\underline\sig}:=\{t\in T_\sig\mid tx=x\ \mbox{for}\ x\in O_\sig\}.$ Then $O_\sig$ is isomorphic to the torus $T_{\sig}/T_{\underline\sig}$ with dual lattice  $\sig^\perp\subset M^Q$.

\begin{lemma} \label{le: p}If $\Ga\subset T_\sig$ acts freely on $X_\sig=X_{\underline\sig}\times O_\sig$ then $$X_{\sig}{/\Ga}=X_{\underline\sig}\times O_\sig/\Ga,$$ where $O_\sig\simeq O_\sig/\Ga$ if $\Ga$ is finite, while $O_\sig/\Ga$ is isomorphic to a torus of dimension $\dim O_\sig-1$ if $\Ga=K^*$.
\end{lemma}

\begin{proof} By assumption $\Ga\cap T_{\underline\sig}$ is trivial. Hence $\Ga$ acts trivially on $X_{\underline\sig}$ and $X_\sig/\Ga=X_{\underline\sig}\times O_\sig/\Ga$. 
\end{proof}

Let $\pi_\sig: (\sigma,N_\sig)\to (\sig^{\Gamma},N_\sig^{\Gamma})$ denote the  projection corresponding to the quotient map $X_\sig\to X_{\sig}//{\Ga_\sig}$.
\begin{lemma}\label{pro}
If $\tau\leq \sig$ then $\pi_\tau(\tau)\simeq\pi_\sig(\tau).$
\end{lemma}

\begin{proof}$X_{\underline\tau}\times O_\tau$ is an open subvariety in $X_\sig$ and $\Ga_\tau$ acts trivially on $O_\tau$. We have
\[(X_{\underline{\tau}}\times O_\tau)/{\Ga_\tau}=X_{\underline{\tau}}/{\Ga_\tau}\times O_\tau=X_{\pi_\tau(\tau)}\times O_\tau.\]
$\Ga_{\sig}/{\Ga_\tau}$ acts freely on $(X_{\underline{\tau}}\times O_\tau)/{\Ga_\tau}=X_{\pi_\tau(\tau)}\times O_\tau$. Thus by the previous lemma 
$X_{\pi_\sig(\tau)}\cong X_{\pi_\tau(\tau)}\times O_{\tau}/{\Ga_\sig}.$
%
%
\end{proof}
For any $\tau\in\Del^\sig$, set $\Ga_\tau:=\{g\in\Ga_\sig\mid \underset{x\in O_{\tau}}{\forall} gx=x\}$.
Similarly one proves:
\begin{lemma}\label{pro2} Let $\Ga\subset \Ga_\sig$ be a  group
containing $\Ga_\tau$, where $\tau\in\Del^\sig$. Let $\pi_\Gamma: \sigma\to\sig^\Gamma$ be the projection corresponding to the quotient $X_\sig\to X_\sig/\Ga$. Then $\pi_\Ga(\tau)\simeq\pi_\sig(\tau).$
\end{lemma}

\begin{lemma} \label{pr} Let $\Gamma$  be a subgroup of $\Gamma_\sigma$, and $\pi_\Gamma: \sigma\to\sig^\Gamma$ be the projection corresponding to the quotient $X_\sig\to X_\sig/\Ga$. For any  $\tau\leq \sig$ and $\tau'\in\overline\tau$  we have
$\tau'=\tau\oplus \la e_1,\ldots,e_k\ra$ where $\la e_1,\ldots,e_k\ra$ is nonsingular and $\pi_\Ga(\tau')=\pi_\Ga(\tau)\oplus \la e_1,\ldots,e_k\ra$.
\end{lemma}

\begin{proof} $X_{\tau'}=X_\tau\times \bbA^k\times O_{\tau'}$ where the action of $\Ga_\tau\cap \Gamma$ on $\bbA^k\times O_{\tau}$ is trivial. Thus $X_{\tau'}/{\Ga_\tau}=X_{\tau}/{\Ga_\tau}\times \bbA^k\times O_{\tau'}$.  Now $\Ga/({\Ga_\tau}\cap \Ga)$ acts freely on $O_{\tau'}\subset s_\tau$ and we use Lemma \ref{le: p}.
\end{proof}

\subsection{Independent and dependent cones}

By Lemma \ref{pro} there exists a lattice isomorphism $j_{\tau\sig}:\pi_\tau(\tau)\to\pi_\sig(\tau)$, where $\tau\leq\sig$. Thus the projections $\pi_\tau$ and $\pi_{\sig}$ are coherent and related: $j_{\tau\sig}\pi_\tau=\pi_{\sig}$.

\noindent {\bf Case 1}: $\Ga_\sig=K^*$. The action of $\ks$ on $X_\sig$ corresponds to a primitive vector $v_\sig\in N_\sig$. The invariant characters $M^\Ga_\sig\subset M_\sig$ are precisely those $F\in M^\Ga_\sig$ such that $F(v_\sig)=0$. The dual morphism is a projection $\pi_\sig :N_\sig\to N_{\sig}/{\bZ\cdot v_\sig}=N^\Ga_\sig$.

\noindent
The quotient morphism of toric varieties $X_\sig\to X_{\sig}/\Ga_\sig$ corresponds to the projection $\sig\to \pi_\sig(\sig)$.

\bigskip\noindent
{\bf Case 2}: $\Ga_\sig\cong \ZZ_n$. The invariant characters $M^\Ga_\sig\subset M_\sig$ form a sublattice of  dimension $\dim(M^\Ga_\sig)=\dim(M_\sig)$, where $M_{\sig}/{M_\sig^\Ga}\simeq \ZZ_n$.
 The dual morphism defines an inclusion $\pi:N_\sig\hookrightarrow N^\Ga_\sig$. The projection $\sig\to \pi_\sig(\sig)$ is a linear isomorphism which does not preserve lattices. This gives 

\begin{lemma}

$X_\tau$ is independent iff $\Ga_\tau$ is finite. $X_\sig$ is dependent iff $\Ga_\sig=\ks$.
\end{lemma}

\begin{definition} Let $\Del^\sig$ be a decomposition of a cone $\sig\in\Sig$. A cone $\tau\in\Del^\sig$ is {\it independent} if $\pi_{\sig{|\tau}}$ is a linear isomorphism. A cone $\tau$ is {\it dependent} if $\pi_{\sig|\tau}$  is not a linear isomorphism.\end{definition}

\medskip
\subsection{Semicomplexes and birational modifications of cobordisms}

By glueing cones $\sig$ corresponding to strata along their faces we construct a {\it semicomplex} $\Sigma$, that is, a partially ordered set of cones  such that for $\sig\leq \sig'$ there exists a face inclusion $i_{\sig\sig'}:\sig\to\sig'$.

\begin{remark} The glueing need not be transitive: for $\sig\leq\sig'\leq\sig''$ we have $i_{\sig'\sig''}i_{\sig\sig'}\neq i_{\sig\sig''}$. Instead, there exists an automorphism $\alpha_\sig$ of $\sig$ such that $i_{\sig'\sig''}i_{\sig\sig'}= i_{\sig\sig''}\alpha_\sig$.
\end{remark}
For any fan $\Sigma$   denote by $\rm{Vert}(\Sig)$ the set of all 1-dimensional faces (rays) in $\Sig$. Denote by $\rm{Aut}(\sig)$ the automorphisms of $\sig$ inducing $\Ga_\sig$-equivariant automorphisms.

\medskip
\begin{definition}\label{de:sub} By a {\it subdivision} of $\Sig$ we mean a collection $\Del=\{\Del^\sig\mid \sig\in\Sig\}$ of subdivisions $\Del^\sig$ of $\sig$ such that 

\item{1$^\circ$} If $\tau\leq \sig$ then the restriction $\Del^\sig_{|\tau}$ of $\Del^\sig$ to $\tau$ is equal to $\Del^\tau$.
\item{2$^\circ$} All rays in $\rm{Vert}(\Del^\sig)\setminus \rm{Vert}(\sig)$ are contained in $\underset{\tau\leq\sig}{\bigcup}\inte(\tau)$.
\item{3$^\circ$} $\Del^\sig$ is $\rm{Aut}(\sig)$-invariant.
\end{definition}

\begin{remark} Condition $3^\circ$ is replaced with a stronger one in the following proposition.
\end{remark}
\begin{lemma} \label{sub} If $\tau'\in\overline{\tau}$, $\tau'\prec \sigma\in\Sigma$ then  $\rm{Vert} (\Del^\sig_{|\tau'})\setminus\rm{Vert}(\tau')\subset\tau$ and thus 
\[\Del^\sig_{|\tau'}=\Del^\sig_{|\tau}\oplus\la e_1,\ldots,e_k\ra=\Del^\tau\times \la e_1,\ldots,e_k\ra.\]
\end{lemma}
 \begin{lemma} \label{maximal2} 
For every point $x\in B\setminus ({B_+\cap B_-})$, $x\in s'$ there exists a toric chart $x\in U_\sig\to X_\sig$, with $\Ga_\sig=K^*$, corresponding to a stratum $s\subset \overline{s'}$ . In particular the maximal cones of $\Sig$ are circuits.
\end{lemma}
\begin{proof}
Let $\tau$ correspond to a stratum $s'\ni x$. By definition of cobordism  $\ds\lim_{t\to 0}tx=x_0$ or $\ds\lim_{t\to \infty}tx=x_0$ exists. The point $x_0$ is $\ks$-fixed and belongs to a stratum $s$, with $\Ga_s=\Ga_\sig=\ks$. Since $U$ is a $\ks$-invariant neighborhood of $x_0$ it contains an orbit $K^*\cdot x$ and the point $x$. Moreover  $\overline{s'}\supset s $ and $\tau\leq \sig$.
\end{proof}
 
 \begin{lemma} Let $\sig$ be the cone corresponding to a stratum $s$ on $B$ and $x\in s$. Then $\widehat X_x=\Spec \widehat O_{x,B}\simeq (X_\sig\times \bbA^{\dim(s)})^\wedge\cong\,\Spec K [[x_1,\ldots,x_k,\ldots, x_n]].$ 
%
%
\end{lemma}
Set $\widetilde X_\sig:=(X_\sig\times \bbA^{\dim(s)})^\wedge$ and let $G_\sig$ denote the group of all $\Ga_\sig$-equivariant autorphisms of $\widetilde X_\sig$.

The subdivision $\Del^\sig$ of $\sig$ defines a toric morphism and induces a proper birational $\Ga_\sig$-equivariant morphism
\[\wtx_{\Del^\sig}:=X_{\Del^\sig} \times_{X_\sig}\wtx_\sig\to\wtx_\sig.\]

\noindent
\begin{proposition} \label{correspondence} Let $\Del=\{\Del^\sig\mid\sig\in\Sig\}$ be a subdivision of $\Sig$ such that:
\begin{enumerate}
\item {For every $\sig\in\Sig$ the morphism $\wtx_{\Del^\sig}\to\wtx_\sig$ is $G_\sig$-equivariant}.
\end{enumerate}
Then $\Del$ defines a $\ks$-equivariant birational modification $f:B'\to B$ such that for every  toric chart $\varp_\sig:U\to X_\sig$ there exists a $\Ga_\sig$-equivariant fiber square
\[\begin{array}{rcccccccc}
U_\sigma \times_{X_{\sigma}}
X_{\Delta^\sigma} &&\simeq& f^{-1}(U_\sigma)
 & \rightarrow &
X_{\Delta^\sigma} &&& \\
&&&\downarrow {\scriptstyle f} & & \downarrow  &&&\\
&& & U_\sig &  \rightarrow & X_{\sigma}&&&(2)\\

\end{array}\] 

\end{proposition}
\begin{definition} A decomposition $\Del$ of $\Sig$ is {\it canonical} if it satisfies  condition (1). 
\end{definition}
\begin{proof} The above diagrams define open subsets $f^{-1}_\sig(U_\sig)$ together with proper birational $\Ga_\sig$-equivariant morphisms $f^{-1}_\sig(U_\sig)\to U_\sig$. Let  $ s'\leq s$ be a stratum corresponding to the cone $\tau\leq\sig$. By Lemma \ref{sub}, the restriction of the diagram (2) defined by $U_\sig\to X_\sig$ to  a neighboorhod $U_\tau$ of $y\in s'$ determines   a diagram defined by the induced toric chart $U_\tau\to X_\tau$ and the decomposition $\Del^\tau$ of $\tau$.
In order to show that the $f^{-1}_\sig(U)$ glue together we need to prove that for $x\in s$ and two different charts of the form $\varp_{1,\sig}:U_{1,\sig}\to X_\sig$ and $\varp_{2,\sig}:U_{2,\sig}\to X_\sig$ where $x\in U_{1,\sig}, U_{2,\sig}$ the induced varieties $V_1:=f^{-1}_{1,\sig}(U_{1,\sig})$ and $V_2:=f^{-1}_{2,\sig}(U_{2,\sig})$ are isomorphic over $U_{1,\sig}\cap U_{2,\sig}$.
\noindent
For simplicity assume that $U_{1,\sig}=U_{2,\sig}=U$ by shrinking $U_{1,\sig}$ and $U_{2,\sig}$ if necessary. The charts $\varp_{1,\sig}, \varp_{2,\sig}:U\to X_\sig$ are defined by the two sets of semiinvariant parameters, $u^1_1,\ldots,u^1_k$ and $u^2_1,\ldots,u^2_k$
 with a nontrivial action of $\Ga_\sig$. These sets can be extended to full sets of parameters $u^1_1,\ldots,u^1_k,u_{k+1},\ldots,u_n$ and $u^2_1,\ldots,u^2_k, u_{k+1},\ldots,u_n$ where $\Ga_\sig$ acts trivially on $u_{k+1},\ldots,u_n$,  and $u_{k+1}\ldots,u_n$ define parameters on the stratum $s$ at $x$.
These two sets of parameters define \'etale morphisms $\varp_{1,\sig},\varp_{2,\sig}:U\to X_\sig\times \bbA^{n-k}$ and fiber squares
\[\begin{array}{rcccccc}
\overline{\varp_{i,\sig}}:& V_i
 & \rightarrow &X_{\Del^\sig}\times \bbA^{n-K}
 &&& \\
&\downarrow & & \downarrow  &&&\\
\varp_{i,\sig}:& U &  \rightarrow & X_\sig \times \bbA^{n-K}&&&\\

\end{array}\] 
\noindent
Suppose the induced $\Ga$-equivariant birational map $f:V_1\dashrightarrow V_2$ is not an isomorphism over $U$.

Let $V$ be the graph of $f$ which is a dominating component of the fiber product $V_1\times_U V_2$. Then either $V\to V_1$ or $V\to V_2$ is not an isomorphism (i.e.~collapses  a curve to a point)
 over some $x\in s\cap U$. Consider an \'etale $\Ga_\sig$-equivariant morphism $e:\whx_x\to U$. Pull-backs of the morphisms $V_i\to U$ via $e$ define 
 two different nonisomorphic $\Ga_\sig$-equivariant
 liftings $Y_i\to {\whx_x}$, since the graph $Y$ of $Y_1 \dashrightarrow Y_2$ (which is a pull-back of $V$) is not isomorphic to at least one $Y_i$.  But these two liftings are defined by two isomorphisms  $\widehat\varp_1,\widehat\varp_2:  {\whx_x}\simeq \wtx_{\sig}$.
These isomorphisms differ by  some automorphism  $g\in G_\sig$, so we have
 $\widehat\varp_1=g\circ \widehat\varp_2$. Since $g$ lifts to the automorphism of $ \wtx_{\Del^\sig} $ we get $Y_1\simeq Y_2\simeq  \wtx_{\Del^\sig}$, which contradicts  the choice of $Y_i$.

 Thus $V_1$ and $V_2$ are isomorphic over any $x\in s$ and $B'$ is well defined by glueing pieces $f^{-1}_\sig(U)$ together. We need to show that the action of $\ks$ on $B$ lifts to the action of $\ks$ on $B'$.

Note that $B'$ is isomorphic to $B$ over the open generic stratum $U\supset B_+\cup B_-$ of points $x$ with $\Ga_x=\{e\}$. By Lemma \ref{maximal2} every point $x\in B\setminus ({B_+\cap B_-})$ is in $U_\sig$, with $\Ga_\sig=K^*$.
Then the diagram (2) defines the action of $\ks$ on $f^{-1}(U_\sig)$.
\end{proof}
\subsection{Simple properties of $\wtx_{\Del^\sig}$}

Recall that $\wtx_\sig=\Spec(K[[x_1,\ldots,x_k,x_{k+1},\ldots,x_n]])$,
where \\$X_\sig=\Spec(K[x_1,\ldots,x_k])$ and $\Ga_\sig$ acts trivially on  $x_{k+1},\ldots,x_n$. This gives us $$X_{\widetilde{\sig}}:=\Spec(K[x_1,\ldots,x_k,x_{k+1},\ldots,x_n])=\Spec(K[x_1,\ldots,x_k])\times \Spec(K[x_{k+1},\ldots,x_n])=X_\sig\times X_\reg(\sig)$$
where ${\widetilde{\sig}}$ and   $\reg(\sig)$ correspond to  $\Spec(K[x_1,\ldots,x_k,x_{k+1},\ldots,x_n])$ and $\Spec(K[x_{k+1},\ldots,x_n]$
 respectively.
We can write $\widetilde{\sig}=\sig\times\reg(\sig)$, 
$\widetilde{\sig}^\vee=\sig^\vee\times\reg(\sig)^\vee$, 
$N_{\widetilde{\sig}}=N_\sig\times N_{\reg(\sig)}$, and  $M_{\widetilde{\sig}}=M_\sig\times M_{\reg(\sig)}$.

Let $\Del^\sig$ be a subdivision of $\sig$. There is a natural morphism
$$j_{\Del^\sig}:\wtx_{\Del^\sig}\to X_{\Del^\sig}$$

\begin{lemma} 
\begin{enumerate}
\item The open cover  $\{X_\tau\mid \tau \in \Del^\sig\}$ of $X_{\Del^\sig}$ defines the open cover of $\{\wtx_\tau\mid \tau \in \Del^\sig\}$ of $\wtx_{\Del^\sig}$, where  $\wtx_\tau:=X_\tau\times_{X_\sig}\wtx_\sig=j_{\Del^\sig}^{-1}(X_\tau)$ and
 $K[\wtx_\tau]=K[\tau^\vee]\otimes_{K[\sig^\vee]}K[[\widetilde{\sig}^\vee]]$
\item The closed orbits $O_\tau\subset X_\tau$ define closed subschemes $\widetilde{O_\tau}:=j_{\Del^\sig}^{-1}(O_\tau)$ of $\wtx_\tau\subset \wtx_{\Del^\sig},$
\noindent where
$K[\widetilde{O_\tau}]={K[\tau^\perp]}\otimes_{K[\sig^\vee\cap\tau^\perp]}K[[({\sig}^\vee\cap\tau^\perp)\times \reg(\sig)^\vee]].$
\item  The local ring $\cO_{\wtx_{\Del^\sig}, \widetilde{O_\tau}}$ at the generic point of $\widetilde{O_\tau}$ contains the residue field $ K(\widetilde{O_\tau})$ (which is a quotient of $K[\widetilde{O_\tau}]$).
The completion of $\cO_{\wtx_{\Del^\sig}, \widetilde{O_\tau}}$ is of the form $$\widehat {\cO_{\wtx_{\Del^\sig}, \widetilde{O_\tau}}} \buildrel{\varp}\over{\simeq} K(\widetilde{O_\tau})
[[{\underline\tau}^\vee]]\supset \cO_{\wtx_{\Del^\sig}, \widetilde{O_\tau}}\supset K(\widetilde{O_\tau})[{\underline\tau}^\vee].$$ 
\item The group $\Gamma_\tau\subset \Ga_\sig$ acts tivially on $K(\widetilde{O_\tau})$.The action of $\Ga_\tau$ on characters $\tau^\vee$ descends to $\underline{\tau}^\vee=\tau^\vee/\tau^\perp$.
In particular if $\tau$ is dependent, the action of $\ks$ on $K(\widetilde{O_\tau})$ is trivial. \end{enumerate}
\end{lemma}

\begin{proof} (1) follows from definition. The elements of $K[\tau^\vee]\otimes_{K[\sig^\vee]}K[[\widetilde{\sig}^\vee]]$ are the finite sums of the form $\sum x_if_i$, where $x_i\in\tau^\vee$ and $f_i\in K[[\widetilde{\sig}^\vee]]=K[[{\sig}^\vee\times\reg(\sig)^\vee]]$ is an infinite power series. Note also that $\sig^\vee\subset \tau^\vee$.
 
 (2) The ideal $I_{O_\tau}\subset K[X_\tau]$ is generated by all characters $x^F$, where $F\in \tau^\vee\setminus\tau^\perp$. 
These characters generate  the ideal $I_{\widetilde{O_\tau}}\subset K[\wtx_\tau]$.
Then the elements of 
$K[\widetilde{O_\tau}]=K[\wtx_\tau]/{I_{\widetilde{O_\tau}}}$ are  the finite sums of the form $\sum x_if_i$, where $x_i\in\tau^\perp$ and $f_i\in K[[(\sig^\vee\cap\tau^\perp)\times\reg(\sig)^\vee]]$ is an infinite power series. We get $K[\wtx_\tau]/{I_{\widetilde{O_\tau}}}=
K[\tau^\perp]\otimes_{K[\sig^\vee\cap\tau^\perp]}K[[(\sig^\vee\cap\tau^\perp)\times\reg(\sig)^\vee]]$.

(3) Note that $K[\widetilde{O_\tau}]$ is a subring of $K[\wtx_\tau]$. The subalgebra generated by $\tau^\vee$ over $K[\widetilde{O_\tau}]$ is equal to $K[\widetilde{O_\tau}][\tau^\vee]=K[\widetilde{O_\tau}][\underline{\tau}^\vee]\subset K[\wtx_\tau]\subset K[\widetilde{O_\tau}][[\underline{\tau}^\vee]]$. Passing to the localizations at $I_{\widetilde{O_\tau}}$ we get inclusions $K(\widetilde{O_\tau})[\underline\tau^\vee]\subset(K[\wtx_\tau])_{\widetilde{O_\tau}}={\cO_{\wtx_{\Del^\sig},\widetilde{O_\tau}}} \subset K(\widetilde{O_\tau})[[\underline\tau^\vee]]=\widehat {\cO_{\wtx_{\Del^\sig},\widetilde{O_\tau}}} $.

(4) The action of $\Ga_\tau$ on $K(O_\tau)=K[\tau^\perp]$ is trivial. Then
$\Ga_\tau$ acts trivially on all characters in $\tau^\perp$ and on $K[\widetilde{O_\tau}]=
K[\tau^\perp]\otimes_{K[\sig^\vee\cap\tau^\perp]}K[[(\sig^\vee\cap\tau^\perp)\times\reg(\sig)^\vee]]$.
\end{proof}

\bigskip
\subsection{Basic properties of valuations}

Let $K(X)$ be the field of rational functions on an algebraic variety or an integral scheme $X$.
A {\it valuation} on  $K(X)$ is a group homomorphism $\mu:K(X)^*\to G$ from the multiplicative group  $K(X)^*$ to a totally ordered group $G$ such that $\mu(a+b)\geq\min(\mu(a),\mu(b))$.
By the {\it center} of a valuation $\mu$ on $X$ we mean an irreducible closed subvariety $Z(\mu)\subset X$ such that for any open affine $V\subset X$, intersecting $Z(\mu)$, the ideal $I_{Z(\mu)\cap V}\subset K[V]$ is generated by all $f\in K[V]$ such that $\mu(f)>0$ and for any $f\in K[V]$, we have $\mu(f)\geq 0$.

Each vector $v\in N^{\QQ}$  defines a linear function on $M$ which
determines a 
valuation $\val(v)$ on a toric variety
$X_{\Del}\supset T$.
For any regular function $f=\sum_{w\in M} a_wx^w\in K[T]$ set
$$\val(v)(f):=\min\{(v,w)\mid a_w\neq 0\}.$$ 

If $v\in\inte(\sig)$, where $\sig\in\Del$, then $\val(v)$ is positive for all $x^F$, where $F\in\sig^\vee\setminus\sig^\perp$. In particular we get
$$Z(\val(v))=\overline{O_\sig}\quad\quad\mbox{iff}\quad\quad v\in\inte\sig.$$

\noindent
If $v\in\sig$ then $\val(v)$ is a valuation on $R=K[X_\sig]=K[\sig^\vee]$, that is, $\val(v)(f)\geq 0$ for all $f\in K[\sig^\vee]\setminus\{0\}$.
We construct ideals for all $a\in \bN$ which uniquely determine $\val(v)$:
\[I_{\val(v),a}=\{f\in R\mid \val(v)(f) \geq a\}=(x^F\mid F\in\sig^\vee, F(v)\geq a)\subset R.\]
%
%
By glueing $I_{\val(v),a}$
%
%
 for all $v\in \sig$ and putting  
$\cI_{\val(v),a|X_\sig}=\cO_{X_\sig}$ if $v\notin\sig$ we construct a coherent sheaf of ideals $\cI_{\val(v),a}$
%
%
 on $X_\Del$.

 Let $\sig\in\Sig$ be a cone of the semicomplex $\Sig$ and $v\in\sig\subset\widetilde{\sig}$. The valuation $\val(v)$ on $K[\sig^\vee]$ extends to  the valuation on $K[[\widetilde{\sig}^\vee]]$. Thus it determines a valution on $\wtx_{\Del^\sig}$, where $\Del^\sig $ is a subdivision of $\sig$. As before we have
\begin{lemma}
 \begin{enumerate}
 \item
$Z(\val(v),\wtx_{\Del^\sig})=\cl{(\widetilde{O_\tau})}\subset \wtx_{\Del^\sig}$,  where $\tau\in\Del$ and  $v\in\inte(\tau)$.
\item There exists a coherent sheaf of ideals $\cI_{\val(v,a),\wtx_{\Del^\sig}}=j_{\Del^\sig}^*(\cI_{\val(v,a),X_{\Del^\sig}})$ on $\wtx_{\Del^\sig}$ such that for every $\del\in\Del
$ containing $v$ and $R=K[\wtx_\del]=K[\del^\vee]\otimes_{K[\sig^\vee]}K[[\widetilde{\sig}^\vee]]$ we have \[I_{\rm{val}\,(v),a}=\{f\in R\mid \rm{val}\,(v)(f) \geq a\}=(x^F\mid F\in\sig^\vee, F(v)\geq a)\subset R.\]
\item The valuation $\val(v)$ on the local ring $\cO_{\whx_{\Del^\sig},\Otau}$, where $v\in\tau$ extends to its completion $\widehat{\cO}_{\whx_\Del,\Otau}=K(\Otau)[[\underline{\tau}^\vee]]$. Moreover $\rm{val}\,(v)_{| K(\Otau)^*}=0$. \end{enumerate}
\end{lemma}

\begin{lemma} \label{closed}

 If $\cl{(\Otau)}\subset \wtx_{\Del^\sig}$ is $G_\sig$-invariant
then $\val(v)$ is $G_\sig$-invariant on $\wtx_{\Del^\sig}$ iff it is $G_\sig$-invariant on $\whx_\tau:=\Spec(K(\Otau)[[\underline{\tau}^\vee]])$.

\end{lemma}
\begin{proof}
$K(\Otau)[[\underline{\tau}^\vee]]$ is faithfully flat over $\cO_{\whx_{\Del^\sig},\Otau}$ and we have $1-1$ correspondence between
ideals $g^*(I_{\val(v),a})$ on $\cO_{\wtx_{\Del^\sig},\Otau}$ and on $K(\Otau)[[\underline{\tau}^\vee]]$.
\end{proof}


\subsection{Blow-ups of toric ideal sheaves}
The sheaf $\cI_{\rm{val}\,(v),a}$ is an example of an $T$-invariant sheaf of ideals  on a toric variety $X_\Del$. It is locally defined by monomial ideals $I_\sig\subset K[X_\sig]$. Any $T$-invariant sheaf of ideals $\cI$ on $X_\Del$ defines a function $\rm{ord}_{\cI}:|\Sig|\to Q$ (see \cite{KKMS}) such that for any $p\in\sig$
\[\rm{ord}_{\cI}(p)=\min\{F(p):x^F\in I_\sig\}\]

\noindent
The function $\rm{ord}_{\cI}$ is concave and piecewise linear on every cone $\sig\in\Del$. If $(x^{F_1},\dots,x^{F_k})= I_\sig$ for $F_1,\ldots,F_k\in\sig^\vee$ then $\rm{ord}_{\cI}(p)=\min(F_1(p),\dots, F_k(p))$.  The cones $\sig_{F_i}:=\{p\in\sig:\rm{ord}_{\cI}(p)=F_i(p)\}$ define a subdivision of $\sig$ and by combining these subdivision together for all $\sig\in \Del$ we get a subdivision $\Del_{\rm{ord}_{\cI}}$ of $\Del$. This is the coarsest subdivision of $\Del$ for which $\rm{ord}_{\cI}$ is linear on every cone.

\begin{lemma}\label{le:kk}
\cite{KKMS} If $\cI$ is an iraviant sheaf of ideals on $X_\Del$ then the normalization of the blow--up of $\cI$ on $X_\Del$ is a toric variety $X_{\Del_{\rm{ord}_{\cI}}}$ corresponding to the subdivision $\Del_{\rm{ord}_{\cI}}$ of $\Del$.
\end{lemma}

\begin{proof} Let $f:X'\to X_\Del$ be the normalized blow--up of $\cI$. Then $X'$ is a toric variety on which $f^*(I)$ is locally invertible. Then $X'$ corresponds to a subdivision $\Del'$ of $\Del$ such that $\rm{ord}_{f^*(\cI)}=\rm{ord}_{\cI}$ is linear on every cone $\Del'$. From the universal property of the blow--up we conclude that $\Del'$ is the coarsest subdivision with this property. Thus $\Del'=\Del_{\rm{ord}_{\cI}}$.
\end{proof}

\noindent

\begin{lemma}\label{blow-up}\cite{KKMS} Given a simplical fan and an integral vector $v\in|\Del|$, there exists a sufficiently divisible natural number $a$, such that
$\Del_{{\rm{ord}_\cI}_{\rm{val}\,(v),a}}=\la v\ra\cdot\Del. $
\end{lemma}

\begin{proof} Let $\sig=\la e_1,\dots,e_k\ra$ be a cone containing $v$ and  assume that $v\in \inte \la e_1,\dots,e_\ell\ra\preceq \sig$, for some $\ell\leq k$. Let $F_j\in\sig^\vee$, for $1\leq j\leq\ell$, be the functional such that $F_j(e_i)=0$ for $i\neq j$ and $F_j(v)=a$. If $a$ is sufficiently divisible then $F_j$ is integral and $x^{F_j}\in I_{\rm{val}\,(v),a}$ for all $1\leq j\leq \ell$. Note that for any $x^F\in I_{\rm{val}\,(v),a}$ we have that  $F(v)\geq a$ and $F(e_i)\geq 0$. This gives $F\geq F_j$ on $\la e_1,\dots,\check e_j,\dots,e_k,v\ra$ and finally ${\rm{ord}_\cI}_{\rm{val}\,(v),a}=F_j$ on $\la e_1,\dots,\check e_j,\dots,e_{k},v\ra$. Note that since $F_j(e_j)>0$, we have that $F_j(p)>F_i(p)$ if $p\in \la e_1,\dots,\check e_i,\dots,e_k,v\ra\setminus\la e_1,\dots,\check e_j,\dots,e_k,v\ra$ so ${\rm{ord}_\cI}_{\rm{val}\,(v),a}=F_j$ exactly on $\la e_1,\dots,\check e_j,\dots,e_k,v\ra$ and $(v)\cdot \sig=\sig_{{\rm{ord}_\cI}_{\rm{val}\,(v),a}}$.
\end{proof}

\subsection{Stable vectors}
Let $g:X\to Y$ be any dominant morphism of integral schemes (that is,
$\overline{g(X)}=Y$) and $\mu$ be a valuation of $K(X)$. Then $g$ induces a valuation $g_*(\mu)$ on $K(Y)\simeq g^*(K(Y))\subset K(X)$:
$g_*\mu(f)=\mu(f\circ g)$.
\begin{definition} Let $\Sig$ be the semicomplex defined for the cobordism $B$. A vector $v\in \rm{int}(\sig)$, where $\sig\in\Sig$, is called \textit{stable} if for every $\sig\leq\sig'$, $\val(v)$ is $G_{\sig'}$-invariant on $\wtx_{\sig'}$.
\end{definition}

\begin{lemma} \label{G} If $\wtx_{\Del^\sig}\to\wtx_{\sig}$ is $G_\sigma$-equivariant and $\val(v)$ is $G_\sigma$-invariant then  $\wtx_{\la v\ra\cdot \Del^\sig}\to\wtx_{\sig}$
is $G_\sigma$-equivariant.
\end{lemma}

\begin{proof}  The morphism $\wtx_{\la v\ra\cdot\Del^\sig}\to\wtx_{\Del^\sig}$ is a pull-back of the morphism $X_{\la v\ra\cdot\Del^\sig}\to X_{\Del^\sig}$. Thus, by Lemma \ref{blow-up},  $\wtx_{\la v\ra\cdot\Del^\sig}\to\wtx_{\Del^\sig}$ is a normalized blow-up of $\cI_{\val(v),a}$ on $\wtx_{\Del^\sig}$. But the latter sheaf is $G_\sigma$-invariant.
\end{proof}

\begin{proposition}\label{blowup} Let $\Del=\{\Del^\sig\mid \sig\in\Sig\}$ be a canonical subdivision of $\Sig$ and $v$ be a stable on $\Sig$. Then $\la v\ra\cdot\Del:=\{\la v\ra\cdot\Del^\sig\mid \sig\in\Sig\}$ is a canonical subdivision of $\Sig$.
\end{proposition}

\subsection{Convexity}

\begin{lemma}\label{convex}  Let $\val(v_1)$ and $\val(v_2)$ be $G_\sig$-invariant valuations on $X_\sig$. Then all valuations $\val(v)$, where $v=av_1+bv_2$, $a,b\geq 0$, $a,b\in \QQ$, are $G_\sig$-invariant.
\end{lemma}

\begin{proof} Let $\Del=\la v_1\ra\cdot\la v_2\ra\cdot\sig$ be a subdivision of $\sig$. Then by Lemma \ref{G}, $\wtx_{\Del^\sig}\to\wtx_\sig$ is $G_\sig$-invariant on  $\wtx_{\Del^\sig}$. The exceptional divisors $D_1$ and $D_2$ are $G_\sig$-invariant and correspond to one-dimensional cones (rays) $\la v_1\ra,\la v_2\ra\in\Del$. The cone $\tau=\la v_1,v_2\ra\in D$ corresponds to the orbit $\widetilde{O_\tau}$ whose closure is $D_1\cap D_2$ and thus the generic point is $G_\sig$-invariant. The action of $G_\sig$ on $\wtx_\sig$ induces an action on the local ring $\wtx_{\Del,\widetilde{O_\tau}}$ at the generic point of $\widetilde{O_\tau}$ and on its completion $K(\widetilde{O_\tau})[[\underline\tau^\vee]]$. Note that for any $v\in\tau$, $\val(v)_{|K(\widetilde{O_\tau})}=0$. For any  $F\in\underline\tau^\vee=\ds{\tau^\vee\over\tau^\perp}$ the divisor $(x^F)$ of  the character $x^F$ on $\whx_\tau:=\Spec  K(\widetilde{O_\tau})[[\underline\tau^\vee]]$ is a combination $n_1D_1+n_2D_2$ for $n_1n_2\in\bZ$. Since $D_1$ and $D_2$ are $G_\sig$-invariant, the divisor $(x^F)=n_1D_1+n_2D_2$ is $G_\sig$-invariant, that is, for any $g\in G$, we have $g x^F=u_{g,F}\cdot x^F$ where $u_{g,F}$ is invertible on $K(\widetilde{O_\tau})[[\underline\tau^\vee]]$.
Thus for every $v\in\tau$ and $g\in G$ we have 
\begin{equation*}
\begin{aligned}
g^*(I_{\val(v),a})&=g^*(x^F | F\in\underline\sig^\vee, F(v)\geq a)
=(u_{g,F}x^F | F\in\sig^\vee, F(v)\geq a)=I_{\val(v),a}.
\end{aligned}
\end{equation*}
Thus $\val(v)$ is $G_\sig$-invariant on $K(\widetilde{O_\tau})[[\underline\tau^\vee]]$ and on its subring $\cO_{\wtx_{\Del,\widetilde{O_\tau}}}$. The latter ring has the same quotient field as $\wtx_\sig$ so $\val(v)$ is $G_\sig$-invariant on $\wtx_\sig$.
\end{proof}

\begin{lemma}\label{convex2}
Let $\sig\in\Sig$ and $v_1,v_2\in\sig$ be stable vectors. Then all vectors $v=av_1+bv_2\in\sigma$, where $a,b\in \QQ_{>0}$, are stable.
\end{lemma}
\subsection{Existence of quotients}
\begin{lemma}\label{le:q} 
Let $\Ga\subset \Ga_\sig$ be a finite subgroup and $\tau\in\Del^\sig$.  
Then $\wtx_{\tau}/\Ga=X_{\tau}/{\Ga}\times_{X_{\sig}/\Ga}\wtx_\sig/\Ga$.
\end{lemma}
\begin{proof} The group $\Ga\simeq \ZZ_n$ acts on characters $x^F$, $F\in M_\sig$, with weights  $a_F$: $t(x^F)=t^{a_F}x^F$ where $t\in\Ga$, and $a_F\in \bZ_n$.  The elements of the ring 
\[K[\wtx_\tau]=K[X_\tau]\otimes_{K[X_\sig]}K[\wtx_\sig]=K[\tau^\vee]\otimes_{K[\sig^\vee]}K[[\widetilde{\sig}^\vee]]\]
are finite sums $\sum x_if_i$ where $f_i\in K[\wtx_\sig]$ is a formal power series and $x_i\in\tau^\vee$ is a character. (Note that $\sig^\vee\subseteq\tau^\vee$ since $\tau\subset \sig$.)
The elements of the ring $K[\wtx_{\tau}]/\Ga=K[X_\tau]^\Ga$ are finite sums $\sum x_if_i$ of weight zero, that is, every $f_i\in K[[\sig^\vee]]$ is a quasihomogeneous power series of weight $a_{f_i}=-a_{x_i}$. The elements of the ring $K[X_\tau]^\Ga\otimes_{K[X_\sig]^\Ga}K[\wtx_\sig]^\Ga$ are of the form $\sum x_if_i$ where $x_i$ and $f_i$ each have  weight zero. We have to prove 
\begin{lemma} \label{le:10} Let $K[\widetilde{\sig}^\vee]=\underset{a\in \ZZ_n}{\oplus}K[\widetilde{\sig}^\vee]^a$ be
a decomposition according to weights. Then $K[\widetilde{\sig}^\vee]^a$ is generated over $K[\widetilde{\sig}^\vee]^0$ by
finitely many monomials.
\end{lemma}
\begin{proof} Note that for any $F\in \widetilde{\sig}^\vee$, the element $nF\in (\widetilde{\sig}^\vee)^a$.
Let $x^{F_1},\dots,x^{F_k}$
generate $K[\widetilde{\sig}^\vee]$. 
Then the elements
$x^{\alpha_1F_1+\cdots+\alpha_kF_k}$, where $$\alpha_1F_1(v_\sig)+\cdots
+\alpha_kF_k(v_\sig)=a\quad  \mbox{and}\quad  0\leq\alpha_i\leq n$$ generate
$K[\widetilde{\sig}^\vee]^a$ over $K[\widetilde{\sig}^\vee]^0$. 
\end{proof}
It follows from Lemma \ref{le:10} that every $f_i\in K[[\widetilde{\sigma}^\vee]]^a$ decomposes as a finite sum $f_i=\sum x_{ij}f_{ij}$, where $f_{ij}\in K[[\widetilde{\sigma}^\vee]]^0$ and $x_{ij}\in(\sigma^\vee)^a$ and
$f=\underset{i}{\sum}x_if_i=\underset{ij}{\sum}x_ix_{ij}f_{ij}\in K[\tau^\vee]^\Ga\otimes_{K[\sig^\vee]^\Ga}K[[\widetilde{\sig}^\vee]]^\Ga $.
\end{proof}
\begin{corollary}\label{quo1}
Let a finite group $\Ga\subset \Ga_\sig$ act on $\wtx_\sig$. Then for a decomposition $\Del^\sig$ of $\sig$ the following quotient exists:
\[\wtx_{\Del^\sig}/\Ga=X_{\Del^\sig}/\Ga\times_{X_{\sig}/\Ga}\wtx_\sig/\Ga=X_{\pi({\Del^\sig})}\times_{X_{\pi(\sig)}}\wtx_{\pi(\sig)}\]
where $\pi:\sig\to \pi(\sig)$ corresponds to the quotient $X_\sig\to X_\sig/\Ga$ and $\wtx_{\pi(\sig)}:=(X_{\pi(\sig)}\times X_{\reg(\sig)})^\wedge=\wtx_{\sig}/\Ga$.
\end{corollary}
As before 
there is a natural morphism
$$j_{\pi(\Del^\sig}):\wtx_{\pi(\Del^\sig)}\to X_{\pi(\Del^\sig)}$$

\begin{lemma} 
\begin{enumerate}
\item The open cover $\{X_{\pi(\tau)}\mid \pi(\tau)\in\pi(\Sig)\}$ of $X_{\pi(\Del^\sig)}$ defines the open cover \\ $\{\wtx_{\pi(\tau)}\mid \pi(\tau)\in\pi(\Sig)\}$ of $\wtx_{\pi(\Del^\sig)}$, where   $\wtx_{\pi(\tau)}:= \wtx_{\tau}/\Ga= 
X_{\pi(\tau)} \times_{X_{\pi(\sig)}}\wtx_{\pi(\sig)}=j_{\pi(\Del^\sig)}^{-1}(X_{\pi(\tau)})$.
\item The closed orbits $O_{\pi(\tau)}\subset X_{\pi(\tau)}$ define closed subschemes $\widetilde{O}_{\pi(\tau)}:=j_{\pi(\Del^\sig)}^{-1}(O_{\pi(\tau)})$ of $\wtx_{\pi(\tau)}\subset \wtx_{\pi(\Del^\sig)},$
\noindent where
$K[\widetilde{O}_{\pi(\tau)}]=K[\widetilde{O}_{\tau}]/\Ga.$
\item The completion of the local ring $\cO_{\wtx_{\pi(\Del^\sig}), \widetilde{O}_{\pi(\tau)}}$ at the generic point of $\widetilde{O}_{\pi(\tau)}$ is of the form $$\widehat {\cO}_{\wtx_{\pi(\Del)}, \widetilde{O}_{\pi(\tau)}}{\simeq} K(\widetilde{O}_{\pi(\tau)})
[[{{\pi(\underline\tau)}}^\vee]].$$  
\end{enumerate}
\end{lemma}

\subsection{Descending of  the group action of $G_\sig$.}
\begin{lemma}\label{21}
If $V\subset \wtx_{\Del^\sig}$ is an open affine $\Ga$--invariant subscheme then for any open affine $\Ga$--invariant subscheme $U\subset V$, we have an open inclusion of schemes $
U/\Ga\subset V/\Ga$.
\end{lemma}

\begin{proof} Let $Z=V\setminus U$ be a closed affine subscheme. Then $I_Z\subset K[U]$ is $\Ga\simeq\ZZ^n$ invariant and generated by a finite number of semiinvariant functions $f_1,\dots, f_k\in I_Z$. Then the functions $g_1=f^n_1,\dots,g_k=f^n_k$ are invariant. Write $U=\underset{i=1}{\bigcup} V_{g_i}$ as the union of open subschemes $V_{g_i}=U_{g_i}$. The algebra $K[V_{g_i}]^\Ga=K[V]^\Ga_{g_i}=K[U]_{g_i}$ the localization of $K[V]$ so there is an open inclusion $V_{g}/\Ga\subset U/\Ga$ and $V_{g}/\Ga\subset V/\Ga$. It follows that $U/\Ga=\underset{i}{\bigcup}V_{i}/\Ga\subset V_{g}/\Ga$.
\end{proof}

\begin{lemma}\label{22}
Any open $\Ga$--equivant embedding of an open affine $\Ga$--invariant subscheme $V$ into $\wtx_{\Del^\sig}$ determines an open
embedding  $V/\Ga\subset \wtx_{\Del^\sig}/\Ga$.
\end{lemma}
\begin{proof} Let $V_\tau=V\cap\wtx_\tau$, where $\tau\in{\Del^\sig}$. Then $V=\underset{\tau\in{\Del^\sig}}{\bigcup}V_\tau$. By the previous lemma $V_\tau/\Ga\subset \wtx_\tau/\Ga$ and $V_\tau/\Ga\subset V/\Ga$ are open embeddings defining an open inclusion $V/\Ga=\underset{\tau\in{\Del^\sig}}{\bigcup}V_{\tau}/\Ga\subset \wtx_{{\Del^\sig}}/\Ga$.
\end{proof}
\begin{lemma}\label{2}
The action of $G_\sig$ descends to $\wtx_{\Del^\sig}/\Ga$ and we have a $G_\sig$--equivariant morphism $\wtx_{\Del^\sig} \to \wtx_{\Del^\sig}/\Ga$.
\end{lemma}

\begin{proof} The lemma is immediate consequence of Lemma \ref{22}. For any $g\in G_\sig$, the morphism $g: \wtx_{{\Del^\sig}}/\Ga\to \wtx_{{\Del^\sig}}/\Ga$ is defined locally by $g:  \wtx_{{\Del^\sig}}/\Ga\supset V/\Ga\to gV/\Ga\subset \wtx_{{\Del^\sig}}/\Ga$.
\end{proof}

\subsection{Basic properties of stable vectors}

\begin{lemma} Let $\Tan_0=\bbA^n=\Tan_0^{a_0}\oplus \Tan_0^{a_1}\oplus\cdots\oplus\ \Tan_0^{a_k}$ denote the tangent space of $\wtx_\sig=\Spec K[[u_1,\dots,u_n]]$ at $0$ and its decomposition according to the weight distribution. Let $d:G_\sig\to \Gl(\Tan_0)$ be the differential morphism defined as $g\mapsto dg$. Then $d(G_\sig)=\Gl(\Tan_0^{a_1})\times\cdots\times \Gl(\Tan_0^{a_k})$.
\end{lemma}
\begin{proof}  The elements of the group $g\in G_\sig$ are defined by $(u_1,\dots,u_n)\mapsto(q_1,\dots,q_n)$ where $g_i=g^*(u_i)$ are quasihomogenious  power series of $\Ga_\sig$-weights $a(g_i)=a(u_i)$. 
\end{proof}
\begin{lemma}\label{le:semiinv}
Let $v\in\sig$, where $\sig\in\Sig$, be an integral  vector such that for any $g\in G_\sig$, there exists  an integral vector $v_g\in \sig$ such that 
$g_*(\val(v))=\val(v_g)$. Then $\val(v)$ is $G_\sig$-invariant on $\widetilde{X}_\sig$.
\end{lemma}

\begin{proof} Set $W=\{v_g\mid g\in G\}$. For any $a\in \NN$, the ideals $I_{\val(v_g), a}$ are generated by monomials. They define the same Hilbert--Samuel function $k\mapsto\dim_K(K[\wtx_\sig]/(I_{\val(v_g),a}+m^k))$, where $m\subset K[\wtx_\sig]$ denotes the maximal ideal. It follows that the set $W$ is finite. On the other hand since $I_{\val(v_g),a}$ are generated by monomials they are uniquely determined by the ideals $\gr(I_{\val(v_g),a})$ in the graded ring
\[\gr(O_{\wtx_\sig})=O_{\wtx_\sig}/{m}\oplus m/m^2\oplus\ldots\]
The connected group $d(G_\sig)$ acts algebraically on $\gr(O_{\wtx_\sig})$ and on the connected component of the Hilbert scheme with fixed Hilbert polynomial. In particular it acts trivially on its finite subset $W$ and consequently $d(G_\sig)$ preserves $\gr(I_{\val(v_g),a})$ and $G_\sig$ preserves $I_{\val(v_g), a}$.
\end{proof}

Let $R\subset K$ be a ring contained in the field. We can order valuations by writing $$\mu_1>\mu_2\quad{\rm if}\quad\underset{a\in R}{\forall} \quad \mu_1(a)\geq \mu_2(a)\quad{\rm and }\quad\mu_1\neq \mu_2.$$ A cone $\sig$ defines a partial ordering: $\quad \quad v_1>v_2 \quad {\rm if}\quad v_1-v_2\in\sig\setminus \{0\}.$ 

\noindent Both orders coincide for $K[X_\sig]\subset K(X_\sig)$:  $v_1>v_2$ iff $\val(v_1)>\val(v_2)$.

\begin{lemma}\label{equiv2}
Let $\sig$ be a cone in $N^{\QQ}_\sig$ with the lattice of 1-parameter subgroups $N_\sig\subset N_\sig^Q$ and the dual lattice of characters $M_\sig$. Let $\mu$ be any integral ({or} rational) valuation centered on $\cl(\Otau)$, where $\tau\preceq\sig$. Then the restriction of $\mu$ to $M_\sig\subset M_\sig\times M_{\reg(\sig)}\subset K(\wtx_\sig)^*$ defines a functional on $\tau^\vee\subseteq M_\sig^Q$ corresponding to a vector $v_\mu\in\inte\tau$ such that 
$F(v_\mu)=\mu(x^F)$ for $F\in M_\sig$ and
$\mu\geq\val(v_\mu)$ on $\wtx_\sig$.
\end{lemma}

\begin{proof} $I_{\mu,a}\supseteq (x^F\mid \mu(x^F)\geq a)=(x^F \mid F(v_\mu)\geq a)=I_{\val(v_\mu),a}$.
\end{proof}
\begin{lemma} \label{group} Let $\Ga\subset \Ga_\sig$ be a finite group acting on $\wtx_\sig$. Let $\pi:N^{\QQ}\to (N^\Ga)^\QQ$ denote the projection corresponding to the geometric quotient $\wtx_\sig\to\wtx_{\pi(\sig)}=\wtx_{\sig}/\Ga$.  Then $\val(v)$ is $G_\sig$-invariant on $\wtx_\sig$ iff $\val(\pi(v))$ is $G_\sig$-invariant on $\wtx_{\pi(\sig)}.$
\end{lemma}

\begin{proof} ($\Rightarrow$) $\val(v)$ is $G_\sig$-invariant on $K[\wtx_\sig]$
and it is invariant on $K[\wtx_\sig]^\Ga$.

\noindent
($\Leftarrow$) Note that $\pi$ defines an inclusion of same dimension lattices $N\hookrightarrow N^\Ga$ and $M^\Ga\hookrightarrow M$.

Assume that $\val(\pi(v))$ is $G_\sig$-invariant. It defines a functional on the lattice $M^\Ga$ and its unique extension to $M\supset M^\Ga$ corresponding to $\val(v)$. Since $g_*(\val(\pi(v)))=\val(\pi(v))$, we have $g_*(\val(v))_{|M^\Ga}=\val(v)_{|M^\Ga}$ and consequently
$g_*(\val(v))_{|M}=\val(v)_{|M}. $
By Lemma \ref{equiv2}, $g_*(\val(v))\geq \val(v)$ for all $g\in G_\sig$. Thus $\val(v)\geq g_*^{-1}(\val(v))$ for all $g^{-1}\in G_\sig$. Finally $g_*(\val(v))= \val(v)$.
\end{proof}

\subsection{Stability of centers from $\para(\pi(\tau))$}

In the following let   $\Del^\sig$ be a decomposition of $\sig\in\Sig$  such that $\wtx_{\Del^\sig}\to \wtx_\sig$ is $G_\sig$-equivariant, $\tau\in \Del^\sig$ be its face and 
$\Ga$ be a finite subgroup of $\Ga_\sig$. Denote by $\pi :(\sig, N_\sig)\to (\sig^\Ga, N^\Ga_\sig)$ the linear isomorphism and the lattice inclusion corresponding to the quotient $X_\sig\to X_{\sig}/\Ga=X_{\pi(\sig)} $.

\begin{lemma}\label{sp} 
 Assume that for any $g\in G_\sig$, there exists a cone $\tau_g\in\Del^\sig$ such that  $g\cdot(\cl(\Otau))=\cl(\Otau_g)$.
%
%
 Let $v\in\inte(\pi(\tau))\subset N^\Ga_\sig\subset N_\sig$ be an integral vector such that $\val(v)$ is not $G_\sig$-invariant on $\wtx_{\sig}{/\Ga}$. Then there exist  integral vectors $v_1\in\inte(\pi(\tau))$ and $v_2\in\pi(\tau)$ such that $$v=v_1+v_2.$$ Moreover if there exists $v_0\in \pi(\sig)$ (not necessarily integral) such that $\rm{val}(v_0)$ is $G_\sig$-invariant and $v>v_0$ on $\pi(\sig)$ then $v_1>v_0$ on $\pi(\sig)$.
 
\end{lemma}
\begin{proof}
If $\rm{val}(v)$ is not $G_\sig$-invariant on $\wtx_\sig/\Ga$ then by Lemma \ref{le:semiinv}
there exists an element $g\in G_\sig$ such that $\mu_g=g_*(\val(v))$ is not a toric valuation. By the assumption $\mu_g$ is centered on $\cl(\widetilde{O}_{\pi(\tau_g)}$. Then by Lemma \ref{equiv2} it defines $v_g\in \inte\pi(\tau_g)$ such that $\mu_g(x^F)=F(v_g)$ for $F\in\sig^\vee$. Moreover $\mu_g>\val(v_g)$. Then
the valuation $g_*^{-1}(\val(v_g))$ is centered on $\cl(\widetilde{O}_{\pi(\tau)}$
 Thus it defines an integral $v_1\in \rm{int}(\pi(\tau))$ such that $v> v_1$ on $\pi(\tau)$ and $v_2:=v-v_1$. 
Then
\[\val(v)=g_*^{-1}(\mu_g)>g_*^{-1}(\val(v_g))\geq\val(v_1). \]

\noindent
Note also that if $v\geq v_0$ then $\mu_g=g_*(\val(v))\geq \val(v_0)$ and $\val(v_g)\geq \val(v_0)$. Thus also $\val(v_1)\geq \val(v_0)$.
\end{proof}
\begin{lemma}\label{equiv}  All valuations $\val(v)$, where $v\in\vr$, $\vr\in \Ver (\Del^\sig)\setminus \Ver (\sig)$, are $G_\sig$-invariant. 
\end{lemma}

\begin{proof} Let $v_\vr$ be the primitive generator of $\vr\in \Ver (\Del^\sig)\setminus \Ver (\sig)$. The ray $\vr$ corresponds to an exceptional divisor $D_\vr$. By the definition there is no decomposition $v_\vr=v_1+v$. Thus by the previous lemma (for $\Ga=\{e\}$), $\val(v)$ is $G_\sig$-invariant.
\end{proof}
\begin{lemma} For any $\tau\subset\sig$, the closure of the orbit $\cl{(\Otau)}\subset \wtx_{\sig}$ is $G_\sig$-invariant.
\end{lemma}
\begin{proof} By Lemma \ref{inclusion}, the ideal of $\cl{(\Otau)}\subset \wtx_\sig$ is generated by all functions with nontrivial $\Ga_\sig$-weights.
\end{proof}
\begin{lemma} \label{stab} The valuations $\val(v)$,  where $v\in \para(\pi(\tau))$, are $G_\sig$-invariant on $\wtx_{\Del^\sig}$. Moreover $v\in\inte(\pi(\sig_0))$, for some $\sig_0\leq\sig$.
\end{lemma}

\begin{proof} Let $v\in\para(\pi(\tau))$, where $\pi(\tau)\in\pi({\Del^\sig})$ is a minimal integral vector such that $\val(v)$ is not $G_\sig$-invariant. We may assume that $v\in\inte(\pi(\tau))$ passing to its face if necessary. Let $\sig'\in\overline{\sig_0}$ be a face of $\sig$ such that $v\in\inte\pi(\sig')$. In particular $\pi(\sig') \supset \pi(\tau) $.
 Then $\pi(\Del^\sig)_{|\pi(\sig')}=\pi(\Del^\sig)_{|\pi(\sig_0)}\oplus\la e_1,\dots,e_k\ra $ by Lemmas \ref{pr} and \ref{sub} and 
 $ v \in {\rm par}{(\pi(\tau))}\subset \pi(\sig_0) $.
 Thus $\sig'=\sig_0$ and $v\in\inte (\pi(\sig_0))$.
Let  $$\pi(\tau)=\la v_1,\dots,v_k,w_1,\dots,w_\ell\ra,$$\noindent where $v_1,\dots,v_k\in\Ver (\pi(\tau))$ and $w_1,\dots,w_\ell\in\Ver (\pi({\Del^\sig}))\setminus \Ver (\pi(\sig))$. By Lemma \ref{equiv}, $\val(w_1),\dots,\val(w_\ell)$ are $G_\sig$-invariant. Write $$v=\alpha_1v_1+\cdots+\alpha_kv_k+\alpha_{k+1}w_1+\cdots+\alpha_{k+\ell}w_\ell,$$\noindent  where $0<\alpha_i<1$. Note that $$v\geq v_0=\alpha_{k+1}w_1+\cdots+\alpha_{k+\ell}w_\ell$$ and $\cl(\widetilde{O}_{\pi(\sig_0)})\subset\wtx_{\pi(\sig)}$ is $G_\sig$-invariant. By Lemma \ref{sp} for $v\in\pi(\sig_0)\leq\pi(\sig)$ and $v>v_0$ we can find integral vectors $v', v''\in\pi(\sig)$ such that $v=v'+v''$, $v'\geq v_0$. Then $$v'':=v-v'\leq v-v_0=\alpha_1v_1+\cdots+\alpha_kv_k.$$ Thus $v''\in \para\la v_1,\ldots,v_k\ra\subseteq \rm{par}(\pi)(\tau)$. 
Write $v'' :=\beta_1v_1+\cdots+\beta_kv_k$, where $\beta_i\leq\alpha_i$. Then $$v'=v-v''=(\alpha_1-\beta_1){v_1}+\cdots+(\alpha_k-\beta_k){v_k}+\alpha_{k+1}w_1+\cdots+\alpha_{k+\ell}\in\para(\pi(\tau)).$$ By the minimality  assumption, $\val(v')$ and $\val(v'')$ are $G_\sig$-invariant and by Lemma \ref{convex}, $\val(v)=\val(v'+v'')$ is $G_\sig$-invariant. 
\end{proof}

\begin{corollary} \label{stab2} Let $\Del=\{\Del^\sig\in\Sig\}$ be a decomposition of $\Sig$. Let $\tau\in\Del^{\sig}$ be an independent face. Then the vectors in $(\pi_{\sig_{|\tau}})^{-1}(\para(\pi_\sig(\tau)))$ are stable.
\end{corollary}

\begin{proof} Put $\Ga=\Ga_\tau$.  Let $\pi: (\sig,N_\sig)\to (\sig, N^\Ga_\sig)$ be the linear isomorphism and a lattice inclusion corresponding to the quotient $X_\sig\to X_\sig/\Ga$.
   Then by Lemma \ref{pro2}, $\pi(\tau)\simeq\pi_\tau(\tau)\simeq\pi_\sig(\tau)$ and by Lemma \ref{stab} vectors in $(\pi_{\sig_{|\tau}})^{-1}(\para(\pi_\sig(\tau)))=\pi^{-1}(\para(\pi(\tau)))$ are stable. 
\end{proof}

\begin{corollary} \label{para}\begin{enumerate} 
\item Assume that for any $g\in G_\sig$, there exists $\tau_g\in \Del^\sig$ such that $g(\cl(\Otau))=\cl(\Otau_g)$. Then $\cl(\Otau)$ is $G_\sig$-invariant.
Moreover  all valuations $\val(v)$, where $v\in \overline{\rm{par}}\,(\tau)\cap \inte(\tau)$, are $G_\sig$-invariant.

\item Let $\tau\in \Del^\sig$ be an independent cone such that $\cl(\Otau)$
%
%
 is $G_\sig$-invariant.  Then for any $v\in\pi^{-1}_\sig(\overline{\rm{par}}\,(\pi(\tau))\cap \inte(\pi(\tau)))$ the valuation $\val(v)$ is $G_\sig$-invariant.
\end{enumerate}
\end{corollary}

\begin{proof}
1. Let $\tau=\la v_1,\dots,v_k\ra$ and  $v=\alpha_1v_1+\cdots+\alpha_kv_k$, where $0<\alpha_i\leq 1$, be a minimal vector in $\inte(\tau)\cap\overline{\rm{par}}\,(\tau)$ such that $\val(v)$ is not $G_\sig$-invariant. Then by Lemma \ref{sp}, the vector $v$ can be written as $v=v'+v''$, where $v',v''<v$, $v'\in \inte(\tau)$, $v''\in\tau$. Thus $v'=\alpha'_1v_1+\cdots+\alpha'_kv_k$ where $0<\alpha'_i\leq\alpha_i\leq 1$
%
%
 and $v''=\alpha''_1v_1+\cdots+\alpha''_kv_k$, where $0\leq\alpha''_i=\alpha_i-\alpha'_i< 1$. Then $v'\in \inte(\tau)\cap\overline{\rm{par}}\,(\tau)$ and $v''\in \para(\tau)$. By Corollary \ref{stab2}, $\val(v'')$ is $G_\sig$-invariant on $\wtx_\sig$. By the minimality assumption $\val(v')$ is $G_\sig$-invariant. Since $v=v'+v''$, the valuation $\val(v)$ is $G_\sig$-invariant on  $\wtx_\sig$ and its center $Z(\val(v))$ equals $\\cl(\widetilde{O}_\tau)$.

2. Let $\pi:N\to N^\Ga$ be the projection corresponding to the quotient $X_\sig\to X_{\sig}/\Ga_\tau$. Then by Lemma \ref{pro2}, we have $\pi(\tau)\simeq\pi_\sig(\tau)$. The proof is now exactly the same as the proof in 1 except that we replace $\wtx_{\Del^\sig}$ with $\wtx_{\Del^\sig}/{\Ga_\tau}$.
\end{proof}
\subsection{Fixed points of the action}
We shall carry over  the concept of fixed point set of the action of $\ks$ to the scheme $\wtx_{\Del^\sig}$. The problem is that $\wtx_{\Del^\sig}$ does not contain enough closed points.

\begin{definition} A point $p\in\wtx_{\Del^\sig}$ is a {\it fixed point} of the action of $\ks$ if $\ks\cdot p=p$ and 
  $\ks$ acts trivially on the residue field $K_p$ of $p$
\end{definition}
\noindent

\begin{lemma} \label{fixed} The set of all fixed points ${\wtx_{\Del^\sig}}^{\ks}$ of the action of $\ks$ is given by the union of the closures of the  orbits $\cl(\widetilde{O}_\del)$ defined by circuits $\del\in\Sig$. The $\cl(\widetilde{O}_\del)$ are maximal irreducible components of the fixed point set.
\end{lemma}
\begin{proof} A point $p$ of $\wtx_{\Del^\sig}$ is a point of a locally closed subscheme defined by the orbit $\widetilde{O}_\tau$.
If $\tau$ is independent then there exists an invertible  character $x^F$, where $F\in\tau^\perp$, on which $\ks$ acts nontrivially.
Then the action on $K_p\ni x^F$ is nontrivial. If $\tau$ is dependent then
the action on $K[\widetilde{O}_\tau]$ and on $K_p$ is trivial so $p\in \widetilde{O}_\tau$ is a fixed point and $p\in \cl({\widetilde{O}_\delta})$, where $\delta\preceq\tau$ is a circuit.

\end{proof}

\begin{corollary} \label{g} Let $\delta\in \Del^\sig$ be a circuit. Then $\cl(\widetilde{O}_\delta)$ is $G_\sig$-invariant.
\end{corollary}

\begin{proof} By Corollary \ref{fixed2}, 
$\cl(\widetilde{O}_\delta)$ is an irreducible component of a $G_\sig$-invariant closed subscheme $\wtx_{\Del^\sig}^{\ks}$. Thus by the Corollary \ref{para}(1) it is $G_\sig$-invariant.
\end{proof}

\subsection{Stability of $\Ctr_+(\sigma)$}

In the sequel $\del=\la v_1,\ldots,v_k\ra\in\Del^\sig$ is a circuit. 
Let $\Ga\subset\Ga_\sig=K^*$ be a finite group.
Denote by $\pi$ (resp. $\pi_\Ga$) the projection corresponding to the quotient $X_\del\to X_\del//K^*$ (resp. $X_\del\to X_\del/\Ga$).
In particular $\pi_\sig(\del)\simeq\pi(\del)$.
Write $\pi_\sig(\del)=\la w_1,\ldots,w_k\ra$ and let  $\sum_{r'_i>0}r'_iw_i =0$ be the unique relation between vectors (**) as in Section \ref{dep}. Set
${\rm Ctr_+}(\del)=\sum_{r'_i>0}w_i \in \overline{\rm{par}}(\pi_\sig(\del_+))\cap\inte(\pi_\sig(\del_+))$, where 
${\del_+}=\la v_i\mid r_i>0\ra$.

Denote by $\whx_\del$ the completion of $\wtx_{\Del^\sig}$ at $\widetilde{O}_\del$. By Corollary \ref{g},  the generic point $\widetilde{O}_\del\in\wtx_{\Del^\sig}$  is $G_\sig$-invariant and thus  $G_\sig$ acts on $\whx_\del$. Moreover  
$K[\whx_\del]=K(\widetilde{O}_\del)[[\del^\vee]]$ is  is faithfully flat over a $\cO_{\wtx_{\Del^\sig,\widetilde{O}_\del}}$.
Also, $\widehat{\cO}_{X_{\pi_\Ga(\Del)},\widetilde{O}_{\pi_\Ga(\del)}}=K(\widetilde{O}_{\pi_\Ga(\del)})[[\underline{\pi_\Ga(\del)}^\vee]]$ is faithfully flat over ${\cO}_{X_{\pi_\Ga(\Del)},\widetilde{O}_{\pi_\Ga(\del)}}$.

 The valuation $\val(v)$ on the local ring $\cO_{X_{\Del^\sig},\widetilde{O}_\del}$ (or ${\cO}_{X_{\pi_\Ga(\Del)},\widetilde{O}_{\pi_\Ga(\del)}}$), where $v\in\del$ extends to its completion $\widehat{\cO}_{X_\Del,\widetilde{O}_\del}=K(\widetilde{O}_\del)[[\underline{\del}^\vee]]$ (respectively$K(\widetilde{O}_{\pi_\Ga(\del)})[[\underline{\pi_\Ga(\del)}^\vee]]$). Moreover $\rm{val}\,(v)_{| K(\widetilde{O}_\del)^*}=0$ and the action of $K^*$ on $K(\widetilde{O}_\del)$ is trivial. 
As in Lemma \ref{closed} we get
\begin{lemma} \label{compa} The valuation $\val(v)$, where $v\in\pi_\Ga(\del)$,  is $G_\sig$-invariant on $\whx_\del/\Ga$ iff it is $G_\sig$-invariant on $\wtx_{\Del^\sig}/\Ga$.
\end{lemma}
\begin{lemma} \label{inva}\begin{enumerate}
\item $\cl( \widetilde{O}_{\del_-}), \cl(\widetilde{O}_{\del_+})\subset\whx_\del$ are $G_\del$-invariant.

\item $\cl(\widetilde{O}_{\del_-}),\cl(\widetilde{O}_{\del_+})\subset\widetilde{\Del^\sig}$ are $G_\del$-invariant.
\end{enumerate}
\end{lemma}
\begin{proof}(1) By Lemmas \ref{sigma} and \ref{fix}, the ideal $I_{\cl{(\widetilde{O}_{\del_+})}}\subset K[\whx_\sig]$ of $\cl(\widetilde{O}_{\del_+})=(\widetilde{O}_\del)^+$ is generated by functions with positive weights. (2) Consider morphisms $\whx_\del\buildrel e\over\to \wtx_{\Del^\sig}\to X_{\Del^\sig}$. The morphism $e$ is $G_\sig$-equivariant and maps the generic points of the orbits $\widetilde{O}_\del^\pm$ on $\whx_\del$ onto the generic points of the corresponding orbits on $ \wtx_{\Del^\sig}$.
\end{proof}

\begin{corollary}\label{quo} There are open $K^*$-equivariant embeddings of schemes $(\whx_\del)_-:=\whx_\del\times_{X_\del}(X_\del)_-\subset \whx_\del$ and $ 
(\whx_\del)_+:=\whx_\del\times_{X_\del}(X_\del)_+\subset \whx_\del$.  
\end{corollary}

\begin{lemma} There exist quotients  
 $$(\whx_\del)_{-}/\ks=\whx_\del/\ks\times_{X_{\del}{/\ks}}(X_\del)_{-}/\ks, \quad   (\whx_\del)_+/\ks=\whx_\del/\ks\times_{X_{\del}{/\ks}}(X_\del)_+/\ks.$$  
\end{lemma}
\begin{proof} The proof id identical with the proof of Lemma \ref{le:q} except for we need to use Lemma \ref{le:100} below instead of Lemma \ref{le:10}
\end{proof}
\begin{lemma} \label{le:100} Let $K(\widetilde{O}_{\del})[\underline{\del}^\vee]=\underset{a\in \ZZ}{\oplus}K(\widetilde{O}_{\del})[\underline{\del}^\vee]^a$ be
a decomposition according to weights. Then $K(\widetilde{O}_{\del})[\underline{\del}^\vee]^a$ is generated over $K(\widetilde{O}_{\del})[\underline{\del}^\vee]^0$ 
by
finitely many monomials.
\end{lemma}
\begin{proof} Let $x^{F_1},\dots,x^{F_k}$
generate $K(\widetilde{O}_{\del})[\underline{\del}^\vee]$. Set
$b:=\max\{|F_1(v_\del)|,\dots,|F_k(v_\del)|\}$. We show that  all the elements
$x^{\alpha_1F_1+\cdots+\alpha_kF_k}$, where $$\alpha_1F_1(v_\del)+\cdots
+\alpha_kF_k(v_\del)=a\quad  \mbox{and}\quad  0\leq\alpha_i\leq k\cdot b^2+|a|,$$ generate
$K[\widetilde{\del}^\vee]^a$ over $K[\widetilde{\del}^\vee]^0$. Without loss of generality assume that $\alpha_i>k\cdot b^2+|a|$ and 
$F_i(v_\del)>0$. Then 
$${}k\cdot b\cdot\max\{\alpha_i\mid F_i(v_\del)<0\}\geq
-\underset{F_i(v_\del)<0}{\sum}\alpha_iF_i(v_\del)
=\underset{F(v_\del)>0}{\sum}\alpha_iF(v_\del)-a\geq kb^2+|a|-a\geq kb^2.$$
Thus there exists $j$ such that $\alpha_j\geq\ds{kb^2\over kb}=b$ and
$F_j(v_\del)<0$. But then
$$x^{\alpha_1F_1+\cdots+\alpha_kF_k}=x^{F_i(v_\del)F_j-F_j(v_\del)F_i}\cdot
x^{\alpha_1F_1+\cdots+\left(\alpha_i+F_j(v_\del)\right)F_i+\cdots+\left(\alpha_j+F_i(v_\del)\right)F_j+\cdots +\alpha_kF_k},$$

\noindent where $x^{F_i(v_\del)F_j-F_j(v_\del)F_i}\in K[\widetilde{\del}^\vee]^0$ and $x^{\alpha_1F_1+\cdots+\left(\alpha_i+F_j(v_\del)\right)F_i+\cdots+\left(\alpha_j-F_i(v_\del)\right)F_j+\cdots +\alpha_kF_k}\in K[\widetilde{\del}^\vee]^a$ with smaller exponents.
\end{proof}

\begin{lemma}\label{3}
 The action of $G_\sig$ on $(\whx_\del)_-$ and $(\whx_\del)_+$ descends to $(\whx_\sig)_-/K^*$ and $(\whx_\sig)_+/K^*$.
\end{lemma}

\begin{proof} 
The proof is almost identical with the proof of Lemma \ref{2} except for the Lemma \ref{21} which shall be replaced with Lemma \ref{31}. 
We replace open affine subsets $V$ with open affine subsets $V$ satisfying the condition (***) below. 
\end{proof}

\begin{lemma}\label{31} Let $V$ be an open affine $\Ga$--invariant subscheme of $\whx_\tau=X_\tau\times_{X_\del}\whx_\del\subset (\whx_\sig)_-$. Then $V$ satisfies the condition.

\noindent
(***) For any open affine $\Ga$--invarient subscheme $U\subset V$ there is an inclusion of open affine subschemes $U/\Ga\subset V/\Ga$.
\end{lemma}

\begin{proof} Let $Z\subset\whx_\del \setminus U$ be a closed subscheme. Then the ideal $I_Z\subset K[\whx_\del]$ is $\ks$--invariant. Let $f=\sum_{F\in\del^\vee}\alpha_Fx^F\in I_Z$ and $f_a:=\sum_{F\in(\del^\vee)^a}\alpha_Fx^F$ be a part of $f$ with weight $a$. Let $m_\del\subset K[\whx_\del]$ be the maximal ideal. Then for any $m^k_\del$, the decomposition of  $[f]:=f+m^k_\del=\sum f_a+m^k_\del=\sum [f_a]\in K[\whx_\del]/m^k_\del$ is  finite. Moreover $t[f]=\sum t^a [f_a]$. It follows that $[f_a]=f_a+m^k_\del\in I_Z+m^k_\del$ and $f_a\in I_Z$. Thus $I_Z$ is generated by semininvariant generators $f_1,\dots,f_k$ with weight $a_1,\dots,a_k$. Note that since $\tau$ is independent, we have that $v_\sig\not\in\rm{span}\,(\tau)$, and $v_\sig$ is not orthogonal to $(\rm{span}\,(\tau))^\perp=\tau^\perp\otimes_\ZZ\QQ$. Thus there exists $F\in \tau^\perp$ such that $F(v_\sig)=a\neq 0$. The corresponding character $x^F$ is invertible on $\whx_\tau$. The function $g_i=f^a_i (x^F)^{-a_i}$ are invariant (with weight $aa_i-aa_i=0$) and 
$(\whx_\tau)_{g_i}=(\whx_\tau)_{f_i}=V_{f_i}$.
Then $U=\underset{i=1}{\bigcup}V_{g_i}=\underset{i=1}{\bigcup}(\whx_\tau)_{g_i}$
and $U_{g_i}/\ks=V_{g_i}/\ks\subset U/\ks, V/\ks$. It follows that $U/\ks\subset V/\ks$.
\end{proof}
Proposition \ref{sigma2}, Lemma \ref{sigma} and the above imply:
\begin{corollary}\label{bira}
  The morphisms $\widehat{\phi}_-: (\whx_\del)_-/\ks\to \whx_\del/\ks$ and
 $\widehat{\phi}_+: (\whx_\del)_+/\ks\to \whx_\del/\ks$
are $G_\sig$-equivariant, proper and birational.
\end{corollary}

\begin{lemma} \label{stable} The vector  $v:=\rm{Mid}\,({\rm Ctr_+}(\del),\del)=\pi_{\sig_{|\del^-}}^{-1}({\rm Ctr_+}(\del))+\pi_{\sig_{|\del^+}}^{-1}({\rm Ctr_+}(\del))$
 is stable.
\end{lemma}

\begin{proof}  Set $ v_-:=\pi_{\sig_{|\del^-}}^{-1}({\rm Ctr_+}(\del)) $   and 
$ v_+:=\pi_{\sig_{|\del^+}}^{-1}({\rm Ctr_+}(\del))$.
 By Lemma \ref{inva}, $\cl(\widetilde{O}_{\del_+})\subset \wtx_{\Del^\sig}$ is $G_\sig$-invariant and, by Corollary  \ref{para}(2) and   Lemma \ref{compa}, $\val(v_+)$ is $G_\sig$-invariant on $\wtx_{\sig}$ and on $\whx_{\del}$. Hence the valuation $\val(v_+)$ descends to a $G_\sig$-invariant valuation $\val(\pi(v_+))$ on $\whx_{\del}{//\ks}=K(\widetilde{O}_\del)[[\underline{\del}^\vee]]^{K^*}$. 
By Corollary \ref{bira},  $\val(\pi(v_-))=\val(\pi(v_+))$
%
%
 is
$G_\sig$-invariant  on 
 $(\whx_{\del})_+/\ks=\whx_{\partial_-(\del)}/K^*=\whx_{\pi(\partial_-(\del))}$. 
 Let $\Gamma\subset K^*$ be the subgroup generated by  all subgroups $\Gamma_\tau\subset K^*$,  where $\tau\in\partial_-(\del)$. 
 Then $K^*/\Ga$ acts freely on $X_{\partial_-(\del)}/\Ga=(X_\del)_+/\Ga$. Let $j: (X_\del)_+/\Ga\to (X_\del)_+/K^*$ be the natural morphism. 
 Let $\pi_\Ga:\del\to \pi_\Ga(\del)$ be the projection corresponding to the quotient $X_\del \to X_\del/\Ga$. 
 By Lemma \ref{le: p}, 
for any $\tau\in\del_-$, the restriction of $j$ to $X_\tau/\Ga\subset (X_\del)_+/\Ga$ is given by $$j:X_\tau/\Ga=
X_{\underline\tau}/\Ga\times O_\tau/\Ga\to X_\tau/K^*=X_{\underline\tau}/\Ga\times O_\tau/K^*.$$ Thus  ${j}^*(\cI_{\val(\pi_\Ga(v_-)),a})=\cI_{\val(\pi(v_-)),a}$.
%
%
 Consider the natural morphisms $i_\Ga:(\whx_\del)_+/\Ga\to(X_\del)_+/\Ga$ and $i_{K^*}:(\whx_\del)_+/\ks\to(X_\del)_+/\ks$. Then
 ${i_\Ga}^*(\cI_{\val(\pi_\Ga(v_-)),a,X_\del/\Ga})=\cI_{\val(\pi_\Ga(v_-)),a,\whx_\del/\Ga}$ and $({i_{\ks}})^*(\cI_{\val(\pi(v_-)),a,X_\del/\ks})=\cI_{\val(\pi(v_-)),a,\whx_\del/\ks}$. Let $\hat{j}:(\whx_\del)_+/\Ga\to (\whx_\del)_+/K^*$ be the natural morphism induced by $j$. The following diagram commutes.
\[\begin{array}{rcccccccc}

&&& (\whx_\del)_+/\Ga
 & \buildrel\widehat{j}\over\rightarrow &(\whx_\del)_+{/\ks}
&&& \\
&&&\downarrow {i_\Ga} & & \downarrow i_{\ks}  &&&\\
&& & (X_\del)_+/\Ga  &  \rightarrow & (\whx_\del)_+{/\ks}.&&&\\

\end{array}\] 
 Thus we get $\hat{j}^*(\cI_{\val(\pi_\Ga(v_-)),a,\whx_\del/K^*})=\cI_{\val(\pi(v_-)),a,\whx_\del/\Ga}$.
Since  the morphism $\widehat j$ is $G_\sig$-equivariant it follows that $\val(\pi_\Ga(v_-))$ is $G_\sig$-equivariant on $(\whx_\del)_+/\Ga$. Since $(\whx_\del)_+\subset \whx_\del$ is an open $G_\sig$-equivariant inclusion and $\Ga$ is finite  we get that the morphism  $(\whx_\del)_+/\Ga\subset (\whx_\del)/\Ga$
is an open $G_\sig$-equivariant inclusion. 
Thus the valuation $\val(\pi(v_-))$ is $G_\sig$-equivariant on $\whx_\del{/\Ga}$ and on $\wtx_{\Del^\sig}/\Ga$ (Lemma \ref{compa}). Finally, by Lemma \ref{group}, $\val(v_-)$ it is $G_\sig$-equivariant on $\wtx_{\Del^\sig}$. By the convexity, $\val(v)=\val(v_++v_-)$ is $G_\sig$-equivariant on $\wtx_{\Del^\sig}$.
\end{proof}
\subsection{Canonical coordinates on $\Sig$}

\medskip\noindent
Note that for any $\sig\in\Sig$ we can order the coordinates according the weights
\[X_\sig\simeq \bbA^k=\bbA_{a_1}\oplus \bbA_{a_2}\oplus\cdots\oplus \bbA_{a_\ell},\]
where $a_1< a_2<\ldots< a_\ell$ and $\Ga_\sig$ acts on $\bbA_{a_i}$ with character $t\to t^{a_i}$, where $t\in\Ga$ and $a_i\in\bZ_n$ if $\Ga_\sig\simeq \bZ_n$ or $a_i\in\bZ$ if $\Ga_\sig\simeq\ks$.( In the first case $a_i$ are represented by integers from $[0,n-1]$. Let us call these coordinates \textit{canonical}.) The canonical coordinates are preserved by the group $\rm{Aut}\, (\sig)$ of all automorphisms of $\sig$ defining $\ks$-equivariant automorphisms of $X_\sig$. Since all stable vectors $v\in\sig$ define $G_\sig$-invariant valuations $\rm{val}\,(v)$ on $\wtx_\sig$, they are in particular $\rm{Ant}\,(\sig)$-invariant. Thus all stable vectors $v\in\sig$ can be assigned the canonical coordinates in a unique way.

\subsection{Canonical $\pi$-desingularization of cones $\sig$ in $\Sig$} 
Given the canonical coordinates we are in position to construct a canonical $\pi$-desingularization of $\sig$ or its subdivision $\Del^\sig$. We eliminate all choices of centers of star subdivisions in the
$\pi$-desingularization algorithm by choosing the center with the smallest canonical coordinates (ordered lexicographically).
\subsection{Canonical $\pi$-desingularization of  $\Sig$}
\medskip
\noindent
Note that the $\pi$-desingularization of an independent $\tau\in\Sig$ is nothing but a desingularization of $\pi(\tau)$.

\noindent
For any cone $\sig=\la v_1,\ldots,v_k\ra \in\Sig$ the vector $v_\sig:=v_1+\ldots+v_k\in\overline{\para}(\sig)$ is stable (Lemma \ref{para}).
Order all cones $\sig\in\Sig$ by their dimension and apply star subdivision at $\la v_\sig\ra\in\sig$ starting from the heighest dimensions to the lowest. 
Note that the result of this  subdivision does not depend on the order of cones of the same dimension. That's because none of two cones of 
$\Sig$ which are of the same dimension are the  faces of the same cone. The cones of higher dimensions were already subdivided and all their proper faces occur in the different cones.
Let $\Del=\{\Del^\sig\mid \sig\in\Sig\}$ denote the resulting subdivision. 
Now we apply the canonical subdivision $\Del^\pi_\sig$ to the subdivided  cones $\Del^\sig$ starting from the lowest dimension to the heighest.
The subdivisions $\Del^\pi_*$ applied to any two (subdivided) cones  of the same dimension in $\Sig$ commute since their faces (of lower dimension) are already $\pi$-nonsingular and thus not affected by further subdivisions. Also as before none of two cones 
which are in different faces of $\Sig$ of the same dimension are the  faces of the same cone. Note also that $\Del^\pi_\sig$ depends only on $\sig$ and is independent of the other faces of $\Sig$.

\subsection{Proof of the Weak Factorization Theorem}
The canonical $\pi$-desingularization $\Delta^\pi$ of $\Sigma$ is obtained by a sequence of star subdivisions at stable centers (Lemmas \ref{stable}, \ref{stab2}). By Propositions \ref{blowup} and  \ref{correspondence}, $\Del$ defines a birational projective modification $f:{B}^\pi\to B$. The modification does not affect points with trivial stabilizers  $B_-=X^-\setminus X$ and $Z^+\setminus Z$ (see Proposition \ref{construction}). This means that $(B^\pi)_-=B_-$ and $(B^\pi)_+=B_+$  and $B^\pi$ is a cobordism between $X$ and $Z$.  Moreover $B^\pi$ admits a projective compactification $\overline{B^\pi}=B \cup X \cup Z$. The   cobordism $\overline{B^\pi}\subset \overline{B^\pi}$ admits a  decomposition into elementary cobordisms $B^\pi_a$, defined by the strictly increasing function
$\chi_B$.  Let $F\in \cC((B_a^\pi)^{K^*})$ be a fixed point component
and $x\in F$ be a point.
By Proposition \ref{correspondence} the modification $f:{B^\pi}\to B$
is locally described 
for a toric chart $\phi_\sig:U\to X_\sig$ by  a smooth $\Ga_\sig $-equivariant morphism $\phi_{\Del^\sig}: f^{-1}(U)\to X_{\Del^\sig}$. Then by Lemma \ref{fixed3}, $\phi_{\Del^\sig}(x)$ is in ${O}_\del$, where $\del\in\Del^\sig$ is  dependent and $\pi$-nonsingular.
In particular the cone $\sig\in\Sig$ is also dependent and $\Ga_\sig=K^*$. The toric chart $\phi_\sig:U\to X_\sig$ can be extended to an $K^*$-equivariant \'etale morphism  $\psi_\sig:U\to X_\sig\times \bbA^r$, where the action of $\Ga_\sig=K^*$ on $\bbA^r$ is trivial. Moreover
since the toric chart $\phi_\sig$ is compatible with a divisor $D$ we can assume that that the all components of $D$ are descibed by some coordinates on $X_\sig\times \bbA^r\simeq \bbA^n$ (see also Section \ref{stratification}).  The morphism $\phi_{\Del^\sig}$ determines
a $K^*$-equivariant \'etale morphism $\psi_{\Del^\sig}: f^{-1}(U)\to X_{\Del^\sig}\times \bbA^r$.

So we locally have a
$K^*$-equivariant e\'tale morphism 
$\psi_\del: V\to X_\del\times \bbA^r,$ where $V\subset \phi_{\Del^\sig}^{-1}$ is an  affine $K^*$-invariant subset of $B^\pi_a$. Similar to Proposition \ref{local} we get   a diagram 
\[\begin{array}{rccccc}
&{(B^\pi_a)}_-/K^*& \supset &{V_x}_-/K^*& \rightarrow
&{X_\del}_-/K^*\times \bbA^r\\
&\uparrow\psi_-&&\uparrow& &\uparrow \phi_{-}\\
&\Ga({(B^\pi_a)}_-/K^*,{(B^\pi_a)}_+/K^*)    &\supset &\Ga({V_x}_-/K^*,{V_x}_+/K^*)    & \rightarrow  &\Ga({X_\del}_-/K^*,{X_\del}_+/K^*) \times \bbA^r\\
& \downarrow\psi_+&& \downarrow & &\downarrow \phi_{+}\\

&{(B^\pi_a)}_+/K^*&\supset&{V_x}_+/K^*&\rightarrow & {X_\del}_+/K^*\times \bbA^r\\

\end{array}\]
with horizontal arrows \'etale induced by
$${(B^\pi_a)}//K^*\supset{V_x}//K^*\rightarrow  {X_\del}//K^*\times \bbA^r.$$
 Here $\Ga(X_-/K^*,X_+/K^*)$ denotes the normalization of the graph of a birational map $X_-/K^*\dashrightarrow X_+/K^*$ for a relevant cobordism $X$.
 We use  functoriality of the graph (a dominated component of the fiber  product  $X_-/K^*\times_{X//K^*}X_+/K^*)$. 
 By Lemma \ref{local} the morphisms $\phi_-$ and $\phi_+$ are blow-ups at smooth centers.
 Thus $\psi_-$ and $\psi_+$ are locally blow-ups at smooth centers so they are  globally blow-ups at  smooth centers. 
The components of $D^\pi:=B^\pi\setminus (U\times K^*)$ are 
either the strict transforms of components of $D$ on $B$ or the exceptional divisors of $B^\pi\to B$. In either case they correspond to toric divisors on $X_{\Del^\sig}\times \bbA^r$ because of the compatibility of charts. Thus the divisor
$D_{a-}: =(D\cap (B_a)_-)/K^*=((B_a)_-/K^*)\,\setminus\, U$   corresponds to a toric dvisor on a smooth toric variety$(X_\del)_-/K^*\times \bbA^r\simeq \bbA^{n-1}$. The center of the blow-up corresponds to a toric subvariety $O_\del=\{0\}\times \bbA^r\subset (X_\del)_-/K^*\times \bbA^r\simeq \bbA^{n-1}$.
This shows that the centers of blow-ups has SNC with with complements of $U$.

Note that every $K^*$-equivariant automorphism of $B$ preserving the divisor $D=B\setminus U$ transforms isomorphically  the strata, the relevant toric charts and the corresponding cones. This induces anautomorphism of $\Sig$ preserving canonical coordinates on the cones. And it lifts  to the $\pi$-desingularization of $\Sig$ and to the corresponding  cobordism $B^\pi\subset\bar{B^\pi}$ constructed via diagrams (2) as in Proposition \ref{correspondence}. The relatively ample divisor for $\bar{B}^\pi \to X$ is a combination of the divisor $X\times\{\infty\}\subset B^\pi$ and the exceptional divisors of the morphism $\bar{B}^\pi \to X\times \bP^1$. Thus it is functorial i.e. invariant with respect to the liftings of automorphisms of $X$ commuting with $X\dashrightarrow Y$ and defines  a decomposition into open invariant subsets $B_a$ and the induced equivariant factorization.
\subsection{The Weak Factorization over an algebraically nonclosed base field.}
For any proper biarational map $\phi: X\dashrightarrow Y$ over a field  $K$ of characterictic zero consider the induced birational map $\phi_{\overline{K}}:=X^{\overline{K}}:=X\times_{\Spec{K}}\Spec{\overline{K}}\dashrightarrow Y^{\overline{K}}:=Y\times_{\Spec{K}}\Spec{\overline{K}}$ over the algebraic closure $\overline{K}$ of $K$. The weak factorization
of   $\phi_{\overline{K}}$ over $\overline{K}$ is $\Gal(\overline{K}/K)$-equivariant and defines the relevant weak factorization of $\phi$ over $K$.

\end{document}